\documentclass[aop,seceqn,MSNbibl,citesort,dvips]{arximspdf}
\usepackage{mathbh}

\doi{10.1214/10-AOP565}
\volume{39}
\issue{2}
\pubyear{2011}
\firstpage{683}
\lastpage{728}

\makeatletter

\newproclaim{Assumption}{Assumption}

\renewcommand{\mid}{\vert}

\newcommand{\ttQ}{\mathtt{Q}}
\newcommand{\ttE}{\mathtt{E}}
\newcommand{\rombc}{\mathrm{bc}}
\newcommand{\romDir}{\mathrm{Dir}}
\newcommand{\romP}{\mathrm{P}}
\newcommand{\ee}{e}
\newcommand{\ti}{\to\infty}
\newcommand{\romper}{\mathrm{per}}

\newcommand{\Z}{\mathbb{Z}}
\newcommand{\N}{\mathbb{N}}
\newcommand{\PP}{\mathbb{P}}
\newcommand{\E}{\mathbb{E}}
\newcommand{\R}{\mathbb{R}}
\newcommand{\RR}{\mathfrak{R}}

\newtheorem{theorem}{Theorem}[section]
\newtheorem{lemma}[theorem]{Lemma}
\newtheorem{prop}[theorem]{Proposition}
\newtheorem{cor}[theorem]{Corollary}

\newcommand{\dd}{{d}}
\newcommand{\eps}{\varepsilon}
\newcommand{\Leb}{\operatorname{Leb}}
\newcommand{\Sym}{\mathfrak{S}}

\newcommand{\Ccal}{{\mathcal C}}
\newcommand{\Dcal}{{\mathcal D}}
\newcommand{\Gcal}{{\mathcal G}}
\newcommand{\Hcal}{{\mathcal H}}
\newcommand{\Lcal}{{\mathcal L}}
\newcommand{\Ncal}{{\mathcal N}}
\newcommand{\Pcal}{{\mathcal P}}
\newcommand{\Tcal}{{\mathcal T}}

\makeatother

\begin{document}
\begin{frontmatter}

\title{A variational formula for the free energy
of an interacting many-particle system\thanksref{T1}}

\runtitle{Free energy of many-particle systems}

\thankstext{T1}{Supported by the DFG-Forschergruppe 718, Analysis and
Stochastics in Complex Physical Systems.}

\begin{aug}
\author[A]{\fnms{Stefan} \snm{Adams}\ead[label=e1]{S.Adams@warwick.ac.uk}},
\author[B]{\fnms{Andrea} \snm{Collevecchio}\thanksref{t1}\ead[label=e2]{collevec@unive.it}}
and
\author[C]{\fnms{Wolfgang} \snm{K\"onig}\corref{}\ead[label=e3]{wolfgang.koenig@wias-berlin.de}}

\runauthor{S. Adams, A. Collevecchio and W. K\"onig}

\affiliation{University of Warwick, Universit\`a Ca' Foscari and
Technical University Berlin and Weierstrass Institute for Applied
Analysis and Stochastics}

\address[A]{A. Adams\\
Mathematics Institute\\
University of Warwick\\
Zeeman Building\\
Coventry CV4 7AL\\
United Kingdom\\
\printead{e1}}
\address[B]{A. Collevecchio\\
Dipartimento di Matematica Applicata\\
Universit\`a Ca' Foscari\\
San Giobbe, Cannaregio 873\\
30121 Venezia\\
Italy\\
\printead{e2}}
\address[C]{W. K\"onig\\
Technical University Berlin\\
Str. des 17. Juni 136\\
10623 Berlin\\
Germany\\
and\\
Weierstrass Institute\\
for Applied Analysis and Stochastics\\
Mohrenstr. 39\\
10117 Berlin\\
Germany\\
\printead{e3}}
\end{aug}

\thankstext{t1}{Supported by Italian PRIN 2007 Grant 2007TKLTSR.}

\received{\smonth{3} \syear{2010}}
\revised{\smonth{5} \syear{2010}}

%
\begin{abstract}
We consider $N$ bosons in a box in $\R^d$ with volume $N/\rho$ under
the influence of
a mutually repellent pair potential. The particle density $\rho\in
(0,\infty)$ is kept
fixed. Our main result is the identification of the limiting free energy,
$f(\beta,\rho)$, at positive
temperature $1/\beta$, in terms of an explicit
variational formula, for any fixed $\rho$ if $\beta$ is sufficiently small,
and for any fixed $\beta$ if $\rho$ is sufficiently small.

The thermodynamic equilibrium is described by the symmetrized trace of
$\ee^{-\beta\Hcal_N}$, where $\Hcal_N$ denotes the corresponding
Hamilton operator. The well-known Feynman--Kac formula reformulates
this trace in terms of $N$ interacting Brownian bridges.
Due to the symmetrization, the bridges are organized in an
ensemble of cycles of various lengths. The novelty of our
approach is a description in terms of a marked Poisson point
process whose marks are the cycles. This allows for an asymptotic
analysis of the system via a
large-deviations analysis of the stationary empirical field.
The resulting variational formula ranges over random shift-invariant
marked point fields and optimizes the sum of the interaction and the relative
entropy with respect to the reference process.

In our proof of the lower bound for the free energy, we drop
all interaction involving ``infinitely long'' cycles, and their possible
presence is signalled by a loss of mass of the ``finitely
long'' cycles in the variational formula. In the proof of the upper
bound, we only keep the mass on the ``finitely long'' cycles.
We expect that the precise relationship between these two bounds
lies at the heart of Bose--Einstein condensation and intend to
analyze it further in future.
\end{abstract}

%

\end{frontmatter}

\section{Introduction and main results}\label{Intro}

In this paper, we study a probabilistic model for interacting
bosons at positive temperature in the thermodynamic limit with positive
particle density. See Section \ref{sec-BEC} for the physical background.

\subsection{The model}\label{sec-model} The main object is
the following symmetrized sum of Brownian bridge expectations:
%
%
\begin{eqnarray}\label{defpartition}
&&
Z_N^{(\rombc)}(\beta,\Lambda)\nonumber\\
&&\qquad=\frac{1}{N!}\sum_{\sigma\in
\Sym_N}\int_\Lambda\dd x_1\cdots\int_\Lambda\dd x_N\\
&&\qquad\quad\hspace*{39pt}{}\times
\bigotimes_{i=1}^N\mu_{x_i,x_{\sigma(i)}}^{(\rombc,\beta)}\biggl[\exp\biggl\{
-\sum_{1\le i<j\le N}\int_0^\beta
v\bigl(\bigl|B_s^{(i)}-B_s^{(j)}\bigr|\bigr) \,\dd
s\biggr\}\biggr].\hspace*{-28pt}\nonumber
\end{eqnarray}
Here $ \mu_{x,y}^{(\rombc,\beta)} $ is the canonical Brownian
bridge measure with boundary condition $\rombc\in\{ \varnothing,
\romper,\romDir\} $, time horizon $ \beta>0 $ and initial point $
x\in\Lambda$ and terminal point $ y\in\Lambda$, and the sum is on
permutations $ \sigma\in\Sym_N $ of $ 1,\ldots,N$. [We write $\mu
(f)$ for the integral of $f$ with respect to the measure $\mu$.] The
\textit{interaction potential} $ v\dvtx\R\to[0,\infty] $ is
measurable, decays sufficiently fast at infinity and is possibly
infinite close to the origin. Our precise assumptions on $v$ appear
prior to Theorem~\ref{thm-mainres} below. We assume that $\Lambda$ is
a measurable subset of $\R^d$ with finite volume.

The boundary condition $\rombc=\varnothing$ refers to the standard
Brownian bridge, whereas for $ \rombc=\romDir$, the expectation
is on those Brownian bridge paths which stay in $ \Lambda$ over the
time horizon $[0,\beta]$. In the case of periodic boundary condition,
$ \rombc=\romper$, we consider Brownian bridges on the torus
$\Lambda= (\R/L\Z)^d $ with side length $L$.

Our main motivation to study the quantity $Z_N^{(\rombc)}(\beta
,\Lambda)$ is the fact that, for both periodic and Dirichlet boundary
conditions, it is related to the $N$-body Hamilton operator,
%
%
\begin{equation}\label{Hdef}
\Hcal_{N,\Lambda}^{(\rombc)} =-\sum_{i=1}^N \Delta_{i}^{(\rombc)}
+\sum_{1\leq i<j\leq N} v(|x_i-x_j|),\qquad x_1,\ldots,x_n\in\Lambda,
\end{equation}
where $\rombc\in\{\romDir, \romper\}$, and $\Delta_{i}^{(\rombc)}$
stands for the Laplacian with bc boundary condition.
More precisely, $Z_N^{(\rombc)}(\beta,\Lambda)$ is equal to the
trace of the projection of the operator $\exp\{-\beta\Hcal
_{N,\Lambda}^{(\rombc)}\}$ to the set of symmetric (i.e.,
permutation invariant) functions $(\R^d)^N\to\R$. This statement is
proven via the Feynman--Kac formula (see \cite{BR97} or \cite{G70}).
Hence, we call $Z_N^{(\rombc)}(\beta,\Lambda)$ a partition function.

It is the main purpose of this paper to derive a variational expression
for the \textit{limiting free energy}
%
%
\begin{equation}\label{limfreeener}
f^{(\rombc)}(\beta,\rho)=-\frac1\beta\lim_{N\to\infty
}\frac{1}{|\Lambda_{L_N}|}\log Z_{N}^{(\rombc)}(\beta,\Lambda_{L_N}),
\end{equation}
where $ |\Lambda_{L_N}|=N/\rho$, for any $\beta,\rho\in(0,\infty
)$, any $d\in\N$ and any $\rombc\in\{\varnothing, \romper,\romDir\}
$. The existence of the thermodynamic limit in (\ref
{limfreeener}) with $ \rombc\in\{ \romper,\romDir\} $ under
suitable assumptions on the interaction potential $v$ can be shown by
standard methods (see, e.g., \cite{Rue69}, Theorem 3.58, and \cite{R71}).
However, to the best of our knowledge, there is no useful
identification or characterization of $f^{(\rombc)}(\beta,\rho
)$ available in the literature. We also give new proofs for the
independence of the value of the free energy on the boundary
conditions, which is another novelty.

Our approach, and the remainder of Section \ref{Intro}, can be
summarized as follows. Since any permutation decomposes into cycles,
and using the Markov property, the family of the $N$ bridges in (\ref
{defpartition}) decomposes into cycles of various lengths, that is, into
bridges that start and end at the same site, which is uniformly
distributed over~$\Lambda$. We conceive these initial-terminal sites
as the points of a standard Poisson point process on $\R^d$ and the
cycles as marks attached to these points (see Section~\ref
{sec-notation} for the relevant notation). In Proposition \ref
{lem-rewrite} below we rewrite $Z_{N}^{(\rombc)}(\beta,\Lambda
)$ in terms of an expectation over a reference process, the marked
Poisson point process $\omega_{\romP}$.

In Section \ref{sec-limfreeen}, we present our results on the
large-$N$ asymptotics of $Z_{N}^{(\rombc)}(\beta,\Lambda)$ when
$\Lambda$ is a centered cube of volume $N/\rho$. Indeed, in
Theorem \ref{thm-mainres}, its exponential rate is bounded from above
and below in terms of two variational formulas that range over marked
shift-invariant point processes and optimize the sum of an energy term
and an entropy term. These bounds are shown to coincide for any fixed
$\rho$ if $\beta$ is sufficiently small, and for any fixed $\beta$
if $\rho$ is sufficiently small. The main value and novelty of these
representations is the explicit description of the interplay between
entropy, interaction and symmetrization of the system. We think that
these formulas, even in the case where our two bounds do not coincide,
are explicit enough to serve as a basis for future deeper
investigations of properties like phase transitions.

The physical interpretation, motivation and relevance are discussed in
Section~ \ref{sec-BEC}.

\subsection{Representation of the partition function}\label{sec-notation}

In this section, we introduce our representation of the
partition function $Z_{N}^{(\rombc)}(\beta,\Lambda)$ for each
boundary condition $ \rombc\in\{\varnothing, \romper,\romDir\} $ in
terms of an expectation over a marked Poisson point process.
The main result of this section is Proposition \ref{lem-rewrite}. We
have to introduce some notation.

We begin with the mark space. The space of marks is defined as
%
%
\begin{equation}
E^{(\rombc)}=\bigcup_{k\in\N} \Ccal^{(\rombc)}_{k,\Lambda},\qquad
\rombc\in\{\varnothing,\romper, \romDir\},
\end{equation}
where, for $k\in\N$, we denote by $\Ccal_k=\Ccal^{(\varnothing
)}_{k,\Lambda}$ the set of continuous functions $f\dvtx[0,k\beta]
\to
\mathbb{R}^d$ satisfying $f(0)=f(k\beta)$, equipped with the topology
of uniform convergence. Moreover,
$ \Ccal^{(\romDir)}_{k,\Lambda} $, respectively, $\Ccal
_{k,\Lambda}^{(\romper)}$, is the space of continuous functions
in $ \Lambda$, respectively, on the torus $\Lambda=(\R\slash L \Z
)^d$, with time horizon $ [0,k\beta] $. We sometimes call the marks
\textit{cycles}. By $\ell\dvtx E^{(\rombc)}\to\N$ we denote the
canonical map defined by $\ell(f)=k$ if $f\in\Ccal^{(\rombc)}_{k,
\Lambda}$. We call $\ell(f)$ the \textit{length} of $f\in E$.
When dealing with the empty boundary condition, we sometimes drop the
superscript $\varnothing$.

We consider spatial configurations that consist of a locally finite set
$\xi\subset\R^d$ of particles, and to each particle $x\in\xi$ we
attach a mark $f_x\in E^{(\rombc)}$ satisfying $f_x(0)=x$. Hence,
a configuration is described by the counting measure
\[
\omega=\sum_{x\in\xi}\delta_{(x,f_x)}
\]
on $\R^d\times E$ for the empty boundary condition, respectively, on
$\Lambda\times E^{(\rombc)}$ for $\rombc\in\{\romper, \romDir\}$.

We now introduce three marked Poisson point processes for the three
boundary conditions. The one for the empty condition will later serve
as a reference process and is introduced separately first.

\subsubsection*{Reference process}

Consider on $\Ccal=\Ccal_1$ the canonical Brownian bridge measure
%
%
\begin{equation}\label{nnBBM}\qquad
\mu^{({\varnothing, \beta})}_{x,y}(A)=\mu^{(\beta)
}_{x,y}(A)=\frac{\PP_x(B\in A;B_\beta\in\dd y)}{\dd y},\qquad
A\subset\Ccal\mbox{ measurable}.
\end{equation}
Here $B=(B_t)_{t\in[0,\beta]}$ is a Brownian motion in $\R^d$ with
generator $\Delta$, starting from $x$ under $\PP_x$. Then $ \mu
^{(\beta)}_{x,y}$ is a regular Borel measure on $\Ccal$ with
total mass equal to the Gaussian density
%
%
\begin{equation}\label{Gaussian}
\mu_{x,y}^{(\beta)}(\Ccal)=g_{\beta}(x,y)=\frac{\PP_x(B_\beta
\in\dd y)}{\dd y}=(4\pi\beta)^{-d/2}\ee^{-1/({4\beta})|x-y|^2}.
\end{equation}
We write $ \PP_{x,y}^{(\beta)}=\mu_{x,y}^{(\beta)}/g_{\beta
}(x,y)$ for the normalized Brownian bridge measure on~$\Ccal$.
Let
\[
\omega_{\romP} = \sum_{x \in\xi_{\romP}} \delta_{(x,B_x)}
\]
be a Poisson point process on $\R^d\times E$ with intensity measure
equal to $ \nu$ whose projection onto $ \R^d\times\Ccal_k $ is
equal to
%
%
\begin{equation}\label{nudef}
\nu_k(\dd x,\dd f)=\frac{1}{k}\Leb(\dd x)\otimes\mu_{x,x}^{
({k\beta})}(\dd f),\qquad k\in\N.
\end{equation}
Alternatively, we can conceive $\omega_{\romP}$ as a marked Poisson
point process on $\R^d$, based on some Poisson point process $\xi
_{\romP}$ on $\R^d$, and a family $(B_x)_{x\in\xi_{\romP}}$ of
i.i.d. marks, given $\xi_{\romP}$. The intensity of $\xi_{\romP}$ is
%
%
\begin{equation}\label{q*def}
\overline{q}=\sum_{k\in\N}q_k\qquad \mbox{with } q_k=\frac
1{(4\pi\beta)^{d/2}k^{1+d/2}},\qquad k\in\N.
\end{equation}
Conditionally given $\xi_{\romP}$, the length $\ell(B_x)$ is an $\N
$-valued random variable with distribution $(q_k/\overline{q})_{k\in
\N}$, and, given $\ell(B_x)=k$, $B_x$ is in distribution equal to a
Brownian bridge with time horizon $[0,k\beta]$, starting and ending at
$x$. Let $\ttQ $ denote the distribution of $\omega_{\romP}$, and
denote by $\ttE$ the corresponding expectation. Hence, $\ttQ $
is a
probability measure on the set $\Omega$ of all locally finite counting
measures on $ \R^d\times E $.

\subsubsection*{Processes for Dirichlet and periodic boundary conditions}

For Dirichlet boundary condition, one restricts the Brownian
bridges to not leaving the set $\Lambda$. Consider the measure
%
%
\begin{equation}\label{nnBBMDir}
\mu^{({\romDir,\beta})}_{x,y}(A)=\frac{\PP_x(B\in A;B_\beta
\in\dd y)}{\dd y},\qquad A\subset\Ccal^{(\romDir)}_{1,\Lambda
}\mbox{ measurable},
\end{equation}
which has total mass
%
%
\begin{equation}
g_\beta^{(\romDir)}(x,y)=\mu_{x,y}^{(\romDir,\beta
)}\bigl(\Ccal^{(\romDir)}_{1,\Lambda}\bigr)=\frac{\PP_x(B_{[0,\beta
]}\subset\Lambda; B_\beta\in\dd y)}{\dd y}.
\end{equation}

For periodic boundary condition, the marks are Brownian bridges on the
torus $\Lambda=(\R\slash L \Z)^d$. The corresponding path measure is
denoted by $\mu_{x,y}^{(\romper,\beta)}$; its total mass is
equal to
%
%
\begin{eqnarray}
g_\beta^{(\romper)}(x,y)&=&\mu_{x,y}^{(\romper,\beta
)}\bigl(\Ccal^{(\romper)}_\Lambda\bigr)=\sum_{z\in\Z^d}g_{\beta
}(x,y+zL)\nonumber\\[-8pt]\\[-8pt]
&=&(4\pi\beta)^{-d/2}\sum_{z\in\Z^d}\ee^{-
{|x-y-zL|^2}/({4\beta})}.\nonumber
\end{eqnarray}
For periodic and Dirichlet boundary conditions (\ref{q*def}) is
replaced by
%
%
\begin{equation}\label{qdefbc}
\overline{q}{}^{(\rombc)}=\sum_{k=1}^N q_k^{(\rombc)}\qquad \mbox
{with } q_k^{(\rombc)}=\frac{1}{k|\Lambda
|}\int_\Lambda\dd x\, g_{k\beta}^{(\rombc)}(x,x).
\end{equation}
Note that this weight depends on $ \Lambda$ and on $N$. We introduce
the Poisson point process $\omega_{\romP} = \sum_{x \in\xi_
{\romP}} \delta_{(x,B_x)}$ on $\Lambda\times E^{(\rombc)}$ with
intensity measure $ \nu^{(\rombc)} $ whose projections on $
\Lambda\times\Ccal_{k,\Lambda}^{(\rombc)} $ with $k\leq N$
are equal to $ \nu_k^{(\rombc)}(\dd x,\dd f)=\frac{1}{k}\Leb
_\Lambda(\dd x)\otimes\mu_{x,x}^{({\rombc,k\beta})}(\dd f)$
and are zero on this set for $k>N$. We do not label $\omega_{\romP}$
nor $\xi_{\romP}$ with the boundary condition nor with $N$; $\xi_{\romP
}$ is a Poisson process on $\Lambda$ with intensity measure
$\overline q^{(\rombc)}$ times the restriction $\Leb_\Lambda$
of the Lebesgue measure to $\Lambda$. By $\ttQ ^{(\rombc)}$
and $\ttE^{(\rombc)}$ we denote probability and expectation
with respect to this process. Conditionally on $\xi_{\romP}$, the
lengths of the cycles $B_x$ with $x\in\xi_{\romP}$ are independent
and have distribution $(q_k^{(\rombc)}/\overline q^{(\rombc)})_{k\in\{
1,\ldots,N\}}$; this process has only marks with lengths
$\leq N$. A cycle $B_x$ of length $k$ is distributed according to
%
%
\begin{equation}\label{pathprobmeasure}
\PP_{x,x}^{({\rombc,k\beta})}(\dd f)=\frac{\mu_{x,x}^{
({\rombc,k\beta})}(\dd f)}{g^{(\rombc)}_{k\beta}(x,x)}.
\end{equation}

We now formulate our first main result, a presentation of the partition
function defined in (\ref{defpartition}) in $\Lambda\subset\R^d$
with $ |\Lambda|<\infty$ and boundary condition $ \rombc\in\{
\varnothing,\romper,\romDir\}$. We write $\langle P,F\rangle$ for
the expectation of a function $F$ with respect to a probability measure
$P$. We introduce a functional on $\Omega$ that expresses the
interaction between particles in $ \Lambda\subset\R^d$, more
precisely, between their marks. Define the \textit{Hamiltonian} $
H_\Lambda\dvtx\Omega\to[0,\infty] $ by
%
%
\begin{equation}\label{Hamiltonian}
H_\Lambda(\omega)=\sum_{x,y\in\xi\cap\Lambda}T_{x,y}(\omega
)\qquad \mbox{where }\omega=\sum_{x\in\xi}\delta_{(x,f_x)}\in
\Omega,
\end{equation}
where we abbreviate, for $\omega\in\Omega,x,y\in\xi$,
%
%
\begin{eqnarray}\label{Tdef}\quad
T_{x,y}(\omega)&=&\frac{1}{2}\sum_{i=0}^{\ell(f_x)-1}\sum
_{j=0}^{\ell(f_y)-1}\mathbh{1}_{\{(x,i)\not=(y,j)\}}\nonumber\\[-8pt]\\[-8pt]
&&\hspace*{68.5pt}{}\times\int_0^\beta
v\bigl(|f_x(i\beta+s)-f_y(j\beta+s)|\bigr) \,\dd s.\nonumber
\end{eqnarray}

The function $H_\Lambda(\omega)$ summarizes the interaction between
different marks of the point process and between different legs of the
same mark; here we call the restriction of a mark $f_x$ to the interval
$[i\beta,(i+1)\beta)]$ with $i\in\{0,\ldots,\ell(f_x)-1\}$ a leg of
the mark. Denote by
%
%
\begin{equation}\label{Nlength}
N^{(\ell)}_{\Lambda}(\omega)=\sum_{x\in\xi\cap\Lambda}\ell(f_x)
\end{equation}
the total length of the marks of the particles in $ \Lambda\subset\R
^d $ (whose marks may be not contained in $\Lambda$).
\begin{prop}[(Rewrite in terms of the marked Poisson process)]\label
{lem-rewrite} Fix $\beta\in(0,\infty)$. Let $ v\dvtx[0,\infty)\to
(-\infty,\infty] $ be measurable and bounded from below, and let $
\Lambda\subset\R^d$ be measurable with finite volume (assumed to be
a torus for periodic boundary condition). Then, for any $N\in\N$, and
$ \rombc\in\{\varnothing,\romper,\romDir\} $,
%
%
\begin{equation}\label{rewrite}
Z_N^{(\rombc)}(\beta,\Lambda)=\ee^{|\Lambda
|\overline
{q}^{(\rombc)}}\ttE^{(\rombc)}\bigl[\ee^{-H_\Lambda(\omega
_{\romP})}\mathbh{1}\bigl\{N^{(\ell)}_\Lambda(\omega_{\romP})=N\bigr\}\bigr].
\end{equation}
\end{prop}

That is, up the nonrandom term $|\Lambda|\overline{q}^{(\rombc)}$,
the partition function is equal to the expectation over the
Boltzmann factor $ \ee^{-H_\Lambda} $ of a marked Poisson process
with fixed total length of marks of the particles.

\subsection{The limiting free energy}\label{sec-limfreeen} In this
section, we present our major result, the identification of the
limiting free energy defined in (\ref{limfreeener}) in terms of an
explicit variational formula (see Theorem \ref{thm-mainres}). We first
introduce some notation.

Define the shift operator $ \theta_y \dvtx\R^d \to\R^d$ as $\theta
_y(x)=x-y$. We extend it to a shift operator on marked configurations by
\[
\theta_y(\omega)=\sum_{x\in\xi}\delta_{(x-y, f_x)}=\sum_{x\in
\xi-y}\delta_{(x, f_{x+y})}\qquad \mbox{for }\omega=\sum_{x\in
\xi}\delta_{(x, f_x)}.
\]
By $\mathcal{P}_{\theta}$ we denote the set of all shift-invariant
probability measures on $\Omega$. The distribution $\ttQ $ of the
above marked Poisson point reference process $\omega_{\romP}$ belongs
to~$\mathcal{P}_{\theta}$.

Define $\Phi_\beta\dvtx\Omega\to[0,\infty]$ by
%
%
\begin{equation}\label{Gammabetadef}
\Phi_\beta(\omega) = \sum_{x \in\xi\cap U} \sum_{y \in\xi
}T_{x,y}(\omega),
\end{equation}
where $T_{x,y}(\omega)$ was defined in (\ref{Tdef}), and $U=[-\frac
{1}{2},\frac{1}{2}]^d$ denotes the centered unit box. The quantity
$\Phi
_\beta(\omega)$ describes all the interactions between different legs
of marks of $\omega$, when at least one of the marks is attached to a
point in~$U$.

Next, we introduce an entropy term. For probability measures $\mu, \nu
$ on some measurable space, we write
%
%
\begin{equation}
H(\mu\mid\nu)=\cases{\displaystyle  \int f\log f \,\dd\nu, &\quad if $\displaystyle f=\frac{\dd
\mu}{\dd\nu}$ exists,
\cr
\infty, &\quad otherwise,}
\end{equation}
for the relative entropy of $\mu$ with respect $\nu$. It will be
clear from the context which measurable space is used. It is easy to
see and well known that $H(\mu\mid\nu)$ is nonnegative and that it
vanishes if and only if $\mu=\nu$. Now we set
%
%
\begin{equation}\label{Idef}
I_\beta(P) = \lim_{N \to\infty} \frac{1}{|\Lambda_{N}|} H(
P_{\Lambda_N} | \ttQ _{\Lambda_N}),\qquad P\in\Pcal
_\theta,
\end{equation}
where we write $P_{\Lambda}$ for the projection of $P$ to $\Lambda$,
that is, the image measure of $P$ under
%
%
\begin{equation}\label{projection}
\omega\mapsto\omega|_{\Lambda}=\sum_{x\in\xi\cap\Lambda}\delta
_{(x, f_x)}\qquad \mbox{for }\omega=\sum_{x\in\xi}\delta_{(x, f_x)}.
\end{equation}
The limit in (\ref{Idef}) is along centered boxes $\Lambda_N$ with
diverging volume.
According to \cite{GZ93}, Proposition 2.6, the limit in (\ref{Idef})
exists, and $I_\beta$ is a lower semicontinuous function with compact
level sets in the topology of local convergence (see Lemma \ref{ldpyn}
below). It turns out there that $I_\beta$ is the rate function of a
crucial large-deviations principle for the family of the stationary
empirical fields, which is one of the important objects of our analysis
and will be introduced at the beginning of Section \ref{sec-LDP}.

Now we introduce two important variational formulas. For any $\beta
,\rho\in(0,\infty)$, define
%
%
\begin{eqnarray}
\label{chiledef}
\chi^{(\le)}(\beta,\rho)&=&\inf\bigl\{I_\beta(P)+\langle
P,\Phi_\beta\rangle\dvtx P\in\Pcal_\theta, \bigl\langle P,N_U^{
(\ell)}\bigr\rangle\le\rho\bigr\},\\
\label{chi=def}
\chi^{(=)}(\beta,\rho)&=&\inf\bigl\{I_\beta(P)+\langle
P,\Phi_\beta\rangle\dvtx P\in\Pcal_\theta, \bigl\langle P,N_U^{
(\ell)}\bigr\rangle=\rho\bigr\}.
\end{eqnarray}
These formulas range over shift-invariant marked processes $P$. They
have three components: the entropic distance $I_\beta(P)$ between $P$
and the reference process $\ttQ $, the interaction term $\langle
P,\Phi_\beta\rangle$ and the condition $\langle P,N_U^{(\ell)
}\rangle=\rho$, respectively, \mbox{$\le$}$\rho$. Obviously, $\chi^{
(\le)}\leq\chi^{(=)}$. Since all the maps $P\mapsto I_\beta
(P)$, $P\mapsto\langle P,\Phi_\beta\rangle$ and $P\mapsto\langle
P,N_U^{(\ell)}\rangle$ are easily seen to be lower semicontinuous
and since the level sets of $I_\beta$ are compact, it is clear that
the infimum on the right-hand side of (\ref{chiledef}) is attained and
is therefore a minimum. However, this is not at all clear for (\ref
{chi=def}); this question lies much deeper and has some relation to the
question about Bose--Einstein condensation (see the discussion in
Section \ref{sec-BEC}).

Now we specify our assumptions on the particle interaction potential $ v$.
\renewcommand{\theAssumption}{(v)}
\begin{Assumption}\label{AssumptionV}
We assume that $ v\dvtx
[0,\infty)\to[0,\infty] $ is measurable and tempered, that is, there
are $h>d , A\ge0 $ and $ R_0>0 $ such that $ v(t)\le At^{-h} $ for $
t\in[R_0,\infty) $. Additionally, we assume that the integral
\[
\alpha(v)= \int_{\R^d}v(|x|) \,\dd x
\]
is finite and that $ \liminf_{r\to0} v(r)>0 $.
\end{Assumption}

We now present variational characterizations for upper and lower bounds
for the exponential rate of the partition function. We denote by
$\Lambda_L=[-\frac L2,\frac L2]^d$ the centered box in $\R^d$ with
volume $L^d$.
\begin{theorem}\label{thm-mainres} Let $L_N=(\frac N\rho)^{1/d}$,
such that $\Lambda_{L_N}$ has volume $N/\rho$. Let $ v$ satisfy
Assumption \ref{AssumptionV}. Denote
%
%
\begin{equation}\label{Dcaldef}
\Dcal_v=\bigl\{(\beta,\rho)\in(0,\infty)^2\dvtx(4\pi\beta
)^{-d/2}\geq\rho\ee^{\beta\rho\alpha(v)}\bigr\}.
\end{equation}
Then, for any $\beta,\rho\in(0,\infty)$, and for $ \rombc\in\{
\varnothing,\romDir, \romper\} $,
%
%
\begin{equation} \label{freeenupper}
\limsup_{N \to\infty} \frac{1}{|\Lambda_{L_N}|} \log Z^{(\rombc)}
_N(\beta,\Lambda_{L_N})\le\frac{\zeta(1+{d/2})}{(4\pi
\beta)^{d/2}}-\chi^{(\le)}(\beta,\rho)
\end{equation}
and
%
%
\begin{eqnarray}\label{freeenlower}
&&\liminf_{N \to\infty} \frac{1}{|\Lambda_{L_N}|} \log Z^{(\rombc)}
_N(\beta,\Lambda_{L_N})\nonumber\\[-8pt]\\[-8pt]
&&\qquad\ge\frac{\zeta(1+{d/2})}{(4\pi\beta)^{d/2}}-
\cases{ \chi^{(\le)}(\beta,\rho), &\quad if
$(\beta,\rho)\in\Dcal_v$,\cr
\chi^{(=)}(\beta,\rho), &\quad if
$(\beta,\rho)\notin\Dcal_v$,}\nonumber
\end{eqnarray}
where $\zeta(m)=\sum_{k=1}^\infty k^{-m}$ denotes the Riemann zeta function.
\end{theorem}

Note that the first term on the right, $\zeta(1+\frac{d}{2})/(4\pi
\beta)^{d/2}$, is equal to the total mass~$ \overline{q}$, the sum of
the $q_k$ defined in (\ref{q*def}). The proof of Theorem \ref
{thm-mainres} is in Sections~\ref{sec-upbound} [proof of (\ref
{freeenupper})] and \ref{sec-lowbound} [proof of (\ref{freeenlower})]
for empty boundary conditions, and in Section \ref{otherbc} for the
other two.

The assumptions $ \int_{\R^d}v(|x|) \,\dd x<\infty$ and $\liminf_{r
\to0} v(r)>0$ are only necessary for our proof of the lower bound in
(\ref{freeenlower}). In the proof of the upper bound in (\ref
{freeenupper}), it is allowed that $v$ takes the value $ +\infty$ on a
set of positive measure (corresponding to hard core repulsion) and also
that $v\equiv0$ (the noninteracting case) (see discussion in
Section \ref{sec-non-interact}).

As an obvious corollary we now identify the free energy defined in
(\ref{limfreeener}) in the high-temperature phase and in the
low-density phase.
\begin{cor}[(Free energy)]\label{limitexists} Fix $(\beta
,\rho)\in\Dcal_v$. Then, for any $
\rombc\in\{\varnothing,\romDir$, $\romper\} $,
the free energy introduced in (\ref{limfreeener})
is given by
%
%
\begin{eqnarray} \label{freeenident}
f(\beta,\rho)&=&f^{(\rombc)}(\beta,\rho)\nonumber\\
&=&-\frac1\beta\frac{\zeta(1+{d/2})}{(4\pi\beta
)^{d/2}}\\
&&{} +\frac
1\beta\min\bigl\{I_\beta(P)+\langle P,\Phi_\beta\rangle\dvtx P\in
\Pcal_\theta, \bigl\langle P,N_U^{(\ell)}\bigr\rangle\le\rho\bigr\}
.\nonumber
\end{eqnarray}
\end{cor}

A by-product of the proof of the lower bound of (\ref{freeenlower})
(see Corollary \ref{cor-upf}) we have the following upper bound on the
free energy.
\begin{lemma} For any $\beta,\rho\in(0, \infty)$, and for $ \rombc\in\{
\varnothing,\romDir, \romper\} $,
%
%
\begin{eqnarray}
f^{(\rombc)}(\beta,\rho)&=&\limsup_{N \to\infty}-\frac1\beta
\frac{1}{|\Lambda_{L_N}|} \log Z^{(\rombc)}_N(\beta,\Lambda
_{L_N})\nonumber\\[-8pt]\\[-8pt]
&\le&\frac{\rho}{\beta} \log( \rho(4 \pi\beta)^{d/2})
+\rho^2\alpha(v).\nonumber
\end{eqnarray}
\end{lemma}

\subsection{Relevance and discussion}\label{sec-BEC}

One of the most prominent open problems in mathematical
physics is the understanding of \textit{Bose--Einstein condensation}
(BEC), a phase transition in a mutually repellent many-particle system
at positive, fixed particle density, if a sufficiently low temperature
is reached. That is, a macroscopic part of the system condenses to a
state which is highly correlated and coherent. The first experimental
realization of BEC was only in 1995, and it has been awarded with a
Nobel prize. In spite of an enormous research activity, this phase
transition has withstood a mathematical proof yet. Only partial
successes have been achieved, like the description of the free energy
of the ideal, that is, noninteracting, system (already contained in
Bose and Einstein's seminal paper in 1925) or the analysis of
mean-field models (e.g., \cite{Toth90,DMP05}) or the analysis of
dilute systems at vanishing temperature \cite{LSSY05} or the proof of
BEC in lattice systems with half-filling~\cite{LSSY05}. However, the
original problem for fixed positive particle density and temperature is
still waiting for a promising attack. Not even a tractable formula for
the limiting free energy was known yet that could serve as a basis for
a proof of BEC. The main purpose of the present paper is to provide
such a formula.

The mathematical description of bosons is in terms of the symmetrized
trace of the negative exponential of the corresponding Hamiltonian
times the inverse temperature. The symmetrization creates long range
correlations of the interacting particles making the analysis an
extremely challenging endeavor. The Feynman--Kac formula gives, in a
natural way, a representation in terms of an expansion with respect to
the cycles of random paths. It is conjectured by Feynman \cite{F53}
that BEC is signaled by the decisive appearance of a macroscopic amount
of ``infinite'' cycles, that is, cycles whose lengths diverge with the
number of particles. This phenomenon is also signaled by a loss of
probability mass in the distribution of the ``finite'' cycles. See \cite
{S93} and \cite{S02} for proofs of this coincidence in the ideal Bose
gas and some mean-field models. A different line of research is
studying the effect of the symmetrization in random permutation and
random partition models (see \cite{Ver96,BCMP05,AD06,AK08,A07},
or in spatial random permutation models going back to
\cite{F91} and extended in \cite{BUe09}).

In the present paper, we address the original problem of a mutually
repellent many-particle system at fixed positive particle density and
temperature and derive an explicit variational expression for the
limiting free energy. More precisely, we prove upper and lower bounds,
which coincide in the high-temperature phase, respectively, low density
phase. The formula yields deep inside in the cycle structure of the
random paths appearing in the Feynman--Kac formula. In particular, it
opens up a new way to analyze the structure of the cycles at any
temperature and density, also in the low-temperature phase, where our
two bounds differ. In future work, we intend to analyze the conjectured
phase transition in that variational formula and to link it to BEC.

The methods used in the present paper are mainly probabilistic. Our
starting point is the well-known Feynman--Kac formula, which translates
the partition function in terms of an expectation over a large
symmetrized system of interacting Brownian bridge paths. In a second
step, which is also well known, we reduce the combinatorial complexity
by concatenating the bridges using the symmetrization. The novelty of
the present approach is a reformulation of this system in terms of an
expectation with respect to a \textit{marked Poisson point process}, which
serves as a reference process. This is a Poisson process in the space
$\R^d$ to whose particles we attach cycles called marks, starting and
ending at that particle. The symmetrization is reflected by an a priori
distribution of cycle lengths. The interaction between the Brownian
particles are encoded as interaction between the marks in an
exponential functional. The particle density is described by a
condition on the total length of the marks in the unit box.

Approaches to Bose gases using point processes have occasionally been
used in the past (see \cite{F91} and the references therein) and also
recently in \cite{R09}, but systems with interactions have not yet
been considered using this technique, to the best of our knowledge.

The greatest advantage of this approach is that it is amenable to a
large-deviations analysis. The central object here is the
\textit{stationary empirical field} of the marked point process, which contains
all relevant information and satisfies a large-deviations principle in
the thermodynamic limit. For some class of interacting systems, this
direction of research was explored in \cite{GZ93,G94}. In the present
paper, we apply these ideas to the more difficult case of the
interacting Bose gas. The challenge here is that the interaction
involves the spatial points and the details of the marks. Modulo some
error terms, we express the interaction and the mark length condition
in terms of a functional of the stationary empirical field. Formally
using Varadhan's lemma, we obtain a variational formula in the limit.

However, due to a lack of continuity in the functionals that describe
the interaction and the mark lengths, the upper and lower bounds
derived in this way, may differ in general. (At sufficiently high
temperature, we overcome this problem by additional efforts and
establish a formula for the limit.) This effect is not a technical
drawback of the method, but lies at the heart of BEC.

In Theorem \ref{thm-mainres}, we formulate the limiting free energy in
terms of a minimizing problem for random shift-invariant marked point
processes with interaction under a constraint on the total length of
the marks per unit volume. Both formulas in our upper and lower bounds
in Theorem \ref{thm-mainres} are formulated in terms of random point
fields having \textit{finitely long} cycles as marks. The concept used in
the present paper is not able to incorporate infinitely long cycles nor
to quantify their contribution to the interaction. In the proof of our
lower bound of the free energy, we drop the interactions involving any
cycle longer than a parameter $R$ that is eventually sent to infinity,
and in our proof of the upper bound we even drop these cycles in the
probability space. As a result, our two formulas register only
``finitely long'' cycles. Their total macroscopic contribution is
represented by the term $\langle P,N_U^{(\ell)}\rangle$, and the
one of the ``infinitely long'' cycles by the term $\rho-\langle
P,N_U^{(\ell)}\rangle$. In this way, the long cycles are only
indirectly present in our analysis: in terms of a ``loss of mass,'' the
difference between the particle density $\rho$ and the total mass of
short cycles. Physically speaking, this difference is the total mass of
a \textit{condensate} of the particles.

The values of the two formulas $\chi^{(\le)}(\beta,\rho)$ and
$\chi^{(=)}(\beta,\rho)$ differ if ``infinitely long'' cycles do
have some decisive contribution in the sense that the optimal point
process(es) $P$ in $\chi^{(\leq)}(\beta,\rho)$ satisfies
$\langle P,N_U^{(\ell)}\rangle< \rho$. We conjecture that the
question whether or not the optimal $P$ in $\chi^{(\leq)}(\beta
,\rho)$ has a loss of probability mass of infinitely long cycles is
intimately related with the question whether or not $\chi^{(\leq)
}(\beta,\rho)=\chi^{(=)}(\beta,\rho)$ and that this question
is in turn decisively connected with the question whether or not BEC
appears. This is in accordance with S\"ut\H o's work \cite{S93,S02}.
The conjecture is that, for given $\beta$ and in $d\geq3$, if $\rho$
is sufficiently small, then it is satisfied, and for sufficiently large
$\rho$ it is not satisfied. The latter phase is conjectured to be the
BEC phase. Future work will be devoted to an analysis of this question.

Here is an abstract sufficient criterion for $\chi^{(\leq)}(\beta
,\rho)=\chi^{(=)}(\beta,\rho)$.
\begin{lemma}Fix $\beta\in(0,\infty)$. If there exists a minimizer
$\widehat P $ of the variational problem $\inf_{P\in\Pcal_\theta
}(I_\beta(P)+\langle P,\Phi_\beta\rangle)$ satisfying $ \widehat
\rho:=\langle\widehat P,N_U^{(\ell)}\rangle<\infty$, then, for
any $ \rho\in(0,\widehat\rho) $,
%
%
\begin{equation}
\chi^{(\leq)}(\beta,\rho)=\chi^{(=)}(\beta, \rho).
\end{equation}
\end{lemma}
\begin{pf}
Pick $ \rho<\widehat\rho$. Let $ P $ be a minimizer in the formula
for $\chi^{(\leq)}(\beta,\rho)$, that is, of $ \inf\{I_\beta
(P)+\Phi_\beta(P)\dvtx P\in\Pcal_\theta,\langle P,N_U^{(\ell)
}\rangle\le\rho\} $. If $ \langle P,N_U^{(\ell)}\rangle$ would
be smaller than $\rho$, then an appropriate convex combination, $
\widetilde P $, of $P$ and $\widehat P$ would satisfy $ \langle
\widetilde P,N_U^{(\ell)}\rangle\in(\langle P,N_U^{(\ell)
}\rangle,\rho] $ and $ I_\beta(\widetilde P)+\Phi_\beta(\widetilde
P)< I_\beta(P)+\Phi_\beta(P) $. This would contradict the minimizing
property of $P$. Hence, $ \langle P,N_U^{(\ell)}\rangle=\rho$,
and therefore $P$ minimizes also the formula for $\chi^{(=)}(\beta
, \rho)$.
\end{pf}

\subsection{The noninteracting case}\label{sec-non-interact}

Let us compare our results to the noninteracting case.
Indeed, \cite{A07}, Theorem 2.1, says that, in the case $v\equiv0$,
the identification of the limiting free energy in (\ref{freeenident})
holds for \textit{any} $\beta,\rho\in(0,\infty)$. To see this, we have
to argue a bit, and we will only sketch the argument.

Explicitly, after applying some elementary manipulations, one sees that
\cite{A07}, Theorem 2.1, amounts to
%
%
\begin{equation}\label{SAcite}
f(\beta,\rho)=-\frac1\beta\frac{\zeta(1+{d/2})}{(4\pi
\beta
)^{d/2}}+\frac1\beta\inf_{\lambda\in\ell^1(\N)\dvtx\sum_k k
\lambda_k\leq1}J(\lambda),
\end{equation}
where we recall that $q$ was defined in (\ref{q*def}), and we put
\[
J(\lambda)=\sum_{k\in\N} q_k+\rho H(\lambda\mid q)+\rho\sum
_{k\in\N} \lambda_k\log\rho-\rho\sum_{k\in\N} \lambda_k.
\]
Now we rewrite the minimum on the right-hand side of (\ref
{freeenident}) in a similar form by splitting $N_U^{(\ell)}$ into
$\sum_{k\in\N}k \Ncal_k$, where
%
%
\begin{equation}\label{NLambdakdef}
\Ncal_{k,\Lambda}(\omega)=\#\{x \in\xi\cap\Lambda\dvtx\ell
(f_x)=k\}
\end{equation}
and $\Ncal_k=\Ncal_{k,U}$ is the number of particles in the unit box
$U$ whose cycles have length $k$ (and are allowed to leave $U$). Then
we may write
\[
\inf\bigl\{I_\beta(P)\dvtx P\in\Pcal_\theta,\bigl\langle P,N_U^{
(\ell)}\bigr\rangle\leq\rho\bigr\}
=\inf_{\lambda\in\ell^1(\N)\dvtx\sum_k k \lambda_k\leq1}\inf
_{P\in\Pcal_\theta\dvtx\lambda(P)=\lambda} I_\beta(P),
\]
where $\lambda(P)=\frac1\rho(\langle P, \Ncal_k\rangle)_{k\in\N
}$. In order to see that (\ref{SAcite}) coincides with (\ref
{freeenident}) for $v=0$, one only has to check that $J(\lambda)=\inf
_{P\in\Pcal_\theta\dvtx\lambda(P)=\lambda} I_\beta(P)$ for any
$\lambda\in\ell^1(\N)$ satisfying $\sum_k k \lambda_k\leq1$.

We do not offer an analytical proof of this fact, but instead a
probabilistic one, which makes use of the large-deviations principle in
Lemma \ref{ldpyn} below for the stationary empirical field $\RR
_{\Lambda_L,\omega_{\romP}}$ introduced in (\ref{RNdef}) with rate
function $I_\beta$. Observe that the mapping $P\mapsto\lambda(P)$ is
continuous as a function from the set of all $P\in\Pcal_\theta$
satisfying $\langle P,N_U^{(\ell)}\rangle\leq\rho$ into the
sequence space $\ell^1(\N)$. Hence, by the contraction principle (see
\cite{DZ98}, Theorem 4.2.1), the sequence $(\lambda(\RR_{\Lambda
_L,\omega_{\romP}}))_{L>0}$ satisfies a large-deviations principle
with rate function $\lambda\mapsto\inf_{P\in\Pcal_\theta\dvtx
\lambda(P)=\lambda} I_\beta(P)$. By uniqueness of rate functions, it
suffices to show that this sequence satisfies the principle with rate
function $J$. We now indicate how to derive this by explicit calculation.

Introduce
\[
M_\Lambda=\biggl\{\lambda\in[0,1]^\N\dvtx\sum_k k\lambda_k\leq1,
\forall k\in\N\dvtx\lambda_k |\Lambda|\rho\in\N_0\biggr\},
\]
and for $\lambda\in M_\Lambda$, we calculate
\begin{eqnarray*}
\ttQ \bigl(\lambda(\RR_{\Lambda,\omega_{\romP}}
)=\lambda\bigr)&=&\ttQ (\forall k\in\N\dvtx\langle\RR
_{\Lambda,\omega_{\romP}},\Ncal_k\rangle=\rho\lambda_k)\\
&=&\ttQ \bigl(\forall k\in\N\dvtx\#\bigl(\xi_{\romP}^{(k)}\cap
\Lambda\bigr)=\rho|\Lambda|\lambda_k\bigr),
\end{eqnarray*}
where $\xi_{\romP}^{(k)}=\{x\in\xi_{\romP}\dvtx f_x\in\Ccal_k\}
$ is the set of those Poisson points with cycle of length $k$. Since
the Poisson processes $\xi_{\romP}^{(k)}$, $k\in\N$, are
independent with intensity~$q_k$, we can proceed with
\begin{eqnarray*}
\ttQ \bigl(\lambda(\RR_{\Lambda,\omega_{\romP}}
)=\lambda\bigr)
&=&\prod_{k\in\N}\ttQ \bigl(\#\bigl(\xi_{\romP}^{(k)}\cap\Lambda
\bigr)=\rho|\Lambda|\lambda_k\bigr)\\
&=&\prod_{k\in\N}\biggl(\ee^{-|\Lambda| q_k}\frac
{(|\Lambda
|q_k)^{\rho|\Lambda|\lambda_k}}{(\rho|\Lambda|\lambda_k)!}\biggr).
\end{eqnarray*}
Using Stirling's formula, we get from here that
\[
\frac1{|\Lambda_L|}\log\ttQ \bigl(\lambda(\RR_{\Lambda
_L,\omega_{\romP}})=\lambda\bigr)\sim-J(\lambda)\qquad
\lambda\in M_{\Lambda_L}\qquad \mbox{as }L\to\infty.
\]
From here, it is easy to finish the proof of the large-deviations
principle for $(\lambda(\RR_{\Lambda_L,\omega_{\romP}}))_{L>0}$ with
rate function $J$. This finishes the proof of (\ref
{freeenident}) for \textit{any} $\beta,\rho\in(0,\infty)$ in the
noninteracting case $v\equiv0$.

The well-known Bose--Einstein phase transition in the free energy was
made explicit in the analysis of the right-hand side of (\ref{SAcite})
in \cite{A07}. It was shown there that
%
%
\begin{eqnarray}\label{phasetransition}\hspace*{32pt}
f(\beta,\rho)&=&
-\frac{1}\beta\frac1{(4\pi\beta)^{d/2}}\nonumber\\[-8pt]\\[-8pt]
&&{}\times
\cases{\displaystyle
\sum_{k\in\N}\frac{\ee^{-\alpha k}}{k^{d/2+1}}+(4\pi
\beta
)^{d/2}\rho\alpha, &\quad if $\displaystyle \rho(4\pi\beta)^{d/2}<\zeta\biggl(\frac
{d}{2}\biggr)$,\cr
\displaystyle \zeta\biggl(1+\frac{d}{2}\biggr), &\quad if $\displaystyle \rho(4\pi\beta)^{d/2}\geq\zeta\biggl(\frac
{d}{2}\biggr)$,}\nonumber
\end{eqnarray}
where $ \alpha$ is the unique root of $\rho=(4\pi\beta)^{-d/2}\sum
_{k\in\N}\frac{\ee^{-\alpha k}}{k^{d/2}} $. Note that
$\zeta(\frac
d2)=\infty$ in $d\in\{1,2\}$, and hence there is no phase transition
in these dimensions. The first line in (\ref{phasetransition})
corresponds to the case where the minimizer $\lambda$ in (\ref
{SAcite}) satifies $\sum_k k \lambda_k=1$, that is, no ``infinitely
long'' cycles contribute to the free energy, and the second line to the
case $\sum_k k \lambda_k<1$. Hence, the Bose--Einstein phase
transition is precisely at the point where the variational formula in
(\ref{SAcite}) with ``$\leq$'' starts differing from the formula with
``$=$.''

\section{Rewrite of the partition function}\label{sec-ProofProp}

In this section, we give the proof of Proposition \ref{lem-rewrite}.

As a first step, we give a representation of $Z_N^{(\rombc)}(\beta
,\Lambda)$ in terms of an expansion with respect to the
cycles of the permutations in (\ref{defpartition}). This is well known
and goes back to Feynman 1955.

We denote the set of all \textit{integer partitions} of $N$ by
%
%
\begin{equation}\label{PNdef}
\mathfrak{P}_N=\biggl\{\lambda=(\lambda_k)_k\in\N_0^\N\dvtx\sum
_k k\lambda_k=N\biggr\}.
\end{equation}
The numbers $ \lambda_k$ are called the \textit{occupation numbers} of
the integer partition $ \lambda$. Any integer partition $\lambda$ of
$N$ defines a conjugacy class of permutations of $1,\ldots,N$ having
exactly $\lambda_k $ cycles of length $k$ for any $k\in\N$. The term
in (\ref{defpartition}) after the sum on $\sigma$ depends only on
this class. Hence, we replace this sum by a sum on integer partitions
$\lambda\in\mathfrak{P}_N$ and count the permutations in that class.
For any of these cycles of length $k$, we integrate out over all but
one of the starting and terminating points of all the $k$ Brownian
bridges belonging to that cycle and use the Markov property to
concatenate them. This gives the $i$th (with $i=1,\ldots,\lambda_k$)
bridge $B^{(k,i)}$ with time horizon $[0,k\beta]$, starting and
terminating at a site, which is uniformly distributed over $\Lambda$.
The family of these bridges $B^{(k,i)}$ is independent, and
$B^{(k,i)}$ has distribution $\PP_{\Lambda}^{({\rombc,k\beta})}$,
where we define
%
%
\begin{equation}\label{muLambdadef}
\PP_{\Lambda}^{({\rombc,\beta})}(\dd f)=\frac{\int_{\Lambda
} \dd x\, \mu^{({\rombc,\beta})}_{x,x}(\dd f)}{\int_\Lambda
\dd x\, g^{(\rombc)}_\beta(x,x)}.
\end{equation}
The expectation will be denoted by $ \E_\Lambda^{({\rombc,\beta})} $.

For $\lambda\in\mathfrak{P}_N$, define
%
%
\begin{eqnarray}\label{GWdef}\qquad
\mathcal{G}_{N,\beta}^{({\lambda})}&=& \frac12 \sum_{k_1, k_2 =1}^N
\sum_{i_1 =1}^{\lambda_{k_1}} \sum_{i_2 =1}^{\lambda_{k_2}} \sum
_{j_1 =0}^{k_1-1}\sum_{j_2= 0}^{k_2-1} \mathbh{1}_{(k_1, i_1,
j_1)\neq(k_2, i_2, j_2)}\nonumber\\
&&\hspace*{114.1pt}{} \times\int_0^\beta\dd s\, v\bigl(\bigl|B^{({k_1, i_1})}(j_1
\beta+ s) \\
&&\hspace*{168.3pt}{} - B^{({k_2, i_2})}(j_2 \beta+ s)
\bigr|\bigr).\nonumber
\end{eqnarray}
In words, $\mathcal{G}_{N,\beta}^{\lambda}$ is the total interaction
between different bridges $B^{({k_1,i_1})}$ and $B^{
({k_2,i_2})}$ and between different legs of the same bridge $B^{(k,i)}$.
\begin{lemma}[(Cycle expansion)]\label{lem-Cycle} For any $N\in\N$,
%
%
\begin{equation}\label{cycle}
Z_N^{({\rombc})}(\beta,\Lambda) = \sum_{\lambda\in\mathfrak{P}_N}
\biggl(\prod_{k\in\N}\frac{[\int_\Lambda\dd x\, g_{k\beta
}^{(\rombc)}
(x,x)]^{\lambda_k}}{\lambda_k! k^{\lambda_k}}\biggr)\bigotimes
_{k\in\N}\bigl(\E_{\Lambda}^{({\rombc,k\beta})}\bigr)^{\otimes\lambda_k}[\ee
^{-\mathcal{G}_{N,\beta}^{({\lambda})}}].\hspace*{-32pt}
\end{equation}
\end{lemma}
\begin{pf}
We are going to split every permutation on the right-hand side of
(\ref{defpartition}) into a product of its cycles. Assume that a
permutation $\sigma\in\Sym_N$ has precisely $\lambda_k$ cycles of
length $k$, for any $k\in\{1,\ldots,N\}$. Then
$\sum_{k=1}^Nk\lambda_k=N$. The corresponding Brownian bridges may be
renumbered $B^{(k,i)}_j$ with $k\in\N$, $i=1,\ldots, \lambda_k$
and $j=1,\ldots,k$.

Then the measure $\int_{\Lambda} \dd x_1 \cdots\int_{\Lambda}\dd
x_N \bigotimes_{i=1}^N \mu_{x_i, x_{\sigma(i)}}^{(\rombc,\beta)}$
splits into an according product, which can be written,
after a proper renumbering of the indices, as
%
%
\begin{equation}\label{prodmeasure}\quad
\prod_{k=1}^N \prod_{i=1}^{\lambda_k} \prod_{j=0}^{k-1} \int
_{\Lambda} \dd x_{k,j+1}^{(i)} \bigotimes_{k\in\N} \bigotimes
_{i=1}^{\lambda_k} \bigotimes_{j=0}^{k-1}\mu_{x_{k,j}^{(i)},
x_{k,j+1}^{(i)}}^{({\rombc,\beta})}\qquad \mbox{where
}x_{k,0}^{(i)}=x_{k,k}^{(i)}.
\end{equation}
Denote by $f_1 \diamond\cdots\diamond f_k$ the concatenation of
$f_1,\ldots,f_k$, that is, $f_1\diamond\cdots\diamond f_k((i-1)\beta
+s)=f_i(s)$ for $s\in[0,\beta]$. Note that the Markov property of the
canonical Brownian bridge measures implies the concatenation formula
%
%
\begin{equation}\label{concatenate}\quad
\mu^{({\rombc,k\beta})}_{x,x}\bigl(\dd(f_1\diamond\cdots\diamond f_k)\bigr)=
\int_{(\Lambda)^{k-1}}\dd x_1\cdots\dd x_{k-1} \bigotimes_{i=1}^k
\mu^{({\rombc,\beta})}_{x_{i-1},x_i}(\dd f_{i}),
\end{equation}
where we put $x_0=x_k=x$.
Now we integrate out over $x_{k,2}^{(i)}, \ldots, x_{k,k}^{(i)}$ for
any $k\in\N$ and $i=1,\ldots,\lambda_k$. In this way, we
obtain that we may replace the bridges $B^{(k,i)}_j$ under the measure
\[
\bigotimes_{k=1}^N \bigotimes_{i=1}^{\lambda_k}\biggl(\int_{\Lambda
}\dd x_k^{(i)} \,\mu_{x_{k}^{(i)}, x_{k}^{(i)}}^{
({\rombc,k\beta})}\biggr)
\]
by the bridges $B^{(k,i)}=B^{(k,i)}_1\diamond\cdots
\diamond B^{(k,i)}_k$ under the measure
\[
\bigotimes_{k=1}^N\biggl[\int_\Lambda\dd x \, g_{k\beta}^{(\rombc
)}(x,x)\biggr]^{\lambda_k}\bigl(\E_\Lambda^{({\rombc,k\beta
})}\bigr)^{\otimes\lambda_k}.
\]
Summarizing, we get
\[
Z_N^{({\rombc})}(\beta,\Lambda)=\sum_{ \lambda\in\mathfrak
{P}_N}\frac{A(\lambda)}{N!} \prod_{k=1}^N\biggl[\int_\Lambda\dd x\,
g_{k\beta}^{(\rombc)}(x,x)
\biggr]^{\lambda_k} \bigotimes_{k\in\N}\bigl(\E_\Lambda^{
({\rombc,k\beta})}\bigr)^{\otimes\lambda_k}
[\ee^{-\mathcal{G}_{N,\beta}^{({\lambda})}}],
\]
where $ A(\lambda)= \#\{\sigma\in\Sym_N\dvtx\sigma$
has $\lambda_k$ cycles of length $k, \forall k\in\N\}$ is the
size of the conjugacy class for the integer partition $\lambda\in
\mathfrak{P}_N $. Standard counting arguments (see \cite
{C02}, Theorem 12.1) give
\[
A(\lambda)=\frac{N!}{\prod_{k=1}^N(\lambda_k!k^{\lambda_k})},
\]
and conclude the proof.
\end{pf}

Now we explain our rewrite of the partition sum in terms of the marked
Poisson point process introduced in Section \ref{sec-notation}, that
is, we prove Proposition \ref{lem-rewrite}. The main idea is to
replace the sum over integer partitions in Lemma \ref{lem-Cycle} by an
expectation with respect to the marked Poisson point process under
conditions on the mark events. We restrict to the case of empty
boundary conditions; the other two require only notational changes.

It will be convenient to write the process $\omega_{\romP}$ as the
superposition
%
%
\begin{equation}
\omega_{\romP}=\sum_{k\in\N}\omega_{\romP}^{(k)}\qquad \mbox
{where }\omega_{\romP}^{(k)}=\sum_{x\in\xi_{\romP}^{(k)}}\delta_{(x,B_x)},
\end{equation}
and $\omega_{\romP}^{(k)}$ is the Poisson process on $\R^d\times
\Ccal_k$ with intensity measure $\nu_k$ defined in (\ref{nudef}).
The processes $\omega_{\romP}^{(k)}$ are independent.
\begin{pf*}{Proof of Proposition \ref{lem-rewrite}}
We start from Lemma \ref{lem-Cycle}. Pick an integer partition
$\lambda\in\mathfrak{P}_N $ with occupation number $ \lambda_k $
satisfying $\sum_{k=1}^N k \lambda_k=N$, and abbreviate the number of
cycles of $\lambda$ by $m=\sum_{k=1}^N \lambda_k$. For any $k\in\N$,
the family $(B^{(k,i)})_{i=1,\ldots, \lambda_k}$ under the measure
$(\PP_{\Lambda}^{({k\beta})})^{\otimes\lambda_k}$ has the same
distribution as the family of marks $(B_x)_{x\in\xi_{\romP}^{(k)}}$
of the conditional Poisson process $\omega_{\romP}^{(k)}$
given $\{\#(\xi_{\romP}^{(k)}\cap\Lambda)=\lambda_k\}$.
Considering the
product measure $\bigotimes_{k\in\N} (\PP_{\Lambda}^{({k\beta
})})^{\otimes
\lambda_k}$ is equivalent to considering the superposition of the
conditional processes $\omega_{\romP}^{(k)}$ with $k\in\N$.

Hence, we have precisely $m$ Poisson points in $\Lambda$.
For any $k\in\N$, conditional on $\{\#(\xi_{\romP}^{(k)}\cap
\Lambda)=\lambda_k\}$, the set $\xi_{\romP}^{(k)}\cap\Lambda$
has the same distribution as the set of starting points, $\{B^{
({k,1})}(0),\ldots,B^{({ k, \lambda_k})}(0)\}$. A comparison of
(\ref{Hamiltonian}) and (\ref{Tdef}) with (\ref{GWdef}) shows that the
interaction term $\Gcal_{N,\beta}^{({\lambda})}$ must be
replaced by the Hamiltonian $H_\Lambda(\omega_{\romP})$. Hence,
\[
\bigotimes_{k\in\N}\bigl(\E_{\Lambda}^{({k\beta})}
\bigr)^{\otimes\lambda_k}[\ee^{-\Gcal_{N,\beta
}^{({\lambda
})}}]
=\ttE\bigl[\ee^{-H_\Lambda(\omega_{\romP})}
| \forall k\in
\N, \#\bigl(\xi_{\romP}^{(k)}\cap\Lambda\bigr)=\lambda_k\bigr].
\]
We see in an elementary way that
%
%
\begin{eqnarray}\label{rewrite3}
&&\ttE\bigl[\ee^{-H_\Lambda(\omega_{\romP})}
|\forall k\in\N
, \#\bigl(\xi_{\romP}^{(k)}\cap\Lambda\bigr)=\lambda_k\bigr]\nonumber\\
&&\qquad=\ttE\bigl[\ee^{-H_\Lambda(\omega_{\romP})}
\mathbh{1}{\bigl\{\forall
k\in\N, \#\bigl(\xi_{\romP}^{(k)}\cap\Lambda\bigr)=\lambda_k\bigr\}}|
\#(\xi_{\romP}\cap\Lambda)=m\bigr]\\
&&\qquad\quad{}\times\frac{\prod_{k\in\N}\lambda_k!}{m!}\overline
q^m\prod_{k\in\N}
(q_k)^{-\lambda_k},\nonumber
\end{eqnarray}
where $\overline q$ and the $q_k$ are defined in (\ref{q*def}). Let us
summarize all the terms involving $\lambda_k$ from (\ref{cycle}) and
(\ref{rewrite3}) [noting that $g_\beta(x,x)=(4\pi\beta k)^{-d/2}$]:
\[
\biggl(\prod_{k\in\N}\frac{(4\pi\beta k)^{-d/2 \lambda
_k}|\Lambda|^{\lambda_k}}{\lambda_k! k^{\lambda_k}}\biggr)\times
\frac{\prod_{k\in\N}\lambda_k!}{m!}\overline q^m \prod_{k\in\N}
(q_k)^{-\lambda_k} = |\Lambda|^m \frac{\overline{q}^{m}}{m!}.
\]
We denote by $\Ncal_{k,\Lambda} (\omega)=\#\{x\in\Lambda\dvtx\ell
(f_x)=k\}$ and $ N_\Lambda(\omega)=\#(\xi\cap\Lambda) $ the number
of particles in $\Lambda$ (whose marks do not have to be contained in
$\Lambda$) with mark length equal to $k$, respectively, with arbitrary
mark length. Then we get
%
%
\begin{eqnarray}
&&Z_N(\beta,\Lambda) \nonumber\\
&&\qquad=\sum_{m=1}^N|\Lambda|^m \frac{\overline
{q}^{m}}{m!}\mathop{\sum_{\lambda\in\mathfrak{P}_N,}}_{\sum_k
\lambda_k=m}\ttE\bigl[\ee^{-H_\Lambda(\omega_{\romP})}
\\
&&\qquad\quad\hspace*{101.5pt}{}\times  \mathbh{1}{\{
\forall k\in\N,
\Ncal_{k,\Lambda}(\omega_{\romP})=\lambda_k\}}|N_\Lambda(\omega_{\romP})=m\bigr].\hspace*{-20pt}\nonumber
\end{eqnarray}
Note that the event $\{N_\Lambda(\omega_{\romP})=m\}$ has probability
$ |\Lambda|^m \frac{\overline{q}^{m}}{m!} \exp\{ - |\Lambda|
\overline{q} \}$. Hence
%
%
\begin{eqnarray}\label{rewrite4}
&&Z_N(\beta,\Lambda) \nonumber\\
&&\qquad=\ee^{|\Lambda| \overline
{q}}\sum
_{m=1}^N\mathop{\sum_{\lambda\in\mathfrak{P}_N,}}_{\sum_k \lambda
_k=m}\ttE\bigl[\ee^{-H_\Lambda(\omega_{\romP})}\mathbh{1}{\{
\forall k\in\N,
\Ncal_{k,\Lambda}(\omega_{\romP})=\lambda_k\}}\\
&&\qquad\quad\hspace*{162.5pt}{}\times
\mathbh{1}\{N_\Lambda(\omega_{\romP})=m\}\bigr].\nonumber
\end{eqnarray}
Note that the events $\{\forall k\in\N, \Ncal_{k,\Lambda}(\omega
_{\romP})=\lambda_k\}\cap\{N_\Lambda(\omega_{\romP})=m\}$ are a
decomposition of the event $\{N_\Lambda^{(\ell)}(\omega_{\romP})=N\}$.
Hence, the assertion in (\ref{rewrite}) follows.
\end{pf*}

\section{\texorpdfstring{Large-deviations arguments: Proof of Theorem \protect\ref{thm-mainres}}%
{Large-deviations arguments: Proof of Theorem 1.2.}}

\label{sec-LDP}

In this section we prove Theorem \ref{thm-mainres} by
applying large-deviations arguments to the representation of the
partition function in Proposition \ref{lem-rewrite}.
In Sections \ref{sec-ldpgeneral}--\ref{sec-lowbound} we carry out the
proof for empty boundary condition, and in Section \ref{otherbc} we
show how to trace the other two boundary conditions back to this case.
In Section \ref{sec-ldpgeneral} we introduce the main object of our
analysis, the stationary empirical field with respect to the marked
Poisson process $\omega_{\romP}$, and we rewrite the partition
function in terms of this field. We also formulate and explain the main
steps of the proof, among which the crucial large-deviations principle
for that field. In Sections \ref{sec-upbound} and \ref{sec-lowbound}
we prove the upper and lower bounds, respectively, for empty boundary condition.

\subsection{The stationary empirical field}\label{sec-ldpgeneral}

Our analysis is based on a large-deviations principle for the
\textit{stationary empirical field}, defined as follows. For any $\xi
\subset\R^d$ and for any centered box $\Lambda\subset\R^d$, let
$\xi_{({\Lambda})}$ be the $\Lambda$-periodic continuation of
\mbox{$\xi\cap\Lambda$}. Analogously, we define the $\Lambda$-periodic
continuation of the restriction of the configuration $\omega$ to
$\Lambda$ as
%
%
\begin{equation}
\omega_{({\Lambda})} = \sum_{z \in\Z^{d}} \sum_{x \in\xi
\cap\Lambda} \delta_{(x+Lz, f_{x})} \qquad\mbox{if } \omega=
\sum_{x \in\xi} \delta_{(x, f_{x})}\in\Omega,
\end{equation}
where $L$ is the side length of the centered cube $\Lambda$.
Then the stationary empirical field is given by
%
%
\begin{equation}\label{RNdef}
\RR_{\Lambda, \omega} = \frac1 {|\Lambda|} \int_{\Lambda} \dd y\,
\delta_{\theta_y ({\omega}_{({\Lambda})})},\qquad \omega\in
\Omega,
\end{equation}
where the shift operator $ \theta_y \dvtx\R^d \to\R^d$ is defined
by $\theta_y(x)=x-y$. It is clear that $\RR_{\Lambda,\omega}$ is a
shift-invariant probability measure on $\Omega$, that is, it is an
element of $\mathcal{P}_{\theta}$.

Now we express $N_\Lambda^{(\ell)}(\omega)$ in terms of $\RR
_{\Lambda,\omega}$. Recall that $U$ denotes the centered unit box.;
we write $\Lambda_L$ for $\Lambda$.
\begin{lemma}\label{lem-totmark} For any centered box $\Lambda\subset
\R^d$ with $|\Lambda|>1$, and any $\omega\in\Omega$,
\[
|\Lambda|\bigl\langle\RR_{\Lambda,\omega}, N_U^{(\ell)}
\bigr\rangle=N_{\Lambda}^{(\ell)}(\omega).
\]
\end{lemma}
\begin{pf}
The assertion follows from \cite{GZ93}, Remark 2.3(1); however, we
give a direct proof without using Palm measures. Let $L>1$ be such that
$\Lambda= \Lambda_L=[-\frac L2,\frac L2]^d$.
We calculate
\begin{eqnarray*}
|\Lambda|\bigl\langle\RR_{\Lambda,\omega}, N_U^{(\ell)}
\bigr\rangle&=&\int_\Lambda\dd z N_U^{(\ell)}\bigl(\theta_z\bigl(\omega
_{({\Lambda})}\bigr)\bigr)
=\sum_{x \in\xi_{({\Lambda})}} \int_\Lambda\dd z\, \mathbh
{1}_{U-x}(z) \ell(f_{x})\\
&=&\mathop{\sum_{x \in\xi_{({\Lambda})}}}_{x \in\Lambda+U}\ell
(f_{x}) |\Lambda\cap(U-x)| \\
&=&N_\Lambda^{(\ell)}(\omega) + \sum_{x \in\xi_{({\Lambda
})} \cap((\Lambda+ U) \setminus\Lambda)}\ell(f_{x}) |\Lambda\cap
(U-x)|\\
&&{} + \sum_{x \in\xi\cap\Lambda}\ell(f_{x}) \bigl(|\Lambda
\cap(U-x)|-1 \bigr).
\end{eqnarray*}
It remains to show that the sum of the two last sums is equal to zero.
Note that the last sum can be restricted to $x \in\xi\cap(\Lambda
\setminus\Lambda_{L-1})$. We use the fact that for each point $x \in
\xi\cap(\Lambda\setminus\Lambda_{L-1})$ there exists a collection
of points in $\xi_{(\Lambda)} \cap( \Lambda_{L+1} \setminus
\Lambda)$, with the same mark of $x$. Indeed, there exists a positive
integer $m(x) \le d$ and a set $\{ x'_1, \ldots, x'_{m(x)}\}$, such
that $x'_i \in\xi_{(\Lambda)} \cap(\Lambda+ U) \setminus
\Lambda$, $x'_i=x+Lz_i$ for some $z_i\in\Z^d$ and $\sum
_{i=1}^{m(x)} |\Lambda\cap(U-x'_i)|= 1 - |\Lambda\cap(U-x)| $.
Notice that
\[
\bigcup_{x \in\xi\cap(\Lambda\setminus\Lambda_{L-1})} \bigcup
_{i=1}^{m(x)} x'_i = \xi_{(\Lambda)} \cap\bigl((\Lambda+
U)\setminus\Lambda\bigr)
\]
and $f_{x} =f_{x'_i }$, for any $ i\le m(x)$.
Hence
\[
\sum_{x \in\xi_{(\Lambda)} \cap((\Lambda+ U) \setminus
\Lambda)}\ell(f_{x}) |\Lambda\cap(U-x)| = \sum_{x \in\xi\cap
\Lambda}\ell(f_{x}) \bigl(1- |\Lambda\cap(U-x)| \bigr).\quad
\]
\upqed\end{pf}

Now we express the interaction Hamiltonian in terms of integrals of the
stationary empirical field against suitable functions; more precisely,
we give lower and upper bounds. In the following lower bound, it is
important that this functional is local and bounded; this will be
achieved up to a small error only.

Fix large truncation parameters $M, R$ and $K$ and introduce $ \xi
^{(\le K)} =\{x\in\xi\dvtx\ell(f_x)\le K\} $ for $\omega
\in\Omega$ and
%
%
\begin{equation}\label{Phitruncdef}
\Phi_\beta^{(R,M,K)}(\omega)=\sum_{x\in\xi^{(\le K)}\cap U}\sum
_{y\in\xi^{({\le K})}\cap
\Lambda_R}T_{x,y}^{(M)}(\omega),
\end{equation}
where $\Lambda_R=[-\frac R2,\frac R2]^d$ and
\[
T_{x,y}^{(M)}(\omega)=\frac{1}{2}\sum_{i=0}^{\ell(f_x)-1}\sum
_{j=0}^{\ell(f_y)-1}\mathbh{1}_{\{(x,i)\not=(y,j)\}}\int_0^\beta
v_M\bigl(|f_x(i\beta+s)-f_y(j\beta+s)|\bigr) \,\dd s,
\]
and where $v_M(r)=(v\wedge M)(r)=\min\{v(r),M\}$. Recall that
$N_\Lambda(\omega)=\#(\xi\cap\Lambda)$ denotes the particle number
in a measurable set $\Lambda\subset\R^d$.
\begin{lemma}[(Hamiltonian bounds)]\label{lem-Hloweresti} Fix any
centred box $\Lambda=\Lambda_L$.
\begin{enumerate}[(ii)]
\item[(i)] For any $M, R, K,S\in(1,\infty)$, and for $L\ge R+2 $,
%
%
\begin{eqnarray}
H_\Lambda(\omega) &\ge&|\Lambda|\bigl\langle\RR_{\Lambda,\omega
},\Phi_\beta^{(R,M,K)}\mathbh{1}\{N_{\Lambda_R}\le S\}
\bigr\rangle\nonumber\\[-8pt]\\[-8pt]
&&{}-CN_{\Lambda_L\setminus\Lambda_{L-R-2}}(\omega),\qquad
\omega\in\Omega,\nonumber
\end{eqnarray}
where $ C=2^d\beta MK^2rS $, and $r$ depends only on $R$ and $d$.

\item[(ii)]
%
%
\begin{equation}
H_\Lambda(\omega)\leq|\Lambda|\langle\RR_{\Lambda,\omega
},\Phi_\beta\rangle,\qquad \omega\in\Omega.
\end{equation}
\end{enumerate}
\end{lemma}
\begin{pf} 
(i) Estimate
%
%
\begin{eqnarray}\label{Hloweresti1}\qquad
&&
|\Lambda|\bigl\langle\RR_{\Lambda,\omega},\Phi_\beta^{
({R,M,K})}\mathbh{1}\{N_{\Lambda_R}\le S\}\bigr\rangle\nonumber\\
&&\qquad=\int_\Lambda\dd z\, \Phi_\beta^{(R,M,K)}\bigl(\theta_z\bigl(\omega
_{(\Lambda)}\bigr)\bigr)\mathbh{1}\bigl\{N_{\Lambda_R}\bigl(\theta_z\bigl(\omega
_{(\Lambda)}\bigr)\bigr)\le S\bigr\}\nonumber\\
&&\qquad \le\int_\Lambda\dd z \sum_{x\in\xi^{(\le K)}_{(\Lambda)}\cap
(U-z)}\sum_{y\in\xi^{(\le K)}_{(\Lambda)}\cap(\Lambda_R-z)}
T_{x,y}^{(M)}\bigl(\omega
_{(\Lambda)}\bigr)\nonumber\\[-8pt]\\[-8pt]
&&\qquad\quad\hspace*{151.5pt}{}\times\mathbh{1}\bigl\{\#\bigl(\xi_{({\Lambda
})}^{(\le K)}\cap(\Lambda_{R}-z)\bigr)\le S\bigr\}\nonumber\\
&&\qquad=\mathop{\sum_{x,y\in\xi^{(\le K)}_{(\Lambda)},
x\in\Lambda+U,}}_{y\in\Lambda+\Lambda_R, x\in\Lambda_{R+1}+y}
T_{x,y}^{(M)}\bigl(\omega_{(\Lambda)}\bigr)\int_{\Lambda\cap
(U-x)\cap(\Lambda_R-y)}\dd z\nonumber\\
&&\qquad\quad\hspace*{78.7pt}{}\times
\mathbh{1}\bigl\{\#\bigl(\xi^{(\le K)}_{(\Lambda
)}\cap(\Lambda_R-z)\bigr)\le
S\bigr\}.\nonumber
\end{eqnarray}
Observe that the integral over $z$ is not larger than one. Now we split
the last sum into the sums on $(x,y)\in\Lambda^2$ and the remainder.
For $(x,y)\in\Lambda^2$, we may replace $T_{x,y}^{(M)}(\omega
_{(\Lambda)})$ by $T_{x,y}^{(M)}(\omega)$ and estimate it
against $T_{x,y}(\omega)$. Hence,
\[
\mbox{left-hand side of (\ref{Hloweresti1})}\leq H_\Lambda(\omega)+\Psi
_\Lambda^{({R,M,K,S})}(\omega),
\]
where the remainder term is
\begin{eqnarray*}
&&\Psi_\Lambda^{(R,M,K,S)}(\omega)\\
&&\qquad=\mathop{\sum_{x,y\in\xi^{(\le K)}_\Lambda, x\in
\Lambda+U,}}_{y\in\Lambda+\Lambda_R, x\in\Lambda
_{R+1}+y,(x,y)\notin\Lambda^2}T_{x,y}^{(M)}\bigl(\omega_{
({\Lambda})}\bigr)\int_{\Lambda\cap(U-x)\cap(\Lambda_R-y)}\dd z \\
&&\qquad\quad\hspace*{115.8pt}{}\times \mathbh{1}\bigl\{\#\bigl(\xi^{(\le K)}_{
({\Lambda})}\cap(\Lambda_R-z)\bigr)\le S\bigr\}\\
&&\qquad\leq\frac12\beta MK^2\\
&&\qquad\quad{}\times\mathop{\sum_{x,y\in\xi^{(\le K)}_\Lambda,
x\in\Lambda+U,}}_{y\in\Lambda+\Lambda_R, x\in
\Lambda_{R+1}+y,(x,y)\notin\Lambda^2}\mathbh{1}\bigl\{\exists z\in\Lambda \cap(U-x)\\
&&\qquad\quad\hspace*{140pt}{}\cap(\Lambda
_R-y)\dvtx\#\bigl(\xi_\Lambda^{(\le K)}\cap(\Lambda_R-z)\bigr)\le
S\bigr\}\\
&&\qquad\le\frac12\beta MK^2\mathop{\sum_{x,y\in\xi^{(\le K)}_{(\Lambda
)}, x\in\Lambda+U,}}_{y\in\Lambda+\Lambda_R,
x\in\Lambda_{R+1}+y,(x,y)\notin\Lambda^2}\mathbh{1}\bigl\{\#\bigl(\xi
_{(\Lambda)}^{(\le K)}\cap(\Lambda_{R-1}+x)\bigr)\le S\bigr\}.
\end{eqnarray*}
The sum over $(x,y)\notin\Lambda^2 $ is split into the sum over $
x\in(\Lambda+U)\setminus\Lambda, y\in\Lambda+\Lambda_R $ and $
x\in\Lambda+U, y\in(\Lambda+\Lambda_R)\setminus\Lambda$. Recall
that $ \Lambda=\Lambda_L $ and that $ L\ge R+1 $. The condition $
x\in\Lambda_{R+1}+y $ implies that in both cases $y $ is summed over
a subset of $ \Lambda_{L+R+2}\setminus\Lambda_{L-R-1} $. Hence,
\begin{eqnarray*}
&&\Psi_\Lambda^{(R,M,K,S)}(\omega)\\
&&\qquad\le\frac12 \beta MK^2\\
&&\qquad\quad{}\times\sum
_{y\in\xi^{(\le K)}_{(\Lambda)}\cap(\Lambda
_{L+R+2}\setminus\Lambda_{L-R-1})}\#\bigl\{x\in\xi_{(\Lambda)}^{(\le K)}\cap
(\Lambda_{R+1}+y)\dvtx\\
&&\qquad\quad\hspace*{131.8pt} \#\bigl(\xi_{(\Lambda)}^{(\le K)}\cap
(\Lambda_{R-1}+x)\bigr)\le S\bigr\}.
\end{eqnarray*}
Now we show that the counting factor is not larger than $ rS $, where $
r $ depends only on $ R $ and the dimension $ d $. Indeed, cover $
\Lambda_{R+1} +y$ with $r $ boxes $\Delta_1,\ldots,\Delta_r $ of
diameter $ (R-1)/2 $, then
\begin{eqnarray*}
&&\#\bigl\{x\in\xi_\Lambda^{(\le K)}\cap(\Lambda
_{R+1}+y)\dvtx\#\bigl(\xi_{(\Lambda)}^{(\le K)}\cap
(\Lambda_{R-1}+x)\bigr)\le S\bigr\}\\
&&\qquad\le\sum_{i=1}^r \#\bigl\{x\in\xi_{(\Lambda)}^{(\le K)}\cap
\Delta_i\dvtx\#\bigl(\xi_{(\Lambda)}^{(\le K)}\cap
(\Lambda_{R-1}+x)\bigr)\le S\bigr)\\
&&\qquad\le\sum_{i=1}^r\#\bigl\{x\in\xi_{(\Lambda)}^{(\le K)}\cap\Delta
_i\dvtx\#\bigl(\xi_{(\Lambda)}^{(\le K)}\cap\Delta_i\bigr)\le S\bigr\}\\
&&\qquad\le rS,
\end{eqnarray*}
since $ \Delta_i\subset\Lambda_{R-1}+x $ if $ x\in\Delta_i $.
This gives
\begin{eqnarray*}
\Psi^{(R,M,K,S)}_\Lambda(\omega)&\le&\tfrac12\beta
MK^2rSN_{\Lambda_{L+R+2}\setminus\Lambda_{L-R-1}}\bigl(\omega_{({
\Lambda})}\bigr)\\
&\le&2^d\beta MK^2rSN_{\Lambda_L\setminus\Lambda_{L-R-2}}(\omega),
\end{eqnarray*}
and finishes the proof of (i).

(ii) 
In a similar way as in (\ref
{Hloweresti1}), one sees that, for any $\omega\in\Omega$,
%
%
\begin{eqnarray}\label{Hloweresti2}\qquad
|\Lambda|\langle\RR_{\Lambda,\omega},\Phi_\beta\rangle
&=&\sum_{x,y\in\xi_{(\Lambda)}} T_{x,y}\bigl(\omega_{({\Lambda
})}\bigr) |\Lambda\cap(U-x)|\nonumber\\
&=&H_\Lambda(\omega)
+\sum_{x,y\in\xi\cap\Lambda} T_{x,y}\bigl(\omega_{(\Lambda)}\bigr)
\bigl(|\Lambda\cap(U-x)|-1\bigr)\\
&&{} +\sum_{x,y\in\xi_{(\Lambda)}\dvtx x\in\Lambda_{L+1},
(x,y)\notin\Lambda^2} T_{x,y}\bigl(\omega_{(\Lambda)}\bigr) |\Lambda
\cap
(U-x)|.\nonumber
\end{eqnarray}
It remains to show that the sum of the two last sums is nonnegative.
Note that the sum on $x$ in the first sum may be restricted to $x\in
\xi\cap(\Lambda\setminus\Lambda_{L-1})$. For each such $x$ and for
any $y\in\xi\cap\Lambda$, there exist a positive integer $m(x) \le
d$ and a set $\{x'_1, y'_1, \ldots, x'_{m(x)}, y'_{m(x)}\}$, such that
$ x'_i \in\xi_{(\Lambda)}\cap\Lambda_{L+1}$, $x'_i=x+Lz_i$ and
$y'_i=y+Lz_i$ for some $z_i\in\Z^d$, and
\[
\sum_{i=1}^{m(x)}|\Lambda\cap(U-x'_i)|=|\Lambda\cap(U-x)|-1.
\]
Then $T_{x,y}(\omega_{(\Lambda)})=T_{x',y'}(\omega_{({
\Lambda})})$ by $\Lambda$-periodicity of $\omega_{(\Lambda)}$.
This shows that the sum of the two last sums in (\ref{Hloweresti2}) is
nonnegative, which finishes the proof of (ii).
\end{pf}

Recall that $ L_{N} = (N/\rho)^{d}$. Applying Lemmas \ref
{lem-totmark} and \ref{lem-Hloweresti}(i) to the representation in
Proposition \ref{lem-rewrite}, we obtain, for any $R,M,K,S>0$, the
upper bound
%
%
\begin{eqnarray}\label{LDPmainstatement}\quad
&&Z_N(\beta,\Lambda_{L_N})\nonumber\\
&&\qquad\leq\ee^{|\Lambda
_{L_N}|\overline q} \ttE\bigl[\exp\bigl\{-|\Lambda_{L_N}|
\bigl\langle\RR_{\Lambda_{L_N},\omega_{\romP}},\Phi_\beta
^{({R,M,K})}\mathbh{1}\{N_{\Lambda_R}\le S\}\bigr\rangle\bigr\}\\
&&\qquad\quad\hspace*{40.24pt}{} \times\exp\{C N_{\Lambda_{L_N}\setminus\Lambda
_{L_N-R-2}}(\omega_{\romP})\}\mathbh{1}\bigl\{\bigl\langle\RR_{\Lambda
_{L_N},\omega_{\romP}},N_U^{(\ell)}\bigr\rangle=\rho\bigr\}\bigr],\nonumber
\end{eqnarray}
for any $N\in\N$, and, using Lemmas \ref{lem-totmark} and \ref
{lem-Hloweresti}(ii), the lower bound
%
%
\begin{equation}\label{LDPmainstatementlower}\qquad
Z_N(\beta,\Lambda_{L_N})\geq\ee^{|\Lambda
_{L_N}|\overline q} \ttE\bigl[\ee^{-|\Lambda_{L_N}|\langle\RR_{\Lambda
_{L_N},\omega_{\romP}},\Phi_\beta\rangle}\mathbh{1}\bigl\{\bigl\langle\RR_{\Lambda
_{L_N},\omega_{\romP}},N_U^{(\ell)}\bigr\rangle=\rho\bigr\}\bigr],
\end{equation}
for any $N\in\N$.

The main point of introducing the stationary empirical field is that
the family $(\RR_{\Lambda_L, \omega_{\romP}})_{L>0}$ satisfies a
large-deviations principle on $\Pcal_\theta$, which is known from the
work by Georgii and Zessin. On $\Pcal_{\theta}$ we consider the
following topology. A~measurable function $g\dvtx\Omega\to\R$ is
called \textit{local} if it depends only on the restriction of $\omega$
to some bounded open cube, and it is called \textit{tame} if $|g| \le c(1
+N_\Lambda)$ for some bounded open cube $\Lambda$ and some constant
$c \in\R^+$. We endow the space $\Pcal_{\theta}$ with the topology
$\tau_{\Lcal}$ of \textit{local convergence}, defined as the smallest
topology on $\Pcal_{\theta}$ such that the mappings $P \mapsto
\langle P,g\rangle$ are continuous for any $g\in\Lcal$, where $\Lcal
$ denotes the linear space of all local tame functions. It is clear
that the map $P\mapsto\langle P,N_U\rangle$ is $\tau_{\Lcal
}$-continuous; however, the map $P\mapsto\langle P,N_U^{({\ell
})}\rangle$ is only lower semicontinuous.
\begin{lemma}[(Large deviations for $\RR_{\Lambda_L, \omega_{\romP
}}$)]\label{ldpyn}
The measures $\RR_{\Lambda_L,\omega_{\romP}}$ satisfy, as $L\to\infty$,
a large-deviations principle in the topology $\tau_\Lcal$ with speed
$|\Lambda_{L}|$ and rate function
$I_\beta\dvtx\Pcal_\theta\to[0,\infty]$ defined in (\ref{Idef}). The
function $I_\beta$ is affine and lower $\tau_\Lcal$-semicontinuous and
has $\tau_\Lcal$-compact level sets.
\end{lemma}
\begin{pf}
This is \cite{GZ93}, Theorem 3.1.
\end{pf}

Our goal is to apply Varadhan's lemma to the expectations on the
right-hand sides of (\ref{LDPmainstatement}) and (\ref
{LDPmainstatementlower}). In conjunction with the large-deviations
principle of Lemma \ref{ldpyn}, this formally suggests that both
(\ref{freeenupper}) and (\ref{freeenlower}) should be valid, as we
explain now. Indeed,
first consider (\ref{LDPmainstatementlower}) and note that the map
$P\mapsto\langle P,\Phi_\beta\rangle$ has the proper continuity
property for the application of the
lower bound half of Varadhan's lemma. If one neglects the fact that the
condition
$\langle P, N_U^{(\ell)}\rangle=\rho$ does not
define an open set of $P$'s, then one easily formally obtains (\ref
{freeenlower})
from (\ref{LDPmainstatementlower}).

Now we consider (\ref{LDPmainstatement}). Assume that the term
$N_{\Lambda_{L_N}\setminus\Lambda_{L_N-R-2}}(\omega_{\romP})$ is a
negligible error term and that taking\vspace*{1pt} the truncation parameters $R,M,K$
and $S$ to infinity will finally turn $\Phi_\beta^{(R,M,K)}\mathbh
{1}\{N_{\Lambda_R}\le S\}$ into $\Phi_\beta
$. The functional $P\mapsto\langle P,\Phi_\beta^{(R,M,K)}\mathbh
{1}\{N_{\Lambda_R}\le S\}\rangle$ has the sufficient
continuity property for the application of the upper bound half of
Varadhan's lemma.
However, the functional $P\mapsto\langle P,N_U^{(\ell)}\rangle$
is not upper semicontinuous. Hence, the \textit{equality} $\langle\RR
_{\Lambda_{L_N},\omega_{\romP}}, N_U^{(\ell)}\rangle=\rho$ is
turned into the \textit{inequality} $\langle P, N_U^{(\ell)}\rangle
\leq\rho$ in the resulting variational formula. Therefore, one easily
formally obtains (\ref{freeenupper})
from (\ref{LDPmainstatement}). In particular, our upper and lower
bounds in Theorem \ref{thm-mainres} may differ. For small $\beta$,
respectively, small $\rho$, we improve the proof in Lemma \ref
{lem-Zquotient} and achieve a coincidence of upper and lower bounds,
but this has nothing to do with large-deviations arguments.

The lack of upper semicontinuity of the functional $P\mapsto\langle
P,N_U^{(\ell)}\rangle$ causes serious technical problems in the
proof of the lower bound, since the condition $\langle P, N_U^{(\ell
)}\rangle=\rho$ must be approximated by some open condition.

In Lemma \ref{lem-Hloweresti}, we already estimated away all the
interaction involving cycles of length $>K$, and in the proof of the
lower bound we will restrict the configuration space to marks with
lengths $\leq K$. This is why our variational formulas spot only the
presence of ``finitely long'' cycles.

\subsection{The upper bound for empty boundary condition}\label{sec-upbound}

In this section, we prove the upper bound in (\ref
{freeenupper}) for $\rombc=\varnothing$. According to (\ref
{LDPmainstatement}), it will be sufficient to prove
%
%
\begin{eqnarray}\label{uppempty}\quad
&&\limsup_{R,M,K,S\to\infty}\limsup_{N\to\infty}\frac1{|\Lambda
_{L_N}|}\log\ttE\bigl[
\exp\bigl\{-|\Lambda_{L_N}|\nonumber\\
&&\qquad\quad\hspace*{191pt}\hspace*{-53pt}{}\times\bigl\langle\RR_{\Lambda
_{L_N},\omega_{\romP}},\Phi_\beta^{(R,M,K)}\mathbh{1}\{
N_{\Lambda_R}\le S\}\bigr\rangle\bigr\}\nonumber\\
&&\qquad\quad\hspace*{108.4pt}{} \times\exp\{C N_{\Lambda_{L_N}\setminus\Lambda
_{L_N-R-2}}(\omega_{\romP})\}\\
&&\qquad\quad\hspace*{108.4pt}\hspace*{39.3pt}{} \times\mathbh{1}\bigl\{\bigl\langle\RR_{\Lambda
_{L_N},\omega_{\romP}},N_U^{(\ell)}\bigr\rangle=\rho\bigr\}
\bigr]\nonumber\\
&&\qquad\leq-\chi^{(\leq)}(\beta,\rho).\nonumber
\end{eqnarray}

An outline of the proof is as follows. We separate first the two
exponential terms from each other with the help of H\"older's
inequality. The latter term will turn out to be a negligible error
term. The functional that appears in the first exponent turns out to be
local and bounded. Since its integral against a probability measure $P$
is a $\tau_{\Lcal}$-continuous and bounded function of $P$,
Varadhan's lemma can be applied and expresses the limit superior in
terms of the variational formula for the truncated versions of the
interaction functionals. The indicator on the event $\{\langle\RR
_{\Lambda_{L_N}, \omega_{\romP}}, N_{U}^{(\ell)} \rangle= \rho
\}$ is estimated against the indicator on its closure, which is the
same set with ``$\leq$'' instead of ``$=$.'' In this way, we obtain an
upper bound against a truncated version of the variational formula
$-\chi^{(\leq)}(\beta,\rho)$. By letting the truncation
parameters go to infinity, this formula converges to $-\chi^{(\leq
)}(\beta,\rho)$.

Let us turn to the details. We abbreviate $\RR_N=\RR_{\Lambda_{L_N},
\omega_{\romP}}$.\vspace*{1pt}

We pick $ \eta\in(0,1) $ and start from (\ref{LDPmainstatement}),
then H\"older's inequality gives
%
%
\begin{eqnarray}\label{hoelder}\quad
Z_N(\beta,\Lambda_{L_N})&\le&\ee^{|\Lambda
_{L_N}|\overline q} \ttE \bigl[ \ee^{-{1}/({1-\eta})|\Lambda
_{L_N}|\langle\RR_{N},\Phi
_\beta^{(R,M,K)}\mathbh{1}\{N_{\Lambda_R}\le S\}\rangle}
\nonumber\\
&&\qquad\quad\hspace*{74.7pt}{}\times \mathbh{1}\bigl\{\bigl\langle\RR_{N}, N_{U}^{(\ell)} \bigr\rangle\leq
\rho\bigr\}\bigr]^{1-\eta}\\
&&{} \times\ttE\bigl[\ee^{{1/\eta}CN_{\Lambda
_{L_N}\setminus\Lambda_{L_N-R-2}}(\omega_{\romP})}\bigr]^\eta;\nonumber
\end{eqnarray}
note that we also estimated ``$=\rho$'' against ``\mbox{$\le$}$\rho$'' in the indicator.
The second term on the right-hand side of (\ref{hoelder}) is easily
estimated using the fact that $N_{\Lambda_{L_N}\setminus\Lambda
_{L_N-R-2}}$ is a Poisson random variable with parameter $\overline
q\times|\Lambda_{L_N}\setminus\Lambda_{L_N-R-2}| $ and that this
parameter is of surface order $ L_N^{d-1}=o(|\Lambda_N|) $. Hence, the
expectation is estimated
\begin{eqnarray*}
&&\ttE\bigl[\ee^{{1/\eta}C N_{\Lambda_{L_N}\setminus
\Lambda_{L_N-R-2}}(\omega_{\romP})}\bigr]^\eta\\
&&\qquad= \ee^{-\eta\overline{q}|\Lambda_{L_N}\setminus\Lambda_{L_N-R-2}|}
\exp\{\eta\ee^{C/\eta}\overline{q}|\Lambda_{L_N}\setminus
\Lambda_{L_N-R-2}|\}\\
&&\qquad\le\ee^{o(|\Lambda_{L_N}|)}.
\end{eqnarray*}

We turn to the first term on the right-hand side of (\ref{hoelder}).
It turns out that $\Phi^{(R,M,K)}_\beta\mathbh{1}\{N_{\Lambda
_R}\le S \}$ is bounded. In fact,
%
%
\begin{eqnarray}\label{boundedtame}
&&\Phi^{(R,M,K)}_\beta(\omega) \mathbh{1}{\{N_{\Lambda
_R}(\omega)\le S \}}\nonumber\\
&&\qquad\leq\frac12 M\beta\biggl[\sum_{x \in U \cap\xi} \ell(f_x)\sum
_{y \in\Lambda_R \cap\xi} \ell(f_y)+\biggl(\sum_{x \in U \cap\xi
} \ell(f_x)\biggr)^2\biggr]\nonumber\\[-8pt]\\[-8pt]
&&\qquad\quad{}\times\mathbh{1}{\{N_{\Lambda_R}(\omega) \le S \}
}\nonumber\\
&&\qquad\leq M\beta K^2S^2 .\nonumber
\end{eqnarray}
Furthermore, it is easily seen that it is also local. Therefore, the map
\[
P\mapsto\bigl\langle P, \Phi^{(R,M,K)}_\beta\mathbh{1}_{\{
N_{\Lambda_R} \le S \}} \bigr\rangle
\]
is bounded and continuous on $\Pcal_\theta$ with respect to the
topology $\tau_\Lcal$. Now we can apply a variant of Varadhan's lemma
\cite{DZ98}, Theorem 4.3.1, in conjunction with the large-deviations
principle of Lemma \ref{ldpyn}, to obtain that
%
%
\begin{eqnarray}\label{ubtrunc}\qquad
&&\limsup_{N\to\infty}\frac1{|\Lambda_{L_N}|}\log\ttE \biggl[
\exp\biggl\{ -\frac1{1-\eta} |\Lambda_{L_N}|
\bigl\langle\RR_{N},\Phi^{(R,M,K)}_\beta\mathbh{1}\{
N_{\Lambda_R} \le S \}\bigr\rangle\biggr\} \nonumber\\
&&\quad\hspace*{201pt}{}\times \mathbh{1}\bigl\{\bigl\langle\RR_{N}, N_{U}^{(\ell)}
\bigr\rangle\le\rho\bigr\}\biggr]\\
&&\qquad\leq-\inf_{P\in\Pcal_\theta\dvtx\langle P, N_{U}^{(\ell)}
\rangle\le\rho}
\biggl( I_\beta(P)+ \frac1{1-\eta}\bigl\langle P,\Phi^{(R,M,K)}_\beta
\mathbh{1}\{ N_{\Lambda_R} \le S \}
\bigr\rangle\biggr),\nonumber
\end{eqnarray}
since the set $\{ P\in\Pcal_\theta\dvtx\langle P, N_{U}^{(\ell)}
\rangle\le\rho\}$ is closed.

It remains to prove that
%
%
\begin{eqnarray}\label{limvar}
&&\liminf_{R,M,K\to\infty,\eta\downarrow0}\liminf_{S\to\infty}
\inf_{P \dvtx\langle P, N_{U}^{(\ell)} \rangle\le\rho} \bigl(
I_\beta(P) + F_{M,R,K,S,\eta}(P)\bigr)\nonumber\\[-8pt]\\[-8pt]
&&\qquad\ge\inf_{ P \dvtx\langle P, N_{U}^{(\ell)} \rangle\le\rho}
\bigl( I_\beta(P) + F(P)
\bigr),\nonumber
\end{eqnarray}
where we used the abbreviations $F(P)=\langle P,\Phi_\beta\rangle$
and $F_{M,R,K,S,\eta}(P)=\frac1{1-\eta}\langle P$, $\Phi^{(R,M,K)}_\beta\mathbh{1}
{\{ N_{\Lambda_R} \le S \}} \rangle$.
Fix $M, R,K >0$ and $\eta\in(0,1)$ and pick a sequence $S_n \to
\infty$ and some $ Q_n $ satisfying $ \langle Q_{n}, N_{U}^{(\ell)}
\rangle\le\rho$ such that
%
%
\begin{eqnarray}\label{discard}
&&
I_\beta( Q_n)+ F_{M,R,K,S_n,\eta}( Q_n)\nonumber\\[-8pt]\\[-8pt]
&&\qquad < \inf_{P \dvtx\langle P,
N_{U}^{(\ell)} \rangle\le\rho} \bigl( I_\beta(P) +
F_{M,R,K,S_n,\eta}(P)\bigr) + \frac1n.\nonumber
\end{eqnarray}
By compactness of the level sets of $I_\beta$, we may assume that the
limiting measure $Q= \lim_{n \to\infty} Q_n$ exists in $\Pcal
_\theta$, where the limit is taken along some suitable subsequence.
Notice further that $ \langle Q ,N_U^{(\ell)} \rangle\leq\rho$
by Fatou's lemma.
Fix any large $ S >0$, then for $n$ sufficiently large,
%
%
\begin{eqnarray}\label{limvar2}
&&
\inf_{ P \dvtx\langle P, N_{U}^{(\ell)} \rangle\le\rho}
\bigl( I_\beta(P) + F_{M,R,K,S_n,\eta}(P) \bigr)\nonumber\\
&&\qquad > I_\beta(Q_n) +F_{M,R,K,S_n,\eta} (Q_n) - \frac1n\\
&&\qquad\ge I_\beta(Q_n) +F_{M,R,K,S,\eta}(Q_n) - \frac1n,\nonumber
\end{eqnarray}
where the second inequality uses the monotonicity of $F_{M,R,K,S,\eta
}$ in $S$.
Now send $n \to\infty$ and use the lower semi-continuity of $I_\beta
$ and the continuity of $F_{M,R,K,S,\eta}$, to get that the limit
inferior of the right-hand side of (\ref{limvar2}) is larger or equal
to $ I_\beta(Q) +F_{M,R,K,S,\eta}(Q)$. Sending $ S\to\infty$ and
using the monotone convergence theorem, we arrive at
%
%
\begin{eqnarray}\label{Slimit}
&&\liminf_{S \to\infty} \inf_{P \dvtx\langle P, N_{U}^{(\ell)}
\rangle\le\rho}\bigl(I_\beta(P) +F_{M,R,K,S,\eta}(P)\bigr)\nonumber\\[-8pt]\\[-8pt]
&&\qquad\ge\inf_{P \dvtx\langle P, N_{U}^{(\ell)} \rangle\le\rho
}\bigl(I_\beta(P) +F_{M,R,K,\infty,\eta}(P)\bigr).
\nonumber
\end{eqnarray}
In a similar way one proves that
\begin{eqnarray*}
&&\liminf_{R,M,K\to\infty,\eta\downarrow0} \inf_{P \dvtx\langle
P, N_{U}^{(\ell)} \rangle\le\rho}\bigl(I_\beta(P) +F_{M,R,K,\infty
,\eta}(P)\bigr)\\
&&\qquad\ge\inf_{ P \dvtx\langle P, N_{U}^{(\ell)} \rangle\le\rho}
\bigl( I_\beta(P) + F(P) \bigr),
\end{eqnarray*}
which implies (\ref{limvar}) and ends the proof of (\ref{uppempty}).

\subsection{The lower bound for empty boundary condition}\label{sec-lowbound}

In this section, we prove the lower bound in (\ref
{freeenlower}) for $\rombc=\varnothing$. According to (\ref
{LDPmainstatementlower}), it will be sufficient to prove
%
%
\begin{eqnarray}
&&\liminf_{N\to\infty}\frac1{|\Lambda_{L_N}|}\log\ttE
\bigl[\ee
^{-|\Lambda_{L_N}|\langle
\RR_{\Lambda_{L_N},\omega_{\romP}},\Phi_\beta\rangle}\mathbh{1}\bigl\{
\bigl\langle\RR_{\Lambda_{L_N},\omega_{\romP}},N_U^{(\ell)}\bigr\rangle
=\rho\bigr\}\bigr]\nonumber\\[-8pt]\\[-8pt]
&&\qquad\geq-\chi^{(=)}(\beta,\rho).\nonumber
\end{eqnarray}

We follow the standard strategy of changing the measure so that
untypical events become typical, and controlling the Radon--Nikodym
density by means of McMillan's theorem. However, for our problem we
have to overcome two major difficulties. First, the map $ P\mapsto
\langle P,\Phi_\beta\rangle$ is not upper semicontinuous, and
second, the set $ \{P\in\Pcal_\theta\dvtx\langle P,N_U^{(\ell
)}\rangle=\rho\} $ appearing in the indicator is not open. This set
induces long-range correlations not only between the points of the
process, but also between their marks. Therefore, the results of \cite
{GZ93} cannot be applied directly, but some ideas of \cite{G94} can be adapted.

We now describe our strategy. In Lemma \ref{lem-makeopen}, we
replace the condition $\langle P,N_U^{(\ell)}\rangle=\rho$ by
the condition $|\langle P,N_U^{(\ell)}\rangle-\rho|<\delta$ for
some small $\delta$ and control the replacement error. This condition
becomes an open condition when restricting the mark space $ E $ to a
cut-off version. A restriction of $\Pcal_\theta$ in Lemma \ref
{lem-Fcont} makes~the map $ P\mapsto\langle P,\Phi_\beta\rangle$
continuous. In order to apply McMillan's theorem to the transformed
point process, an ergodic approximation is carried out in Lemma \ref
{lem-ergodappr}.

Let us turn to the details. First, we prepare for relaxing the
condition ``\mbox{$=$}$\rho$'' to ``\mbox{$\approx$}$\rho$'' in the following step,
which is of independent interest. Bounding the quotient $Z_{N+1}/Z_N$
of partition functions is often the key step to prove the equivalence
of the canonical ensemble with the grand canonical ensemble, where the
particle number is not fixed but governed by the mean. In the
following, we give a lower bound in our case, which will also imply a
nontrivial upper bound for the limiting free energy. Our proof is
carried out in the setting of the cycle expansion introduced in
Section \ref{sec-ProofProp} and is independent of the reformulation in
terms of the marked Poisson point process.
\begin{lemma}\label{lem-Zquotient}
For any $N\in\N$ and any measurable set $\Lambda\subset\R^d$,
%
%
\begin{equation}\label{Zquotient}
\frac{Z_{N+1}(\beta,\Lambda)}{Z_{N}(\beta,\Lambda)}\geq(4\pi
\beta)^{-d/2}\frac{|\Lambda|}{N+1} \ee^{-N\beta
\alpha
(v)/|\Lambda|},
\end{equation}
where we recall that $\alpha(v)=\int_{\R^d} v(|x|) \,\dd x$.
\end{lemma}
\begin{pf}
The strategy is as follows. We start with the cycle expression for the
partition function $ Z_l $. We then add a particle, that is, an
additional cycle of length one, and control the changes in the
combinatorial factor and in the energy. Here our assumption $
\int_{\R^d} v(|x|) \,\dd x <\infty$ allows to bound the additional
interaction energy.

We abbreviate $Z_N(\beta,\Lambda)$ by $Z_N$ in this proof. Recall
(\ref{PNdef}). According to Lem\-ma~\ref{lem-Cycle}, the cycle
representation of the partition function reads
%
%
\begin{equation}
\label{cycleexpression}
Z_N=\sum_{\lambda\in\mathfrak{P}_N} F_1(\lambda)F_2(\lambda),
\end{equation}
with the combinatorial and interaction part
\begin{eqnarray*}
F_1(\lambda)&=&\prod_{k=1}^N\frac{(4\pi\beta k)^{-d \lambda
_k/2}|\Lambda|^{\lambda_k}}{\lambda_k!k^{\lambda_k}},\\
F_2(\lambda)&=&\Biggl(\bigotimes_{k=1}^N\bigl(\E_{\Lambda}^{(k\beta)}
\bigr)^{\otimes\lambda_k}\Biggr)[\ee^{-\Gcal_{N,\beta
}^{(\lambda)}}].
\end{eqnarray*}
Define the injection
\[
T\dvtx\mathfrak{P}_N\to\mathfrak{P}_{N+1},\qquad
T(\lambda)=\widetilde\lambda\qquad\mbox{with } \widetilde
\lambda_k =
\cases{
\lambda_1+1, &\quad if $k=1$,\cr
\lambda_k, &\quad if $k\ge2$.
}
\]
All the terms in (\ref{cycleexpression}) are nonnegative, and hence we
may estimate
%
%
\begin{eqnarray}\label{lowerest}
Z_{N+1}&\ge&\sum_{\widetilde\lambda\in\mathfrak{P}_{N+1}\dvtx
\widetilde\lambda_1\ge1} F_1(\widetilde\lambda)F_2(\widetilde
\lambda)\nonumber\\
&=& \sum_{\lambda\in\mathfrak{P}_N} F_1(T(\lambda))F_2(T(\lambda
))\\
&=&\sum_{\lambda\in\mathfrak{P}_N} \frac{F_1(T(\lambda
))}{F_1(\lambda)}\frac{F_2(T(\lambda))}{F_2(\lambda)} F_1(\lambda
)F_2(\lambda).\nonumber
\end{eqnarray}
The first quotient on the right-hand side of (\ref{lowerest}) is
bounded from below as follows:
%
%
\begin{equation}\label{F1esti}
\frac{F_1(T(\lambda))}{F_1(\lambda)}=(4\pi\beta)^{-d/2}\frac
{|\Lambda|}{\lambda_1+1}\geq(4\pi\beta)^{-d/2}\frac{|\Lambda|}{N+1}.
\end{equation}
The second quotient is estimated via Jensen's inequality as follows.
Recall that $ B_{(j-1)\beta+s}^{(k,i)} $ is the Brownian bridge
of the $j$th leg of the $i$th cycle of length $k$, $1\le i\le\lambda
_k $,\vspace*{-1pt}
%
%
\begin{eqnarray}\label{energyest}
F_2(T(\lambda))&=&\E_{\Lambda}^{(\beta)}\otimes
\Biggl(\bigotimes_{k=1}^N\bigl(\E_{\Lambda}^{(k\beta)}
\bigr)^{\otimes\lambda_k}\Biggr)\nonumber\\
&&{}\times\Biggl[\ee^{-\Gcal_{N,\beta
}^{(\lambda)}}
\exp\Biggl\{-\sum_{k\in\N}\sum_{i=1}^{\lambda
_k}\sum_{j=1}^k\int_0^\beta v\bigl(\bigl|B_s-B_{(j-1)\beta+s}^{(k,i)}\bigr|\bigr) \,\dd
s\Biggr\}\Biggr]\nonumber\\[-8pt]\\[-8pt]
&\ge&\Biggl(\bigotimes_{k=1}^N\bigl(\E_{\Lambda}^{({k\beta
})}\bigr)^{\otimes\lambda_k}\Biggr)\nonumber\\
&&{}\times \Biggl[\ee^{-\Gcal
_{N,\beta}^{(\lambda)}}\exp\Biggl\{-\sum_{k\in\N}\sum_{i=1}^{\lambda_k}
\sum_{j=1}^k\int_0^\beta\E_{\Lambda}^{(\beta)}
\bigl[v\bigl(\bigl|B_s-B_{(j-1)\beta+s}^{(k,i)}\bigr|\bigr)\bigr] \,\dd
s\Biggr\}\Biggr].\nonumber
\end{eqnarray}
Given $ \lambda\in\mathfrak{P}_N$ and $ k\in\N, i\in\{1,\ldots,
\lambda_k\}, j\in\{1,\ldots, k\} $, we write $ f(s):=B_{(j-1)\beta
+s}^{(k,i)}$, and we estimate the expectation in the exponent as
follows:
%
%
\begin{eqnarray}\label{energyesti}\qquad
&&\E_{\Lambda}^{(\beta)}\bigl(v\bigl(|B_s-f(s)|\bigr)\bigr)\nonumber\\
&&\qquad=\frac{1}{|\Lambda|}\int_{\Lambda} \dd x\int_{\Lambda} \dd y\,
\frac{g_s(x,y)v(|y-f(s)|)g_{\beta-s}(y,x)}{g_\beta(x,x)}\vadjust{\goodbreak}\\
&&\qquad=\frac{1}{|\Lambda|}\int_{\Lambda} \dd y\, v\bigl(|y-f(s)|\bigr)\int
_{\Lambda} \dd x \biggl(\frac{g_{\beta-s}(y,x)g_s(x,y)}{g_\beta
(y,y)}\biggr)\frac{g_\beta(y,y)}{g_\beta(x,x)}\nonumber\\
&&\qquad=\frac{1}{|\Lambda|}\int_{\Lambda} \dd y\, v\bigl(|y-f(s)|\bigr),\nonumber
\end{eqnarray}
since, because of $ g_\beta(x,x)=g_\beta(y,y) $, the integral over $
x $ is exactly $1$.
An upper bound follows easily because the interaction potential is
nonnegative, that is,
%
%
\begin{eqnarray}\label{energyesti2}\qquad\quad
\E_{\Lambda}^{(\beta)}\bigl(v\bigl(|B_s-f(s)|\bigr)\bigr) &=&\frac
{1}{|\Lambda|}\int_{\Lambda} \dd y\, v\bigl(|y-f(s)|\bigr)\le\frac
{1}{|\Lambda|}\int_{\R^d} v(|x|) \,\dd x\nonumber\\[-8pt]\\[-8pt]
&=& \frac{1}{|\Lambda|}\alpha(v).\nonumber
\end{eqnarray}
Using this in (\ref{energyest}), we get
\begin{eqnarray*}
F_2(T(\lambda))&\ge&\Biggl(\bigotimes_{k=1}^N\bigl(\E_{\Lambda
}^{(k\beta)}\bigr)^{\otimes\lambda_k}\Biggr)\bigl[\ee^{-\Gcal
_{N,\beta}^{(\lambda)}}\ee^{-\sum_{k\in\N}\sum
_{i=1}^{\lambda_k}\sum_{j=1}^k\beta{1}/{|\Lambda|}\alpha
(v)}\bigr]\\
&=& F_2(\lambda)\ee^{-{N\beta}/{|\Lambda|}\alpha(v)}.
\end{eqnarray*}
Using this and (\ref{F1esti}) in (\ref{lowerest}), the assertion follows.
\end{pf}

Now we draw two corollaries. First, we give an upper bound for the free
energy, introduced in (\ref{limfreeener}). Recall that $\Lambda
_{L_N}$ is the centered box with volume $N/\rho$.
\begin{cor}[(Upper bound for the free energy)]\label{cor-upf} For any
$\beta, \rho\in(0, \infty)$,
\[
\limsup_{N \to\infty}-\frac1\beta\frac{1}{|\Lambda_{L_N}|} \log
Z_N(\beta,\Lambda_{L_N}) \le\frac{\rho}{\beta} \log( \rho
(4 \pi\beta)^{d/2} ) + \rho^2\alpha(v).
\]
\end{cor}
\begin{pf}
We use Lemma \ref{lem-Zquotient} iteratively, to get
\begin{eqnarray*}
Z_N(\beta, \Lambda_{L_N}) &=& \prod_{l=0}^{N-1} \frac{Z_{l+1}(\beta
, \Lambda_{L_N})}{Z_l(\beta, \Lambda_{L_N})}
\ge\prod_{l=0}^{N-1} \biggl( (4 \pi\beta)^{-d/2} \frac1{\rho}
\ee^{- \beta\alpha(v) \rho} \biggr)\\
&=& \biggl( (4 \pi\beta)^{-d/2} \frac1{\rho} \ee^{- \beta\alpha
(v) \rho}\biggr)^N.
\end{eqnarray*}
The assertion follows by taking $\limsup_{N \ti}-\frac1\beta\frac
1{|\Lambda_{L_N}|} \log$.
\end{pf}
\begin{cor}\label{cor-ZNreduction} Fix $(\beta,\rho)\in\Dcal_v$.
Then, for any $N,\widetilde N\in\N$ satisfying $\widetilde N\leq N$,
\[
\ttE\bigl[ \ee^{-H_{\Lambda_{L_N}}(\omega_{\romP})}
\mathbh{1}\bigl\{
N^{(\ell)}_{\Lambda_{L_N}}(\omega_{\romP})= N\bigr\} \bigr]
\geq\ttE\bigl[ \ee^{-H_{\Lambda_{L_N}}(\omega_{\romP})}
\mathbh
{1}\bigl\{N^{(\ell)}_{\Lambda_{L_N}}(\omega_{\romP})= \widetilde N\bigr\}
\bigr].
\]
In particular, the map $\widetilde N\mapsto Z_{\widetilde N}(\beta
,\Lambda_{L_N})$ is increasing in $\widetilde N\in\{1,\ldots,N\}$.
\end{cor}
\begin{pf}
Observe that, for $l<N$, by Lemma \ref{lem-Zquotient},
\[
\frac{Z_{l+1}(\beta,\Lambda_{L_N})}{Z_{l}(\beta,\Lambda_{L_N})}
\geq(4\pi\beta)^{-d/2}\frac{|\Lambda_{L_N}|}{l+1}
\ee
^{-l\beta\alpha(v)/|\Lambda_{L_N}|}
\geq(4\pi\beta)^{-d/2}\frac1\rho\ee^{-\beta
\rho\alpha(v)}
\geq1,
\]
where the last step follows from $(\beta,\rho)\in\Dcal_v$. Hence,
for any $\widetilde N\in\N$ satisfying $\widetilde N\leq N$, we have
$Z_{N}(\beta,\Lambda_{L_N})\geq Z_{\widetilde N}(\beta,\Lambda
_{L_N})$. Now use Proposition \ref{lem-rewrite} to finish.
\end{pf}

\subsubsection*{Openness}

As we already mentioned, some of the technical difficulties
for the application of Varadhan's lemma come from the fact that the set
$\{P\in\Pcal_\theta\dvtx\langle P,N_U^{(\ell)}\rangle=\rho\}
$ is not open. This problem will be taken care of in the following
lemma: we derive a lower bound for the right-hand side in (\ref
{LDPmainstatementlower}) in terms of the same expectation, where the
strict condition $=\rho$ is replaced by the condition $\in(\rho
-\delta,\rho+\delta)$, for some $\delta> 0$. Though this set is not
open in $\Pcal_{\theta}$, it will be open after restricting $\Omega$
to some cut-off version $\Omega^{(K,R)}$, which we will
introduce a bit later.
\begin{lemma}\label{lem-makeopen}
Fix $\beta,\rho\in(0,\infty)$. We abbreviate $\RR_{N}(\omega) =
\RR_{\Lambda_{L_N},\omega}$ for $\omega\in\Omega$. Fix $ \delta
\in(0, \rho) $. Then for any $N\in\N$,
%
%
\begin{eqnarray}\label{lowbound1}\qquad
&&\ttE\bigl[ \ee^{-H_{\Lambda_{L_N}}(\omega_{\romP})}
\mathbh{1}\bigl\{
\bigl\langle\RR_{N}(\omega_{\romP}), N^{(\ell)}_{U}\bigr\rangle= \rho\bigr\}
\bigr] \nonumber\\
&&\qquad\ge\frac{(C_1\wedge C_2)^{\delta|\Lambda_{L_N}|}}{2\delta
|\Lambda_{L_N}|+2}
\\
&&\qquad\quad{}\times\ttE \bigl[ \ee^{-|\Lambda_{L_N}| \langle
\RR_{N}(\omega
_{\romP}),\Phi_\beta\rangle} \mathbh{1}\bigl\{\bigl\langle\RR_{N}(\omega
_{\romP}), N^{(\ell)}_{U}\bigr\rangle\in
(\rho-\delta,\rho+\delta)\bigr\}\bigr],\nonumber
\end{eqnarray}
where $C_1= 1 \wedge(\ee^{-(\rho+\delta)\beta
\alpha(v)}(4\pi
\beta)^{-d/2}\frac{1}{\rho+\delta})$ and $ C_2= \ee^{-
{\overline q}/({\rho-\delta})}$.
\end{lemma}
\begin{pf}
Define the subset
\[
\Pcal_l=\biggl\{P\in\Pcal_\theta\dvtx\bigl\langle P,N_U^{(\ell)}\bigr\rangle
=\frac{l}{|\Lambda_{L_N}|}\biggr\}
\]
of probability measures. Abbreviate
%
%
\begin{eqnarray}
Y^{(1)}_l &=& \ttE\bigl[ \ee^{-H_{\Lambda
_{L_N}}(\omega_{\romP})} \mathbh{1}_{\Pcal_l}(\RR_{N}(\omega_{\romP}))
\bigr], \\
Y^{(2)}_l &=&\ttE\bigl[ \ee^{-|\Lambda_{L_N}|
\langle\RR
_{N}(\omega_{\romP}),\Phi_\beta\rangle} \mathbh{1}_{\Pcal_l}(\RR
_{N}(\omega_{\romP})) \bigr].
\end{eqnarray}
Notice that, since $N/|\Lambda_{L_N}|=\rho$, the left-hand side of
(\ref{lowbound1}) is equal to $Y^{(1)}_N$, while the expectation
on the right-hand side is equal to
\[
\sum_{l\in\N\dvtx(\rho-\delta)|\Lambda_{L_N}|< l< (\rho+\delta
)|\Lambda_{L_N}|}Y^{(2)}_l.
\]
We now estimate the quotients $Y^{(1)}_{l+1}/Y^{(1)}_l$,
respectively, $Y^{(2)}_{l+1}/Y^{(2)}_l$, from below and above.
More precisely, we show, for any $l\in\N_0$,
%
%
\begin{equation}\label{quoZ1}
Y^{(1)}_{l+1}\ge C_1 Y^{(1)}_l \qquad\mbox{if } (\rho-\delta
)|\Lambda_{L_N}|< l\leq\rho|\Lambda_{L_N}|
\end{equation}
and
%
%
\begin{equation}\label{quoZ2}
Y^{(2)}_{l}\ge C_2 Y^{(2)}_{l+1} \qquad\mbox{if } \rho
|\Lambda_{L_N}|\leq l< (\rho+\delta)|\Lambda_{L_N}|.
\end{equation}

The proof of (\ref{quoZ1}) follows from Lemma \ref{lem-Zquotient},
combined with Proposition \ref{lem-rewrite}. Now we prove (\ref{quoZ2}).

We find a map\vspace*{1pt} $ \Tcal\dvtx\Pcal_{l+1}\to\Pcal_l$ that describes a
thinning procedure with the parameter $ p=\frac{l}{l+1} $. To this
end, we introduce a probability kernel $K$ from $\Omega$ to $\Omega$
by putting $K(\omega,\cdot)$ equal to the distribution of $\omega
^{(\eta)}=\sum_{x\in\xi}\eta_x\delta_{(x,f_x)}=\sum_{x\in
\xi^{(\eta)}}\delta_{(x,f_x)}$, where $\omega=\sum_{x\in\xi
}\delta_{(x,f_x)}\in\Omega$, and, given $\omega$, $(\eta_x)_{x\in
\xi}$ is a Bernoulli sequence with parameter $p$. The mapping
%
%
\begin{equation}
\Tcal\dvtx\Pcal_{l+1}\to\Pcal_l,\qquad \Tcal(P)=PK,
\end{equation}
describes the distribution of what is left from a configuration with
distribution $P$ after deleting each particle independently with
probability $p$. Given $P\in\Pcal_{l+1} $, it follows, writing $\ttE
_\eta$ for the expectation with respect to $(\eta_x)_{x\in\xi}$,
\begin{eqnarray*}
\bigl\langle\Tcal(P),N_U^{(\ell)}\bigr\rangle&=&\int_\Omega P(\dd\omega
)\int_\Omega K(\omega, \dd\widetilde\omega) N_U^{(\ell
)}(\widetilde\omega)\\
&=&\int_\Omega P(\dd\omega) \ttE_\eta\biggl[N_U^{(\ell)}
\biggl(\sum_{x\in\xi}\eta_x\delta_{(x,f_x)}\biggr)\biggr]\\
&=&\int_\Omega P(\dd\omega) \ttE_\eta\biggl[\sum_{x\in\xi\cap
U}\eta_x\ell(f_x)\biggr]\\
&=&p\bigl\langle P,N_U^{(\ell)}\bigr\rangle
=\frac{l}{l+1}\bigl\langle P,N_U^{(\ell)}\bigr\rangle
=\frac{l}{|\Lambda_{L_N}|},
\end{eqnarray*}
which shows that $\Tcal\dvtx\Pcal_{l+1} \to\Pcal_{l} $ is well defined.
Since $\Tcal$ removes particles, and therefore energy, the estimate
%
%
\begin{equation}\label{energyineq}
\langle P, \Phi_\beta\rangle\ge\langle\Tcal(P),\Phi_\beta
\rangle,\qquad P\in\Pcal_{l+1},
\end{equation}
follows easily.
Inequality (\ref{energyineq}) gives the estimate
%
%
\begin{eqnarray}\label{upperboundZ}
Y^{(2)}_{l+1}&\le&\ttE\bigl[\ee^{-|\Lambda
_{L_N}|\langle\Tcal
(\RR_N(\omega_{\romP})),\Phi_\beta\rangle}\mathbh{1}_{\Pcal
_l}(\Tcal(\RR_N(\omega_{\romP})))\bigr]\nonumber\\[-8pt]\\[-8pt]
&=&\int_{\Pcal_l} \ee^{-|\Lambda_{L_N}|\langle P,\Phi
_\beta\rangle
}\frac{\dd\ttQ \circ\RR_N^{-1}\circ\Tcal^{-1}}{\dd\ttQ \circ
\RR_N^{-1}}(P) \ttQ \circ\RR_N^{-1}(\dd
P),\nonumber
\end{eqnarray}
where we recall that $\ttQ $ and $\ttE$ are the distribution of
and expectation with respect to the marked Poisson process $\omega
_{\romP}$, and we conceive $\RR_N$ as a map $\Omega\to\Pcal_\theta
$; note that $\RR_N$ depends only on the configuration in $\Lambda_{L_N}$.

Now we identify the corresponding Radon--Nikodym density $\varphi
_N=\dd\ttQ \circ\RR_N^{-1}\circ\Tcal^{-1}/\dd\ttQ \circ\RR
_N^{-1}$ on the image $\RR_N(\Omega)$. We claim that
%
%
\begin{equation}\label{density}
\varphi_N(\RR_{N}(\omega))=p^{\#(\xi\cap\Lambda
_{L_N})}\ee
^{(1-p)\overline q |\Lambda_{L_N}|},\qquad \omega\in\Omega.
\end{equation}
This is shown as follows. Note that $\varphi_N$ is the density of
$\Tcal(\RR_N(\omega_{\romP}))$ with respect to $\RR_N(\omega_{\romP})$
and that $\Tcal(\RR_N(\omega_{\romP}))$ has the distribution of
$\RR_N(\omega_{\romP}^{(\eta)})$. Recall that the particle
process $\xi_{\romP}\cap\Lambda_{L_N}$ is a standard Poisson process
on $\Lambda_{L_N}$ with intensity $\overline q |\Lambda_{L_N}|$, and
$\xi_{\romP}^{(\eta)}\cap\Lambda_{L_N}$ has intensity
$p\overline q |\Lambda_{L_N}|$. It is standard that the right-hand
side of (\ref{density}) is the density of $\xi_{\romP}^{(\eta)}\cap
\Lambda_{L_N}$ with respect to $\xi_{\romP}\cap\Lambda
_{L_N}$. But this implies that (\ref{density}) holds, as we have, for
any nonnegative measurable test function $g\dvtx\Pcal\to[0,\infty]$,
\begin{eqnarray*}
\int g(P) \ttQ \circ\Tcal(\RR_N)^{-1}(\dd P)
&=&\ttE[g(\Tcal(\RR_N(\omega_{\romP})))]
=\ttE\bigl[\ttE_\eta\bigl[g\bigl(\RR_N\bigl(\omega_{\romP}^{(\eta
)}\bigr)\bigr)\bigr]\bigr]\\
&=&\ttE\bigl[p^{\#(\xi_{\romP}\cap\Lambda_{L_N})}\ee
^{(1-p)\overline q |\Lambda_{L_N}|}g(\RR_N(\omega_{\romP}))\bigr]\\
&=&\int_\Omega p^{\#(\xi\cap\Lambda_{L_N})}\ee^{(1-p)\overline q
|\Lambda_{L_N}|}g(\RR_N(\omega)) \ttQ (\dd\omega).
\end{eqnarray*}

Note that, for $(\rho-\delta)|\Lambda_{L_N}|<l\leq\rho|\Lambda_{L_N}|$,
\[
\varphi_N(\RR_{N}(\omega))\le\ee^{(1-p)\overline q
|\Lambda
_{L_N}|}= \ee^{{\overline q}/({l+1})|\Lambda_{L_N}|}
\le\ee^{
{\overline q}/({\rho-\delta})},\qquad \omega\in\Omega.
\]
Hence, from (\ref{upperboundZ}) we have
\[
Y^{(2)}_{l+1}\le\ee^{{\overline q}/({\rho-\delta
})} \int
_{\Pcal_l} \ee^{-|\Lambda_{L_N}|\langle P,\Phi_\beta
\rangle}\ttQ \circ\RR_N^{-1}(\dd P)
=\ee^{{\overline q}/({\rho-\delta})} Y^{(2)}_l,
\]
and thus the estimate (\ref{quoZ2}).


Now we finish the proof of the lemma subject to (\ref{quoZ1}) and
(\ref{quoZ2}). By Lem\-ma~\ref{lem-Hloweresti}(ii), we have
$Y^{(1)}_N\geq Y^{(2)}_N$ and therefore
\begin{eqnarray*}
&&\mbox{left-hand side of (\ref{lowbound1})}\\
&&\qquad= Y^{(1)}_N
\geq
\frac{1}{2\delta|\Lambda_{L_N}|+2}\\
&&\qquad\hphantom{= Y^{(1)}_N
\geq}{}\times\biggl(\sum_{(\rho-\delta
)|\Lambda_{L_N}|< l\le\rho|\Lambda_{L_N}|}Y^{(1)}_N
+\sum_{\rho|\Lambda_{L_N}|< l< (\rho+\delta)|\Lambda_{L_N}|}
Y^{(2)}_N\biggr).
\end{eqnarray*}
For $ (\rho-\delta)|\Lambda_{L_N}|< l\le\rho|\Lambda_{L_N}| $ the
estimate (\ref{quoZ1}) gives
\[
Y^{(1)}_N\ge C_1 Y^{(1)}_{N-1}\ge\cdots\ge C_1^{\delta
|\Lambda_{L_N}|} Y^{(1)}_l\ge C_1^{\delta|\Lambda_{L_N}|}
Y^{(2)}_l,
\]
because $C_{1} \le1$, where we again used Lemma \ref
{lem-Hloweresti}(ii). On the other hand, for $ \rho|\Lambda_{L_N}|<
l<(\rho+\delta)|\Lambda_{L_N}| $ the estimate (\ref{quoZ2}) gives
\[
Y^{(2)}_N\ge C_2 Y^{(2)}_{N+1}\ge\cdots\ge C_2^{\delta
|\Lambda_{L_N}|} Y^{(2)}_l,
\]
where we used $C_{2} < 1$.
Therefore
%
%
\begin{eqnarray}
Y^{(1)}_N&\ge&\frac{(C_1\wedge C_2)^{\delta|\Lambda
_{L_N}|}}{2\delta|\Lambda_{L_N}|+2}\sum_{(\rho-\delta)|\Lambda
_{L_N}|< l< (\rho+\delta)|\Lambda_{L_N} } Y^{(2)}_l\nonumber\\[-8pt]\\[-8pt]
&=&\mbox
{right-hand side of (\ref{lowbound1})},\nonumber
\end{eqnarray}
which finishes the proof of the lemma.
\end{pf}

As a conclusion of Lemma \ref{lem-makeopen} we have the following
lower bound for any sufficiently large $N\in\N$:
%
%
\begin{eqnarray}\label{Zlowbound3}
Z_N(\beta, \Lambda_{L_N}) &\ge&\ee^{|\Lambda_{L_N}|
(\overline q -
C \delta)} \nonumber\\
&&{}\times\ttE \bigl[ \ee^{-|\Lambda_{L_N}|
\langle\RR
_{N}(\omega_{\romP}),\Phi_\beta\rangle} \\
&&\hspace*{19.8pt}{}\times \mathbh{1}\bigl\{\bigl\langle\RR_{N}(\omega_{\romP}), N^{(\ell
)}_{U}\bigr\rangle\in
(\rho-\delta,\rho+\delta)\bigr\}\bigr],\nonumber
\end{eqnarray}
for any $\delta\in(0, \frac\rho2)$ and some $C$ depending only on
$\beta, \rho$ and $v$.
Furthermore, if $(\beta, \rho) \in\Dcal_v$, then we can combine
Lemma \ref{lem-makeopen} with
Corollary \ref{cor-ZNreduction} to get, for any $\widetilde\rho\in
(0, \rho]$ and any $\delta\in(0, \frac{\widetilde\rho}2)$, for
any sufficiently large $N\in\N$,
%
%
\begin{eqnarray}\label{Zlowbound4}
Z_N(\beta, \Lambda_{L_N}) &\ge& \ee^{|\Lambda_{L_N}|
(\overline q -
C \delta)} \nonumber\\
&&{}\times\ttE \bigl[ \ee^{-|\Lambda_{L_N}|
\langle\RR
_{N}(\omega_{\romP}),\Phi_\beta\rangle} \\
&&\hspace*{19.8pt}{}\times \mathbh{1}\bigl\{\bigl\langle\RR_{N}(\omega_{\romP}), N^{(\ell)}_{U}\bigr\rangle
\in
(\widetilde\rho-\delta,\widetilde\rho+\delta)\bigr\}\bigr].\nonumber
\end{eqnarray}
Hence, in order to prove both bounds in (\ref{freeenlower}), it is
enough to prove
%
%
\begin{eqnarray}\label{mainaim}\quad
&&\liminf_{\delta\downarrow0}\liminf_{N\to\infty}\frac1{|\Lambda
_{L_N}|}\log\ttE \bigl[ \ee^{-|\Lambda_{L_N}|
\langle\RR
_{N}(\omega_{\romP}),\Phi_\beta\rangle} \nonumber\\
&&\hspace*{114.40pt}{}\times \mathbh{1}\bigl\{\bigl\langle\RR_{N}(\omega_{\romP}), N^{(\ell
)}_{U}\bigr\rangle\in(\rho-\delta,\rho+\delta)\bigr\}\bigr]\\
&&\qquad\geq-\chi^{(=)}(\beta,\rho),\nonumber
\end{eqnarray}
for any $\beta,\rho\in(0,\infty)$, since $\chi^{(\leq)}(\beta
,\rho)=\inf_{\widetilde\rho\in(0,\rho)}\chi^{(=)}(\beta
,\rho)$.

\subsubsection*{Restriction of the mark space}

We will approximate the mark space $E$ by the cut-off version
\[
E^{(K,R)}:= \bigcup_{k=1}^K \Ccal_{k,R}\qquad \mbox{where }
\Ccal_{k,R}:= \Bigl\{ f \in\Ccal_k \dvtx{\sup_{s \in[0, k \beta]}}
|f(s)-f(0)| \le R\Bigr\}.
\]
Let\vspace*{1pt} $\Omega^{(K,R)}$ denote the set of locally finite point
measures on $\R^d\times E^{({K,R})}$. Define the canonical
projection $\pi_{K,R}\dvtx\Omega\to\Omega^{(K,R)}$ by
\[
\pi_{K,R}(\omega)=\omega^{(K,R)}=\sum_{x\in\xi\dvtx f_x\in
E^{(K,R)}}\delta_{(x,f_x)}.
\]
On $\Omega^{(K,R)}$ we consider the Poisson point process
%
%
\begin{equation}
\omega_{\romP}^{(K,R)}=\pi_{K,R}(\omega_{\romP})=\sum_{x \in
\xi_{\romP}\dvtx B_x\in E^{(K,R)}} \delta_{(x, B_x)}
\end{equation}
as the reference process. The distribution of $\omega_{\romP}^{({
{K,R}})}$ is denoted $\ttQ ^{(K,R)}$, its intensity measure is
$\nu^{(K,R)}=\sum_{k=1}^K\nu_k^{(K,R)}$, where $\nu
_k^{(K,R)}$ is the restriction of $\nu_k$ to $\Omega^{(K,R)}$; see
(\ref{nudef}). By $I_{\beta}^{(K,R)}$ we denote
the rate function with respect to $\omega_{\romP}^{(K,R)}$, that
is, $I_{\beta}^{(K,R)}$ is defined as $I_\beta$ in (\ref{Idef}) with
$\omega_{\romP}$ replaced by $\omega_{\romP}^{(K,R)}$. If
there is no confusion possible, we identify the set $\Pcal
_\theta(\Omega^{(K,R)})$ of shift-invariant marked random point
fields on $\Omega^{(K,R)}$ with the set of those $P\in\Pcal
_\theta=\Pcal_\theta(\Omega)$ that are concentrated on $\Omega
^{(K,R)}$.
A variant of Lemma \ref{ldpyn} gives that $(\RR_{\Lambda_L,\omega
_{\romP}^{(K,R)}})_{L>0}$ satisfies the large-deviations
principle with rate function $I_{\beta}^{(K,R)}$. Observe that
$\RR_{\Lambda_L,\omega_{\romP}^{(K,R)}}=\RR_{\Lambda
_L,\omega_{\romP}}\circ\pi_{K,R}^{-1}$. Hence, according to the
contraction principle, we have the identification
%
%
\begin{equation}\label{IKident}
I_{\beta}^{(K,R)}(P)=\inf\{I_\beta(Q)\dvtx Q\in\Pcal_\theta
, Q\circ\pi_{K,R}^{-1}=P\},
\end{equation}
since the map $Q\mapsto Q\circ\pi_{K,R}^{-1}$ is continuous.

For a while, we keep $K$ and $R$ fixed. Now we work on the expectation
on the right-hand side of (\ref{LDPmainstatementlower}). We obtain a
lower bound by requiring that $\RR_{\Lambda_{L_N}, \omega_{\romP}}$
be concentrated on $\Omega^{(K,R)}$. On this event, we may
replace $\RR_{\Lambda_{L_N}, \omega_{\romP}}$ by $\RR_{\Lambda
_{L_N}, \omega_{\romP}^{(K,R)}}$, and we may replace the
expectation $\ttE $ with respect to the Poisson process $\omega
_{\romP}$ by the expectation $\ttE^{(K,R)} $ with respect to
$\omega_{\romP}^{(K,R)}$. We write $ \RR_N $ for $ \RR
_{\Lambda_{L_N}, \omega_{\romP}^{(K,R)}} $ in the following.
Hence, we can extend (\ref{Zlowbound3}) by
%
%
\begin{eqnarray}\label{reduc1}
Z_N(\beta,\Lambda_{L_N})
&\ge&\ee^{|\Lambda_{L_N}| (\overline q - C \delta)}
\nonumber\\
&&{}\times\ttE^{(K,R)} \bigl[ \ee^{-|\Lambda_{L_N}| \langle
\RR_{N},\Phi_\beta
\rangle} \\
&&\hspace*{43.3pt}{}\times\mathbh{1}\bigl\{\bigl\langle\RR_{N}, N^{(\ell)}_{U}\bigr\rangle
\in(\rho- \delta, \rho+ \delta)
\bigr\}\bigr].\nonumber
\end{eqnarray}
Notice that $\{P \in\Pcal_\theta(\Omega^{(K,R)}) \dvtx
\langle P, N^{(\ell)}_{U}\rangle\in(\rho- \delta, \rho+
\delta) \}$ is an open set.
In order to apply the lower bound of Varadhan's lemma to the right-hand
side, we need to have that the map $P\mapsto\langle P,\Phi_{\beta
}\rangle$ is upper semicontinuous. This will be achieved by a further
restriction procedure.

\subsubsection*{Continuity}

We prove the continuity of the map $ P\mapsto\langle P,\Phi
_\beta\rangle$ on the following suitable subset of measures. For
$r\in(0,\infty)$, put
%
%
\begin{eqnarray}
\Gamma_r&=&\biggl\{\omega\in\Omega^{(K,R)}\dvtx T_{x,y}(\omega
)\leq r\ \forall x,y \in\xi, \nonumber\\[-8pt]\\[-8pt]
&&\hspace*{6.9pt} \mbox{and } |x-y|\geq\frac1r \mbox{ for all
distinct }x,y\in\xi
\biggr\},\nonumber
\end{eqnarray}
where $T_{x,y}(\omega)$ was defined in (\ref{Tdef}). Denote
\[
\Pcal_{\theta, r} := \bigl\{ P \in\Pcal_{\theta}\bigl(\Omega^{(K,R)}\bigr)
\dvtx P(\Gamma_{r} ) =1 \bigr\}.
\]
In the following lemma we use that the map $ t \mapsto t^{d-1} \sup_{s
\ge t - 2R} v(s)$ is integrable, which easily follows from the
temperedness assumption in Assumption \ref{AssumptionV}.
\begin{lemma}\label{lem-Fcont} For any $r>0$, the map $P\mapsto
\langle P, \Phi_\beta\rangle$ is continuous on the set $\Pcal
_{\theta, r}$.
\end{lemma}
\begin{pf}
We adapt the proof of the lower bound in \cite{G94}, Theorem 2.
Recall that $\pi_{n} \dvtx\Omega\to\Omega_{2n}$ denotes the
projection $\pi_{n}(\omega) = \sum_{x \in\xi\cap\Lambda_{2n}}
\delta_{(x, f_{x})}$ on the box $\Lambda_{2n} = [-n , n]^{d}$. For
any $P$ let $P_{n} := P \circ\pi^{-1}_{n}$. Let $P$ and a net
$(P^{(\alpha)})_{\alpha\in D}$ be in $\Pcal_{\theta, r}$ such
that $P^{(\alpha)}$ converges to $P$ (in the topology $\tau_\Lcal
$). Then we have, for any $n\in\N$ and $\alpha\in D$,
%
%
\begin{eqnarray}\label{FvO}
&&\bigl| \langle P, \Phi_{\beta} \rangle- \bigl\langle P^{(\alpha)}, \Phi
_{\beta} \bigr\rangle\bigr| \nonumber\\
&&\qquad\le|\langle P, \Phi_{\beta}- \Phi_{\beta} \circ\pi_{n }\rangle
| + \bigl| \bigl\langle P^{(\alpha)}-P, \Phi_{\beta}\circ\pi_n \bigr\rangle
\bigr|\nonumber\\[-8pt]\\[-8pt]
&&\qquad\quad{} + \sup_{\alpha\in D} \bigl|\bigl\langle P^{(\alpha)}, \Phi
_{\beta} - \Phi_{\beta} \circ\pi_{n} \bigr\rangle\bigr|\nonumber\\
&&\qquad\leq\bigl| \bigl\langle P^{(\alpha)}-P, \Phi_{\beta}\circ\pi_n
\bigr\rangle\bigr| + 2\sup_{\widetilde P\in\Pcal_{\theta,r}}\langle
\widetilde P, |\Phi_{\beta}- \Phi_{\beta} \circ
\pi_{n}|\rangle.\nonumber
\end{eqnarray}
Observe that the last term on the right-hand side vanishes as $n \to
\infty$ since $\Phi_{\beta} \circ\pi_{n}$ converges to $\Phi
_\beta$ uniformly on $\Gamma_r$. Indeed, for $\omega\in\Gamma_{r}$ estimate
%
%
\begin{eqnarray}
\Phi_{\beta}(\omega) - \Phi_{\beta}(\pi_{n}(\omega))
&=& \sum_{x \in U \cap\xi} \sum_{y \in\xi\cap\Lambda_{2n}^{c}}
T_{x,y}(\omega)\nonumber\\[-8pt]\\[-8pt]
& \le&\frac12 \sum_{x \in U \cap\xi} \sum_{y \in\xi\cap\Lambda
_{2n}^{c}} K^{2} \beta\sup_{s \ge|x- y| - 2R }
v(s),\nonumber
\end{eqnarray}
where we also used that $\ell(f_{x}) \le K$ and ${\sup_{s \in[0,\beta
\ell(f_x)]}} |f_{x}(s)-f_x(0)|\le R$ for any $ x \in\xi$, since $
\omega\in\Omega^{(K,R)}$. Since $|x-y| \ge\frac1r $ for any
distinct $x,y \in\xi$, the upper bound is not larger than
\[
K^{2} \beta C_{r,R} \int_{n}^{\infty} t ^{d-1} \sup_{s \ge t - 2R}
v(s) \,\dd t,
\]
for some constant $C_{r,R}$ that depends only on $r$ and $R$. Now use
that mapping $ t \mapsto t^{d-1} \sup_{s \ge t - 2R} v(s)$ is integrable.

For any $n$, the first term on the right-hand side of (\ref{FvO})
vanishes asymptotically since the net $(P^{(\alpha)})_{\alpha\in
D}$ converges to $P$, and $\Phi_\beta\circ\pi_n$ is local and
bounded on $\Gamma_r$.
\end{pf}

\subsubsection*{Ergodic approximation}

As a preparation for the construction of an ergodic
approximation, we now show that any $P$ with finite energy is tempered,
that is, the expectation of the square of the mean-particle density is
finite. Here we use the assumption that $\liminf_{r\downarrow
0}v(r)>0$, which is part of Assumption \ref{AssumptionV}. Hence, we may
pick $R^*>0 $
and $ \zeta> 0$ such that $ v(|x|)\ge\zeta$ for all $ |x|\le R^* $.
\begin{lemma}[(Temperedness)]\label{lem-Temp} Fix $K,R \in\N
$, and let $ P\in\Pcal_\theta(\Omega^{(K,R)}) $ with $
\langle P,\Phi_\beta\rangle<\infty$. Then
\[
\langle P,N_U^2\rangle<\infty\quad\mbox{and}\quad \bigl\langle
P,\bigl(N_U^{(\ell)}\bigr)^2\bigr\rangle<\infty.
\]
\end{lemma}
\begin{pf}
We may assume that $R^*<\frac12$. Therefore, we obtain a lower bound
for $ \langle P, \Phi_\beta\rangle$ by restricting the sums on $x,y
$ to $x,y\in\Lambda_{R^*/4}=[-\frac{R^*}4, \frac{R^*}4]^d$ and by
dropping all the parts of the cycles except for the first one,
%
%
\begin{eqnarray}\label{temp1}\quad
\langle P,\Phi_\beta\rangle
&=&\frac{1}{2}\int P(\dd\omega) \sum_{x\in\xi\cap U,y\in\xi
}\sum_{i=0}^{\ell(f_x)-1}\sum_{j=0}^{\ell(f_y) -1}\mathbh{1}_{\{
(x,i)\not=(y,j)\}}\nonumber\\
&&\hspace*{159pt}{}\times \int_0^\beta
v\bigl(|f_x(i\beta+s)\nonumber\\[-8pt]\\[-8pt]
&&\hspace*{200.7pt}{}-f_y(j\beta+s)|\bigr) \,\dd s\hspace*{-15pt}\nonumber\\
&\ge&\frac{1}{2}\int P(\dd\omega) \sum_{x,y\in\xi\cap\Lambda
_{R^*/4}} \mathbh{1}\{x \not=y\}\int_0^\beta v\bigl(|f_x(s)-f_y(s)|\bigr) \,\dd
s.\nonumber
\end{eqnarray}
Define, for any $\omega\in\Omega^{(K,R)}$ and $ x\in\xi$,
%
%
\begin{equation}
\tau_x(\omega)=\inf\{s\in[0,\beta]\dvtx|f_x(s)-x| > R^*/4\}\wedge
\delta.
\end{equation}
Note that $|x-y|\leq R^*/2$ on the right-hand side of (\ref{temp1}).
Since $ v(|x|)\ge\zeta$ for all $ |x|\le R^* $, each integral on the
right-hand side of (\ref{temp1}) can be estimated from below as
follows:
\begin{eqnarray*}
\int_0^\beta v\bigl(|f_x(s)-f_y(s)|\bigr) \,\dd s&\ge&\int_0^{\tau_x(\omega
)\wedge\tau_y(\omega)} v\bigl(|f_x(s)-f_y(s)|\bigr) \,\dd s\\
&\geq&\zeta\bigl(\tau_x(\omega)\wedge\tau_y(\omega)\bigr),\qquad x\in\xi
^{(k)}, y\in\xi^{(k')}.
\end{eqnarray*}
We get a further lower bound in (\ref{temp1}) by inserting the
indicator on the event $\{\tau_x=\delta=\tau_y\}$
\[
\langle P,\Phi_\beta\rangle\geq\frac{\delta\zeta}2 \int P(\dd
\omega) \#\{(x,y)\in(\xi\cap\Lambda_{R^*/4}
)^2\dvtx x\not=y,\tau_x=\delta=\tau_y\}.
\]
Since the event $\{\tau_x=\delta\}$ is decreasing for decreasing
$\delta$ and its probability tends to one as $\delta\downarrow0$,
the above counting variable tends to the number of distinct pairs in
$\xi\cap\Lambda_{R^*/4}$. Hence, for some sufficiently small $\delta
>0$, we have
\[
\langle P,\Phi_\beta\rangle\geq\frac{\delta\zeta}4
\int P(\dd\omega) \#\{(x,y)\in(\xi\cap\Lambda
_{R^*/4})^2\dvtx x\not=y\}
\geq\frac{\delta\zeta}8 \langle P, N_{\Lambda_{R^*/4}}^2
\rangle.
\]
Hence, if $\langle P,\Phi_\beta\rangle$ is finite, then, by
shift-invariance of $P$, also $\langle P, N_{\Lambda}^2\rangle$ is
finite for any bounded box $\Lambda$. Since $P$ is concentrated on
configurations with bounded leg length, also $\langle P, (N_{\Lambda
}^{(\ell)})^2\rangle$ is finite for any bounded box $\Lambda$.
\end{pf}

Now we approximate any probability measure on $\Omega^{(K,R)}$
with an ergodic measure. Define
%
%
\begin{equation}
\psi_{R}(t):=
\cases{\displaystyle
\sup_{s \ge t - 2R} v(s), &\quad if $t \ge3R$,\vspace*{2pt}\cr
v(R), &\quad if $t \in[0,3R]$.
}
\end{equation}
Recall from Assumption \ref{AssumptionV} that $\psi_{R}(t) = O(t^{-h})$
for some $h>d$.
\begin{lemma}[(Ergodic approximation)]\label{lem-ergodappr}
Fix $K,R \in\N$ and $\eps>0$. Then, for any $P\in\Pcal_\theta
(\Omega^{(K,R)})$ satisfying $ I_{\beta}^{(K,R)}(P) +
\Phi_\beta(P) < \infty$ and for any neighborhood $V$ of $P$ in
$\Pcal_\theta(\Omega^{(K,R)})$, there exists an ergodic
measure $\widetilde P\in V$ and some $r>0$ such that $\widetilde
P(\Gamma_r)=1$, and $\langle\widetilde P,\Phi_\beta\rangle\leq
\langle P,\Phi_\beta\rangle+\eps$ and $I_{\beta
}^{(K,R)}(\widetilde P)\leq I_{\beta}^{(K,R)}(P)+\eps$.
\end{lemma}
\begin{pf}
This is similar to \cite{G94}, Lemma 5.1. Recall that $P_{n} $ denotes
the projection of $P$ on $\Omega_{n}$, the configuration space on the
box $\Lambda_{2n}= [-n,n]^{d}$. Since $\langle P, \Phi_{\beta}
\rangle<\infty$, and as $\Phi_{\beta} \ge0$, we have $\langle
P_{n}, \Phi_{\beta} \rangle< \infty$. Hence $\lim_{r \ti
}P_{n}(\Gamma_{r}) =1$, for any $n \in\N$. Therefore, we can choose
a sequence $r(n) \ti$ such that $\lim_{n \ti} P_{n}(\Gamma_{r(n)})
=1$. Set $ m = n + 3R$. Denote by $\widehat{P}^{(n)}$ the
probability measure under which the particle configurations in the (up
to the boundary, disjoint) boxes $ \Lambda_{m} + 2mk$, with $k \in\Z
^{d}$, are independent and distributed as $P'_n := P_{n}(\cdot\mid
\Gamma_{r(n)})$. In particular, no points are contained in the
corridors $(\Lambda_{m}\setminus\Lambda_{n}) + 2mk$.

We now put
\[
P^{(n)}=\frac1{|\Lambda_{m}|}\int_{\Lambda_{m}}\widehat
P^{(n)}\circ\theta_{z} \,\dd z.
\]
It is then clear that $P^{(n)}\in\Pcal_{\theta}$. A standard
argument shows that $P^{(n)} $ is ergodic (see, e.g., \cite
{G88}, Theorem 14.12). Since $\Gamma_{r(n)}$ is shift invariant and
$\widehat{P}^{(n)}(\Gamma_{r(n)}) =1$, it also follows that
$P^{(n)}(\Gamma_{r(n)}) =1$.
We claim that $\widetilde P=P^{(n)}$ with $n$ sufficiently large,
satisfies the requirements. For this, we have to show that
(1) $\limsup_{n \ti} I_\beta(P^{(n)} ) \le I_\beta(P)$, $ $(2)
$ \limsup_{n \ti} \langle P^{(n)}, \Phi_{\beta} \rangle\le
\langle P, \Phi_{\beta} \rangle$, and finally (3) the net
$(P^{(n)})_{n\in\N}$ converges to $P$ (in the topology $\tau_\Lcal$).

The proof of (1) can be found in the proof of \cite{G94}, Lemma 5.1.

Now we turn to the proof of (2). First note that
%
%
\begin{equation}\label{inspirata1}
\bigl\langle P^{(n)}, \Phi_{\beta} \bigr\rangle= \frac1{|\Lambda_{m}|}
\int_{\Lambda_{m}} \dd z \int\widehat{P}^{(n)} (\dd\omega
) \sum_{x \in\xi\cap(U -z)} \sum_{y \in\xi} T_{x,y}(\omega),
\end{equation}
where we recall the notation in (\ref{Tdef}). The sum on $y$ in (\ref
{inspirata1}) will be split in the sum over $y \in\xi\cap\Lambda
_{n}$ and the remainder. The first sum is handled as follows. As $x,y$
both belong to $\Lambda_{n}$, the measure $\widehat{P}^{(n)}$ can
be replaced by $P_n'$. Furthermore, since $T_{x,y}(\omega) \ge0$, the
integration with respect to $P_n'$ may be estimated against the
integration with respect to $P(\cdot)/P_{n}(\Gamma_{r(n)})$. This gives
\begin{eqnarray*}
&&\frac1{|\Lambda_{m}|} \int_{\Lambda_{m}} \dd z \int\widehat
{P}^{(n)} (\dd\omega) \sum_{x \in\xi\cap(U -z)} \sum_{y \in
\xi\cap\Lambda_{n}} T_{x,y}(\omega)\\
&&\qquad\le\frac1{P_{n}(\Gamma_{r(n)})} \frac1{|\Lambda_{m}|} \int
_{\Lambda_{m}} \dd z \int{P}(\dd\omega) \sum_{x \in\xi\cap(U
-z)} \sum_{y \in\xi}T_{x,y}(\omega).
\end{eqnarray*}
Now use the shift invariance of $P$ and recall that $\lim_{n \ti
}P_{n}(\Gamma_{r(n)}) =1$ to see that the last expression approaches
$\langle P, \Phi_{\beta} \rangle$.

Now we consider the remainder sum in (\ref{inspirata1}), where $y$ is
summed over $\xi\cap\Lambda^{\mathrm{c}}_{m}$. Observe that $|x-y| \ge
3R$, hence we may estimate
\begin{eqnarray*}
T_{x,y}(\omega)&\le&\beta K^{2} \psi_{R}(|x-y|) \le\beta K^{2} \sup
_{x \dvtx|x| \le|z| +1} \psi_{R} (|x-y|)\\
& \le&\beta K^{2} \psi_{R} (|y| - |z| -1),
\end{eqnarray*}
where in the last inequality we used the fact that $|x-y| \ge|x| -
|y|$ and that $\psi_{R}(\cdot) $ is nonincreasing. Now we
distinguish to which of the boxes $\Lambda_{n} + 2km$, with $k \in\Z
^{d}$, the point $y$ belongs (recall that the configurations in these
boxes are independent). Hence for any $z \in\Lambda_{m}$, we have that
\begin{eqnarray*}
&&\int\widehat{P}^{(n)} (\dd\omega) \sum_{x \in\xi\cap(U-z)}
\sum_{y \in\xi\cap\Lambda^{\mathrm{c}}_{m}}T_{x,y}(\omega) \\
&&\qquad\le\beta K^{2} \sum_{k \in\Z^{d}\setminus\{0\}} \int_{\Omega
_{n}} P_n'\bigl(\dd\omega^{(1)}\bigr) \int_{\Omega_{n}} P'_n \bigl(\dd\omega
^{(2)}\bigr) \#\bigl(\xi^{(1)} \cap(U-z)\bigr) \\
&&\qquad\quad\hspace*{54.2pt}{}\times \sum_{y \in(\xi^{(2)} \cap\Lambda_{n}) + 2k m} \psi
_{R}(|y| - |z| -1)\\
&&\qquad\le\frac{\beta K^{2}}{P_{n}(\Gamma_{r(n)})^{2}} \langle P , N_{U}
\rangle\langle P , N_{\Lambda_{n}} \rangle\sum_{k \in\Z^{d}
\setminus\{0\}} \psi_{R} (2|k| m - m - |z| -1),
\end{eqnarray*}
where we estimated integrals with respect to $P_n'$ against integrals
with respect to $ P/P_{n}(\Gamma_{r(n)})$ twice, and used the shift
invariance of $P$. Now we use Assumption~\ref{AssumptionV} and obtain a
constant $C$
(depending only on $R$) such that $\psi_{R}(t)\leq C t^{-h}$ for any
$t\geq0$. Using this in the last display gives that
\begin{eqnarray*}
&&\int\widehat{P}^{(n)}(\dd\omega) \sum_{x \in\xi\cap(U-z)}
\sum_{y \in\xi\cap\Lambda^{\mathrm{c}}_{m}} T_{x,y}(\omega)\\
&&\qquad\le\frac{\beta K^{2}C2^{d}}{P_{n}(\Gamma_{r(n)})^{2}} \langle P ,
N_{U} \rangle^{2} n^{d}\sum_{k \in\Z^{d} \setminus\{0\}} (2|k| m -
m - |z| -1)^{-h}.
\end{eqnarray*}
Now add the factor $1/|\Lambda_{m}|$ and integrate over $z\in\Lambda
_{m}$. Pick some $l=l(n)$ such that $l\sim n$ and $n^d(n-l)^{-h}\to0$
as $n\to\infty$ and split the integral on $z\in\Lambda_{m}$ into
the integrals on $z\in\Lambda_{l}$ and on the remainder. Then it is
easy to see that
\[
\lim_{n \ti} \frac1{|\Lambda_{m}|} \int_{\Lambda_{m}} \dd z
\int\widehat{P}^{(n)} (\dd\omega) \sum_{x \in\xi\cap(U
-z)} \sum_{y \in\xi\cap\Lambda^{\mathrm{c}}_{m}} T_{x,y}(\omega)=0.
\]
Now we have shown (2), that is, that $\limsup_{n \ti} \langle
\widehat{P}^{(n)}, \Phi_{\beta} \rangle\le\langle P, \Phi
_{\beta} \rangle$.

For the proof of (3), we pick $ f\in\Lcal$. Using an affine
transformation, if necessary, we may assume that $ f=f(\cdot\cap
\Delta) $ and $|f|\le N_\Delta$ for some bounded measurable $ \Delta
\subset\R^d $. To estimate the difference of $ |P^{(n)}(f)-P(f)| $ we
integrate over the box $ \Lambda_m $ and get
%
%
\begin{eqnarray}\label{Pnconverges}\quad
&&\bigl|P^{(n)}(f)-P(f)\bigr| \nonumber\\
&&\qquad\le\frac{1}{|\Lambda_m|}\int_{\Lambda_m}\dd x\, \mathbh{1}\{
x+\Delta\subset\Lambda_m\} \bigl|P_n\bigl(f\circ\theta_x\mid\Gamma
_{r(n)}\bigr)-P(f\circ\theta_x)\bigr|\\
&&\qquad\quad{} +\frac{1}{|\Lambda_m|}\int_{\Lambda_m}\dd x\, \mathbh{1}\{
x+\Delta\not\subset\Lambda_m\}\bigl|\widehat
P^{(n)}(N_{\Delta+x})+P(N_{\Delta+x})\bigr|.\nonumber
\end{eqnarray}
Now $ P(N_{\Delta+x})\le\frac{|\Delta|\mu(P)}{P_n(\Gamma_{r(n)})}
$, where $ \mu(P) <\infty$ is the intensity of $P$. In the same way
we obtain
\[
\widehat P^{(n)}(N_{\Delta+x})=P_n\bigl(N_{\Delta+x \operatorname{mod}
2m+1}\mid\Gamma_{r(n)}\bigr)\le\frac{|\Delta|\mu(P)}{P_n(\Gamma_{r(n)})}.
\]
Hence the second term on the right-hand side of (\ref{Pnconverges}) is
not larger than the volume of $\{x\in\Lambda_m\dvtx x+\Delta\not
\subset\Lambda_m\}$ (which is of surface order of $\Lambda_m$) times
$ O(|\Lambda_m|^{-1}) $, that is, it vanishes. Concerning the first
term on the right-hand side of (\ref{Pnconverges}), we estimate
\begin{eqnarray*}
&&\bigl|P_n\bigl(f\circ\theta_x\mid\Gamma_{r(n)}\bigr)-P(f\circ\theta_x)
\bigr|\\
&&\qquad\leq\biggl|\frac{1}{P_n(\Gamma_{r(n)})}-1\biggr| P_n\bigl(N_{\Delta+x};
\Gamma_{r(n)}\bigr)+P_n\bigl(N_{\Delta+x};\Gamma_{r(n)}^{\mathrm{c}}\bigr)\\
&&\qquad\le|\Delta| \mu(P)\biggl|\frac{1}{P_n(\Gamma_{r(n)})}-1
\biggr|+P(N_\Delta^2)^{1/2}\bigl(1-P_n\bigl(\Gamma_{r(n)}\bigr)\bigr)^{1/2}.
\end{eqnarray*}
By Lemma \ref{lem-Temp}, $P(N_\Delta^2)$ is finite, hence the
right-hand side vanishes as $n\to\infty$. Therefore, also the first
term on the right-hand side of (\ref{Pnconverges}) vanishes, and we
conclude that (3) holds.
\end{pf}

\subsubsection*{Final step: Proof of the lower bound in (\protect\ref{freeenlower})}

Now we can finish the proof of the lower bound in (\ref
{freeenlower}). Recall that it is sufficient to prove (\ref{mainaim})
for any $\beta,\rho\in(0,\infty)$, to get both lower bounds in
(\ref{freeenlower}). Fix $K, R\in\N$ and $\delta\in(0,\rho) $. We
start from the right-hand side of (\ref{reduc1}). Fix $\eps>0 $, and
pick $P\in\Pcal_\theta(\Omega^{(K,R)})$ satisfying $I_{\beta
}^{(K,R)}(P)+\langle P,\Phi_\beta\rangle<\infty$ and $
|\langle P,N_U^{(\ell)}\rangle-\rho|<\delta$. By Lemma \ref
{lem-ergodappr}, we may fix some $r>0$ and some ergodic measure
$\widetilde P\in\Pcal_\theta(\Omega^{(K,R)})$ satisfying $
|\langle P,N_U^{(\ell)}\rangle-\rho|<\delta$ and $\langle
\widetilde P,\Phi_\beta\rangle\leq\langle P,\Phi_\beta\rangle
+\eps$ and $I_\beta^{(K,R)}(\widetilde P)\le I_\beta^{(K,R)}(P)+\eps
$ and $\widetilde P(\Gamma_r)=1$. Since $I_{\beta
}^{(K,R)}(\widetilde P)<\infty$, for $N$ large enough there is a
density $f_{N}^{(K,R)}$ of the projection $\widetilde P_{L_N}$ of
$\widetilde P$ to $\Omega_{L_N}^{(K,R)}$ with respect to the
projection $\ttQ _{L_N}^{(K,R)}$ of the restricted marked
Poisson point process $\ttQ ^{(K,R)}$ to $\Omega_{L_N}$, where
we recall that $\Omega_{L_N}$ is the set of restrictions of
configurations in $\Omega$ to $\Lambda_{L_N}$, and $\Omega
_{L_N}^{(K,R)}$ is defined analogously.
We conceive $ \RR_{N}$ as a map $\RR_{N,\cdot}\dvtx\Omega_{L_N}\to
\Pcal_\theta(\Omega^{(K,R)}) $. Now introduce the event
%
%
\begin{eqnarray}
C_N&=&\biggl\{\omega\in\Omega_{L_N}^{(K,R)}\dvtx\langle\RR
_{N,\omega},\Phi_\beta\rangle\le\langle\widetilde P,\Phi_\beta
\rangle+\eps,\nonumber\\[-8pt]\\[-8pt]
&&\hspace*{15.1pt} \frac1{|\Lambda_{L_N}|}\log f_{N}^{(K,R)}(\omega)\le
I_{\beta}^{(K,R)}(\widetilde
P)+\eps\biggr\}.\nonumber
\end{eqnarray}
Then we can estimate
%
%
\begin{eqnarray}\label{Zlowbound5}\quad
&&\ttE^{(K,R)}\bigl[ \ee^{-|\Lambda_N| \langle\RR
_{N},\Phi
_\beta\rangle}\mathbh{1}\bigl\{\bigl|\bigl\langle\RR_{N},N_U^{(\ell)}\bigr\rangle
-\rho\bigr|<\delta\bigr\} \bigr]\nonumber\\
&&\qquad=\int_{\Omega_{L_N}^{(K,R)}} \dd\ttQ _{L_N}^{(K,R)}
\ee^{-|\Lambda_{N}|\langle\RR_{N},\Phi_\beta\rangle}\mathbh
{1}\bigl\{\bigl|\bigl\langle\RR_{N},N_U^{(\ell)}\bigr\rangle-\rho\bigr|<\delta\bigr\}
\nonumber\\
&&\qquad\geq\int_{C_N}\widetilde P_{L_N}(\dd\omega) \frac
1{f_N^{(K,R)}(\omega)}\ee^{-|\Lambda_{N}| \langle\RR
_{N},\Phi_\beta
\rangle}\mathbh{1}\bigl\{\bigl|\bigl\langle\RR_{N},N_U^{(\ell)}\bigr\rangle-\rho
\bigr|<\delta\bigr\}\\
&&\qquad\ge\ee^{-|\Lambda_{L_N}|(I_\beta^{(K,R)}(\widetilde
P)+\eps
)} \ee^{-|\Lambda_{L_N}|(\langle\widetilde P,\Phi
_\beta\rangle
+\eps)}\nonumber\\
&&\qquad\quad{} \times\widetilde P_{L_N}\bigl(C_N\cap\bigl\{\omega\in\Omega
_{L_N}^{(K,R)}\dvtx\bigl|\bigl\langle
\RR_{N},N_U^{(\ell)}\bigr\rangle-\rho\bigr|<\delta\bigr\}\bigr).\nonumber
\end{eqnarray}
The continuity of the map $ P\mapsto\langle P,\Phi_\beta\rangle$
(see Lemma \ref{lem-Fcont}),
the law of large numbers and McMillan's theorem imply that
\begin{eqnarray*}
&&\widetilde P_{L_N}\biggl(\biggl\{\omega\in\Omega_{L_N}^{(K,R)}\dvtx
\bigl|\bigl\langle\RR_{N,\omega},N_U^{(\ell)}\bigr\rangle-\rho\bigr|<\delta,
\langle\RR_{N,\omega},\Phi_\beta\rangle\le\langle\widetilde
P,\Phi_\beta\rangle+\eps,\\
&&\qquad\hspace*{124.4pt}\frac{1}{|\Lambda_{L_N}|}\log f_N^{(K,R)}(\omega)\le I_\beta
^{(K,R)}(\widetilde P)+\eps\biggr\}\biggr) \to1,
\end{eqnarray*}
as $N\to\infty$. Using this in (\ref{Zlowbound5}) and this in (\ref
{reduc1}), we arrive at
%
%
\begin{equation}\label{Zlowbound6}
\liminf_{N\to\infty}\frac1{|\Lambda_{L_N}|}\log Z_N(\beta,\Lambda_{L_N})
\geq\overline q-\delta- I_\beta^{(K,R)}(\widetilde P)-\eps
-\langle\widetilde P,\Phi_\beta\rangle-\eps.\hspace*{-34pt}
\end{equation}
Now recall that $\langle\widetilde P,\Phi_\beta\rangle\le\langle
P,\Phi_\beta\rangle+\eps$ and $I_\beta^{(K,R)}(\widetilde
P)\le I_\beta(P)+\eps$. Now we can let $\eps\downarrow0$ and take
the infimum over $P$, to obtain
\begin{eqnarray*}
&&\liminf_{N\to\infty}\frac1{|\Lambda_{L_N}|}\log Z_N(\beta
,\Lambda_{L_N})\\
&&\qquad\geq\overline q -\delta-\inf_{P\in\Pcal_\theta(\Omega^{(K,R)})
\dvtx|\langle P, N^{(\ell)}_{U}\rangle-\rho|<\delta
}\bigl\{ I_{\beta}^{(K,R)}(P) +\langle P, \Phi_\beta\rangle
\bigr\}.
\end{eqnarray*}

Our last step is to approach the variational formula $\chi^{
(=)}(\beta,\rho)$ on the right-hand side of (\ref{freeenlower}) by the
finite-$K$ and finite-$R$ versions.
\begin{lemma}[(Removing the cut-off)]\label{lem-Kapprox}For
any $\delta\in(0,\rho)$,
%
%
\begin{eqnarray}\label{Kapproxvarform}
&&\limsup_{K, R \to\infty}\inf_{P\in\Pcal_\theta(\Omega^{(K,R)})
\dvtx|\langle P, N^{(\ell)}_{U}\rangle-\rho|<\delta
}\bigl\{ I_{\beta}^{(K,R)}(P) +\langle P, \Phi_\beta\rangle
\bigr\}\nonumber\\[-8pt]\\[-8pt]
&&\qquad\le\inf_{P\in\Pcal_\theta(\Omega) \dvtx\langle P, N^{(\ell
)}_{U}\rangle= \rho}\{ I_{\beta}(P) + \langle P, \Phi_\beta
\rangle\}= \chi^{(=)}(\beta,\rho).\nonumber
\end{eqnarray}
\end{lemma}
\begin{pf}
Fix $P\in\Pcal_\theta$ satisfying $\langle P, N_{U}^{(\ell)}
\rangle= \rho$ and $I_\beta(P) + \Phi_\beta(P)<\infty$. For $K,R
\in\N$, consider $P_{K,R}=P\circ\pi_{K,R}^{-1}$. Then we have $
P_{K,R}(\Omega^{(K,R)})=1 $ and $ \langle P_{K,R},N_U^{(\ell)}\rangle
=\langle P,\pi_{K,R}\circ N_U^{(\ell)}\rangle
\uparrow\langle P,N_U^{(\ell)}\rangle$ for $ K,R \to\infty$ by
the monotonous convergence theorem. Hence, for $K$ and $R$ sufficiently
large, $ |\langle P_{K,R},N_U^{(\ell)}\rangle-\rho|<\delta$.
Observe that $ \langle P_{K,R},\Phi_\beta\rangle\le\langle P,\Phi
_\beta\rangle$ since $ \Phi_\beta\ge0 $. By (\ref{IKident}), we
have $ I_\beta^{(K,R)}(P_{K,R})\le I_\beta(P) $. Finally,
observe that the infimum over $P$ such that $|\langle P, N^{(\ell
)}_{U}\rangle-\rho|<\delta$ is obviously not larger than the infimum
over $P$ satisfying $\langle P, N^{(\ell)}_{U}\rangle= \rho$.
\end{pf}

\subsection{\texorpdfstring{Proof of Theorem \protect\ref{thm-mainres} for Dirichlet
and periodic boundary conditions}%
{Proof of Theorem 1.2 for Dirichlet and periodic boundary conditions}}

\label{otherbc}

In this section, we show how to adapt the proof of
Theorem \ref{thm-mainres} for empty boundary conditions to obtain the
proof for Dirichlet and periodic boundary conditions. Let us make a
couple of obvious observations. First, the restriction of the
periodized Brownian bridge measure on paths that do not leave the box
$\Lambda$ equals the Brownian bridge measure with Dirichlet boundary
conditions, that is,
\[
\mu_{x,x}^{({\romper,k\beta})}|_{\Ccal_{k,\Lambda}^{
(\romDir)}}=\mu_{x,x}^{({\romDir,k\beta})}.
\]
Hence, it is easy to see that $ \overline{q}^{({\romDir})}\le
\overline{q}^{({\romper})} $ and that
%
%
\begin{equation}\label{Zbccompare}
Z_N^{(\romDir)}(\beta,\Lambda)\le Z_N(\beta,\Lambda)\le
Z_N^{(\romper)}(\beta,\Lambda),
\end{equation}
since the Feynman--Kac formula for $Z_N^{(\romDir)}$ contains
only those paths that stay in $\Lambda$ all the time with the same
distribution as under which they appear in the formula for $Z_N^{
(\romper)}$. Hence, it will be sufficient to prove the upper bound in
(\ref{freeenupper}) for $Z_N^{(\romper)}$ and the lower bound in
(\ref{freeenlower}) for $Z_N^{(\romDir)}$ only.

We start with the representation of $Z_N^{(\romDir)}$ and
$Z_N^{(\romper)}$ given in Proposition \ref{lem-rewrite}. The
first step is to show that the weights $ \overline{q}^{(\rombc)}
$ converge to $\overline{q}=\sum_{k\in\N}q_k $. For notational
reasons, we now write $\overline{q}^{(\rombc)}_{\Lambda}$ for
$\overline{q}^{(\rombc)}$; however notice that it depends on
$N$. Recall that $\Lambda_{L_N}$ is the centered box with side length
$L_N=(N/\rho)^{1/d}$.
\begin{lemma}
Let $\rombc\in\{\romDir, \romper\} $. Then
%
%
\begin{equation}
\lim_{N\to\infty}\overline{q}^{(\rombc)}_{\Lambda
_{L_N}}=\overline{q}.
\end{equation}
\end{lemma}
\begin{pf}
(a) First we consider periodic boundary conditions. Then we have
%
%
\begin{equation}\label{perbc1}
\overline{q}^{(\romper)}_{\Lambda_{L_N}}=(4\pi\beta
)^{-d/2}\sum_{k=1}^N\frac{1}{k^{1+d/2}}\sum_{z\in\Z^d}\ee^{-
{|z|^2}/({4k\beta})L_N^2}.
\end{equation}
Since\vspace*{-1pt} the summand for $z=0$ converges toward $(4\pi\beta)^{-d/2}\sum
_{k=1}^\infty\frac{1}{k^{1+d/2}}=\overline q$, we only\vspace*{1pt} have to show
that $\sum_{k=1}^N\frac{1}{k^{1+d/2}}\sum_{z\in\Z^d\setminus\{0\}
}\ee^{-{|z|^2}/({4k\beta})L_N^2}$ vanishes as $N\to\infty$.

Using an approximation with an integral, one sees that, for some $c\in
(0,\infty)$, only depending on $d$,
\[
\sum_{z\in\Z^d\setminus\{0\}}\ee^{-a|z|^2}\le c a^{-d/2}
\qquad\mbox{for all } a\in(0,\infty).
\]
Using this with $ a=L_N^2/(4\beta k) $, we see that $\sum_{z\in\Z
^d\setminus\{0\}}\ee^{-{|z|^2}/({4k\beta})L_N^2}$ is of order
$k^{d/2}L_N^{-d}$. Using that $N$ is of order $L_N^d$ and applying the
harmonic series, we see that $\sum_{k=1}^N\frac{1}{k^{1+d/2}}\sum
_{z\in\Z^d\setminus\{0\}}\ee^{-{|z|^2}/({4k\beta})L_N^2}$ is
of order $L_N^{-d}\log L_N$ and therefore vanishes as $N\to\infty$.

(b) Now we consider Dirichlet boundary conditions. For any $
M\in\N$ and $ \delta\in(0,1)$, we get, for any sufficiently large $N$,
%
%
\begin{eqnarray}\label{bcdir1}
\overline{q}_{\Lambda_{L_N}}^{({\romDir})}
&=&\frac{1}{|\Lambda_{L_N}|}\sum_{k=1}^N\frac{1}{k}\int_{\Lambda
_{L_N}}\dd x\, \mu_{x,x}^{(k\beta)}\bigl(
B_{[0,k\beta]}\subset\Lambda_{L_N}\bigr)\nonumber\\[-8pt]\\[-8pt]
&\ge&\sum_{k=1}^M\frac{1}{k}\frac{1}{|\Lambda_{L_N}|} \int
_{(1-\delta)\Lambda_{L_N}}\dd x\,
\mu_{x,x}^{(k\beta)}\bigl(B_{[0,k\beta]}\subset\Lambda
_{L_N}\bigr).\nonumber
\end{eqnarray}
It is easy to see that, in the limit $N\to\infty$, the integrand $
\mu_{x,x}^{(k\beta)}(B_{[0,k\beta]}\subset\Lambda_{L_N})$
tends to $\mu_{0,0}^{(k\beta)}(\mathbh{1})= (4\pi k\beta
)^{-d/2}$, uniformly in $x\in(1-\delta)\Lambda_{L_N}$ and $k\in\{
1,\ldots,M \}$. Hence,
\[
\liminf_{N\to\infty}
\overline{q}_{\Lambda_{L_N}}^{({\romDir})}\geq\sum
_{k=1}^M\frac{1}{k} (4\pi k\beta)^{-d/2}\frac{|(1-\delta)\Lambda
_{L_N}|}{|\Lambda_{L_N}|},
\]
which tends to $\overline q$ as $M\to\infty$ and $\delta\downarrow0$.
\end{pf}

\subsubsection*{Proof of the upper bound for periodic boundary condition}

We continue to write $\Lambda$ for $\Lambda_{L_N}$, where
$L_N=(N/\rho)^{1/d}$. We adapt the proof of the upper bound in
Section \ref{sec-upbound} for periodic boundary conditions. The main
idea is to drop all the paths that reach the boundary of the box
$\Lambda$ and to use that their distribution is equal to the one under
the free Brownian bridge measure. Let us introduce, for parameters
$r\in(0,1)$ and $\widetilde R\in(0,\infty)$, the random variable
%
%
\begin{equation}
N^{({\ell,\widetilde R})}_{r\Lambda}(\omega) =\sum_{x\in\xi
\cap r\Lambda}\ell(f_x)\mathbh{1}\Bigl\{{\sup_{s\in[0,\beta\ell
(f_x)]}}|f_x(s)-f_x(0)|\leq\widetilde R\Bigr\},
\end{equation}
the total length of the marks of particles starting in $r\Lambda$ that
stay within distance $\leq\widetilde R$ from their starting sites.
Furthermore, let
\begin{eqnarray*}
H_{r\Lambda}^{({\widetilde R})}(\omega)&=&\sum_{x,y\in\xi\cap
r\Lambda}T_{x,y}(\omega)\mathbh{1}\Bigl\{{\sup_{s\in[0,\beta\ell
(f_x)]}}|f_x(s)-f_x(0)|\leq\widetilde R\Bigr\}\\
&&\hspace*{36.8pt}{}\times\mathbh{1}\Bigl\{{\sup_{s\in[0,\beta\ell
(f_y)]}}|f_y(s)-f_y(0)|\leq\widetilde R\Bigr\},
\end{eqnarray*}
be the Hamiltonian in (\ref{Hamiltonian}) restricted to paths starting
in $ r\Lambda$ and traveling no further than $\widetilde R$. Note
that, for $N$ large enough (depending only on $r$ and $\widetilde R$),
such paths will never reach the boundary of $\Lambda$ and therefore
have the same distribution under the periodized Brownian bridge measure
as under the free one or the one with Dirichlet boundary condition.
Hence, we estimate
%
%
\begin{eqnarray}\label{periodicbcupper}
&&
\ttE^{(\romper)}\bigl[\ee^{-H_\Lambda(\omega_{\romP})}\mathbh{1}\bigl\{
N_\Lambda^{({\ell})}(\omega_{\romP})=N\bigr\}\bigr]\nonumber\\
&&\qquad\le\ttE^{(\romper)}\bigl[\ee^{-H^{({\widetilde
R})}_{r\Lambda}(\omega_{\romP})}\mathbh{1}\bigl\{N_{r\Lambda}^{
({\ell,\widetilde R})}(\omega_{\romP})\le
N\bigr\}\bigr]\nonumber\\[-8pt]\\[-8pt]
&&\qquad=\ttE^{(\romDir)}\bigl[\ee^{-H^{({\widetilde
R})}_{r\Lambda}(\omega_{\romP})}\mathbh{1}\bigl\{N_{r\Lambda}^{
({\ell,\widetilde R})}(\omega_{\romP})\le N\bigr\}\bigr]\nonumber\\
&&\qquad\le\ttE\bigl[\ee^{-H_{r\Lambda}^{({\widetilde R})}(\omega
_{\romP})}\mathbh{1}\bigl\{N_{r\Lambda}^{({\ell,\widetilde
R})}(\omega_{\romP})\le
N\bigr\}\bigr],\nonumber
\end{eqnarray}
where ``$(\romper)$'' and ``$(\romDir)$'' refer to the box $\Lambda$.
Therefore, we can use the same method as in Section \ref{sec-upbound},
the only\vspace*{2pt} two differences being that $\Lambda$ is replaced by $r\Lambda
$ and that we deal solely with paths that do not travel further than
$\widetilde R$. That is, we have two additional truncation parameters
$r$ and $\widetilde R$. It is straightforward to see that adapted
versions of Lemmas \ref{lem-totmark} and \ref{lem-Hloweresti} hold
and that the proof\vspace*{1pt} given in Section~\ref{sec-upbound} applies verbatim
as well. Finally, one takes the limits $\widetilde R\to\infty$ and
$r\uparrow1$ in the resulting variational formula, which is the same
as the proof of (\ref{Slimit}).

\subsubsection*{Proof of the lower bound for Dirichlet boundary conditions}

We continue to write $\Lambda$ for $\Lambda_{L_N}$, where
$L_N=(N/\rho)^{1/d}$. The strategy for Dirichlet boundary conditions
is as follows. First we pick some $ \eps\in(0,\frac12) $ and
consider $ \widetilde\Lambda=(1-\eps)\Lambda$ and $ \partial
\Lambda=\Lambda\setminus\widetilde\Lambda$. The idea is to require
that $\partial\Lambda$ receives no particle and that the marks of all
particles in $\widetilde\Lambda$ have length $\leq K$ and spatial
extension $\leq R$. In\vspace*{1pt} this way, we get a lower estimate against the
truncated version of the Poisson process on $\widetilde\Lambda$
rather than on $L$. The only difference to the proof for empty boundary
condition is then that Lemma \ref{lem-makeopen}, which was given
before the introduction of the truncation, now has to be proved with
the presence of the truncation, which requires some adaptation. Every
other step of the proof is literally the same for $\Lambda$ instead of
$\widetilde\Lambda$, which means that in the end of the proof, the
parameter $\eps$ has to be sent to $0$, which is extremely simple.

Let us come to the details. We first show that there exist $c>0$ and
$C_{K,R}>0$ such that, for any $N, R,K\in\N$,
%
%
\begin{eqnarray}\label{fact1lower}
&&\ttE^{(\romDir)}\bigl[\ee^{-H_\Lambda(\omega_{\romP})}\mathbh{1}\bigl\{
N_\Lambda^{({\ell})}(\omega_{\romP})=N\bigr\}\bigr]
\nonumber\\[-8pt]\\[-8pt]
&&\qquad\ge\ee^{-\eps c|\Lambda|}\ee^{-C_{K,R}|\Lambda|} \ttE
^{(K,R)}\bigl[\ee^{-H_{\widetilde\Lambda}(\omega_{\romP})}\mathbh{1}\bigl\{
N_{\widetilde\Lambda}^{({\ell})}(\omega_{\romP})=N\bigr\}\bigr],\nonumber
\end{eqnarray}
where $ C_{K,R}\to0 $ as $ R\to\infty$ and afterward $K\to\infty
$. This is done as follows. Estimate
%
%
\begin{eqnarray}\label{lowerestDir}
&&\ttE^{(\romDir)}\bigl[\ee^{-H_\Lambda(\omega_{\romP})}\mathbh{1}\bigl\{
N_\Lambda^{({\ell})}(\omega_{\romP})=N\bigr\}
\bigr]\nonumber\\
&&\qquad=\ttE \bigl[\ee^{-H_\Lambda(\omega_{\romP})}\mathbh{1}\bigl\{
N_\Lambda^{({\ell})}(\omega_{\romP})=N\bigr\}\nonumber\\
&&\hspace*{7.1pt}\qquad\quad{} \times \mathbh{1}\{\forall
x\in\xi_{\romP}\cap\Lambda\dvtx
B_x([0,\beta\ell(B_x)])\subset\Lambda\}\bigr]\nonumber\\
&&\qquad\ge\ttE \bigl[\ee^{-H_\Lambda(\omega_{\romP})}\mathbh{1}\bigl\{
N_\Lambda^{(\ell)}(\omega_{\romP})=N\bigr\}\mathbh{1}\bigl\{\forall
x\in\xi_{\romP}\cap\widetilde\Lambda\dvtx B_x\in E^{(K,R)}\bigr\}
\\
&&\hspace*{7.1pt}\qquad\quad{} \times\mathbh{1}\{\forall x\in\xi_{\romP}\cap\Lambda\dvtx
B_x([0,\beta\ell(B_x)])\subset\Lambda\}\mathbh{1}\{N_{\partial
\Lambda}(\omega_{\romP})=0\}\bigr] \nonumber\\
&&\qquad=\ttE\bigl[\ee^{-H_{\widetilde\Lambda}(\omega_{\romP})}\mathbh{1}\bigl\{
N_{\widetilde\Lambda}^{(\ell)}(\omega_{\romP})=N\bigr\}\mathbh{1}\{
N_{\partial\Lambda}(\omega_{\romP})=0\}\nonumber
\\
&&\hspace*{70.2pt}\qquad\quad{} \times\mathbh{1}\bigl\{\omega_{\romP}\bigl(\widetilde\Lambda
\times\bigl(E^{(K,R)}\bigr)^{
\mathrm{c}}\bigr)=0\bigr\}\bigr].\nonumber
\end{eqnarray}
Independence of the events in the indicators gives
%
%
\begin{eqnarray}
&&\mbox{right-hand side of (\ref{lowerestDir}) }\nonumber\\
&&\qquad=\ttE^{(K,R)}
\bigl[\ee^{-H_{\widetilde\Lambda}(\omega_{\romP})}\mathbh{1}\bigl\{
N_{\widetilde\Lambda}^{(\ell)}(\omega_{\romP})=N\bigr\}\bigr]\ttQ
\bigl(N_{\partial\Lambda}(\omega_{\romP})=0\bigr)\nonumber\\[-8pt]\\[-8pt]
&&\qquad\quad{} \times\ttQ \bigl(\omega_{\romP}\bigl(\widetilde\Lambda\times
\bigl(E^{(K,R)}\bigr)^{\mathrm{c}}\bigr)=0\bigr)\nonumber\\
&&\qquad=\ttE^{(K,R)} \bigl[\ee^{-H_{\widetilde\Lambda}(\omega
_{\romP})}\mathbh{1}\bigl\{N_{\widetilde\Lambda}^{(\ell)}(\omega
_{\romP})=N\bigr\}\bigr]\ee^{-\overline{q}|\partial\Lambda|}\ee^{-\nu(\widetilde
\Lambda\times(E^{(K,R)})^{
\mathrm{c}})},\nonumber
\end{eqnarray}
since $N_{\widetilde\Lambda}(\omega_{\romP})$ and $\omega_{\romP
}(\widetilde\Lambda\times(E^{(K,R)})^{\mathrm{c}})$ are
Poisson distributed with respective parameters $\overline{q}|\partial
\Lambda|$ and $\nu(\widetilde\Lambda\times(E^{(K,R)})^{\mathrm{c}})$.
We estimate $\overline{q}|\partial\Lambda|\leq c \eps|\Lambda|$
for some $c>0$ and
%
%
\begin{eqnarray}
&&\nu\bigl(\widetilde\Lambda\times\bigl(E^{(K,R)}\bigr)^{
\mathrm{c}}\bigr)\nonumber\\
&&\qquad\leq|\widetilde\Lambda|\sum_{k=K+1}^\infty\frac{q_k}{k}
+|\widetilde\Lambda|\sum_{k=1}^K\mu_{0,0}^{(k\beta)}
\Bigl({\max_{s\in[0, \beta k]}}|B_s|> R\Bigr)\\
&&\qquad\leq|\Lambda| C_{K,R},\nonumber
\end{eqnarray}
with some $C_{K,R}$ that vanishes as $R\to\infty$ and afterward
$K\to\infty$. Hence, we have got (\ref{fact1lower}).

Now we need a version of Lemma \ref{lem-makeopen} for truncated point
processes, that is, we need to show that, for any $R,K\in\N$ and for
any $\delta\in(0,\rho)$, for all sufficiently large~$N$,
%
%
\begin{eqnarray}\label{fact2lower}\qquad
&&\ttE^{(K,R)}\bigl[ \ee^{-H_{\Lambda}(\omega_{\romP})}
\mathbh{1}\bigl\{\bigl\langle\RR_{N}(\omega_{\romP}), N^{(\ell)}_{U}\bigr\rangle
= \rho\bigr\} \bigr]\nonumber\\
&&\qquad\ge\frac{(C_1\wedge C_2)^{\delta|\Lambda|}}{2\delta|\Lambda|+2}
\nonumber\\[-8pt]\\[-8pt]
&&\qquad\quad{}\times\ttE^{(K,R)} \bigl[ \ee^{-|\Lambda| \langle
\RR
_{N}(\omega_{\romP}),\Phi_\beta\rangle} \nonumber\\
&&\qquad\quad\hspace*{43.2pt}{}\times\mathbh{1}\bigl\{\bigl\langle\RR
_{N}(\omega_{\romP}), N^{(\ell)}_{U}\bigr\rangle\in
(\rho-\delta,\rho+\delta)\bigr\}\bigr],\nonumber
\end{eqnarray}
where $C_1$ and $C_2$ may depend on $R$ and $K$.

Since Lemma \ref{lem-Zquotient} was used in the proof of Lemma \ref
{lem-makeopen}, we first need a truncated version of Lemma \ref
{lem-Zquotient}. For this we consider the truncated version of
$Z_N(\beta,\Lambda)$,
%
%
\begin{eqnarray}
&&Z_N^{(K,R)}(\beta,\Lambda)\nonumber\\
&&\qquad=\sum_{\lambda\in\mathfrak
P_N\dvtx\sum_{k=1}^Kk\lambda_k=N} \prod_{k=1}^K\frac{(\overline
{q}^{(R)}_{k,\Lambda})^{\lambda_k}|\Lambda|^{\lambda
_k}}{\lambda_k !k^{\lambda_k}}\\
&&\qquad\quad\hspace*{93.2pt}\times{}\bigotimes_{k=1}^K\bigl(\E_\Lambda
^{(R,k\beta)}\bigr)^{\otimes\lambda_k}[\ee^{-\Gcal
_{N,\beta}^{(\lambda)}}],\nonumber
\end{eqnarray}
where
\[
\overline{q}^{(R)}_{k,\Lambda}=\frac{1}{|\Lambda|}\int
_\Lambda\dd x\, \mu_{x,x}^{(k\beta)}\Bigl({\max_{s\in[0,
\beta k]}}|B_s-B_0|\le R\Bigr),
\]
and where $ \E_\Lambda^{(R,k\beta)} $ is the expectation with
respect to the probability measure
\[
\PP_\Lambda^{(R,k\beta)}(\dd f)=\frac{\int_\Lambda\dd
x\,\mu_{x,x}^{(k\beta)}(\dd f\mathbh{1}\{{\max_{ s\in[0,
\beta k]}}|f_s-f_0|\le R\})}{|\Lambda|\overline{q}_\Lambda
^{(R)}}.
\]
All steps in the proof of Lemma \ref{lem-Zquotient} are easily
adapted, but the estimate in (\ref{energyesti2}) needs a slightly
different argument. We now estimate
\begin{eqnarray*}
\\[-17pt]
&&\E_\Lambda^{(R,\beta)}\bigl(v\bigl(|B_s-f(s)|\bigr)\bigr)\\[-2.5pt]
&&\qquad=\frac{1}{\overline{q}^{(R)}_\Lambda|\Lambda|}\int_\Lambda
\dd x \,\E_x\Bigl[v\bigl(|B_s-f(s)|\bigr)
\mathbh{1}\Bigl\{{\max_{0\le s\le\beta}}|B_s-B_0|\le R\Bigr\}
,B_\beta\in\dd x\Bigr]\big/\dd x\\[-2.5pt]
&&\qquad \le\frac{(4\pi\beta)^{-d/2}}{\overline{q}{}^{(R)}_\Lambda
|\Lambda|}\int_\Lambda\dd x\int_\Lambda\dd y\, \frac
{g_s(x,y)v(|y-f(s)|)g_{\beta-s}(y,x)}{g_\beta(x,x)}.
\\[-11pt]
\end{eqnarray*}
Now we can proceed as in (\ref{energyesti}), (\ref{energyesti2}) and
obtain that $\E_\Lambda^{(R,\beta)}(v(|B_s-f(s)|)\le\frac
{\alpha(v)(4\pi\beta)^{-d/2}}{\overline{q}^{(R)}_\Lambda
|\Lambda|}$.
Hence, we get the following truncated version of Lemma \ref{lem-Zquotient}:
%
%
\begin{eqnarray}\label{lower-Dir}
\nonumber\\[-12.5pt]
\frac{Z_{N+1}^{(K,R)}(\beta,\Lambda)}{Z_N^{(K,R)}(\beta
,\Lambda)}&\ge&\frac{|\Lambda|}{N+1} \exp\biggl(-\frac{N\beta\alpha
(v)(4\pi\beta)^{-d/2}}{|\Lambda|\overline{q}^{(R)}_\Lambda
}\biggr).\\[-11pt]
\nonumber
\end{eqnarray}
Using this instead of Lemma \ref{lem-Zquotient} in the proof of
Lemma \ref{lem-makeopen}, we get the truncated version (\ref
{fact2lower}) of Lemma \ref{lem-makeopen} with $C_2$ as before and with
$C_1$ replaced by
\begin{eqnarray*}
\\[-15.5pt]
C^{(R)}_1&=&1\wedge\frac{\overline{q}^{(R)}_\Lambda}{\rho
+\delta}\exp\biggl(-\frac{(\rho+\delta)\beta\alpha(v)(4\pi\beta
)^{-d/2}}{\overline{q}_\Lambda^{(R)}}\biggr).
\\[-11.5pt]
\end{eqnarray*}

The remaining proof of the lower bound is exactly as in the case of
empty boundary condition, with $\widetilde\Lambda$ instead of
$\Lambda$. This slight difference vanishes in the end when taking
$\eps\downarrow0$.\vspace*{-1.5pt}

\section*{Acknowledgment}
We thank an anonymous referee whose detailed comments helped us to fix
two technical points in the proofs.\vspace*{-1.5pt}

%

%
\printaddresses


\begin{thebibliography}{23}

\bibitem{A07}
%
\begin{bmisc}[vtex]
\bauthor{\bsnm{Adams},~\bfnm{Stefan}\binits{S.}}
(\byear{2009}).
\bhowpublished{Large deviations for empirical path
measures in cycles of integer partitions. Preprint}.
\end{bmisc}
%
\endbibitem

\bibitem{AD06}
%
\begin{barticle}[vtex]
\bauthor{\bsnm{Adams},~\bfnm{Stefan}\binits{S.}} \AND
\bauthor{\bsnm{Dorlas},~\bfnm{Tony}\binits{T.}}
(\byear{2008}).
\btitle{Asymptotic {F}eynman--{K}ac formulae for large symmetrised
systems of
random walks}.
\bjournal{Ann. Inst. H. Poincar\'e Probab. Statist.}
\bvolume{44}
\bpages{837--875}.
\bid{doi={10.1214/07-AIHP132}, mr={2453847}}
\end{barticle}
%
\endbibitem

\bibitem{AK08}
%
\begin{barticle}[mr]
\bauthor{\bsnm{Adams},~\bfnm{Stefan}\binits{S.}} \AND
\bauthor{\bsnm{K{\"o}nig},~\bfnm{Wolfgang}\binits{W.}}
(\byear{2008}).
\btitle{Large deviations for many {B}rownian bridges with symmetrised
initial-terminal condition}.
\bjournal{Probab. Theory Related Fields}
\bvolume{142}
\bpages{79--124}.
\bid{doi={10.1007/s00440-007-0099-5}, mr={2413267}}
\end{barticle}
%
\endbibitem

\bibitem{BCMP05}
%
\begin{barticle}[mr]
\bauthor{\bsnm{Benfatto},~\bfnm{Giuseppe}\binits{G.}},
\bauthor{\bsnm{Cassandro},~\bfnm{Marzio}\binits{M.}},
\bauthor{\bsnm{Merola},~\bfnm{I.}\binits{I.}} \AND
\bauthor{\bsnm{Presutti},~\bfnm{E.}\binits{E.}}
(\byear{2005}).
\btitle{Limit theorems for statistics of combinatorial partitions with
applications to mean field {B}ose gas}.
\bjournal{J. Math. Phys.}
\bvolume{46}
\bpages{033303, 38}.
\bid{doi={10.1063/1.1855933}, mr={2125575}}
\end{barticle}
%
\endbibitem

\bibitem{BUe09}
%
\begin{barticle}[mr]
\bauthor{\bsnm{Betz},~\bfnm{Volker}\binits{V.}} \AND
\bauthor{\bsnm{Ueltschi},~\bfnm{Daniel}\binits{D.}}
(\byear{2009}).
\btitle{Spatial random permutations and infinite cycles}.
\bjournal{Comm. Math. Phys.}
\bvolume{285}
\bpages{469--501}.
\bid{doi={10.1007/s00220-008-0584-4}, mr={2461985}}
\end{barticle}
%
\endbibitem

\bibitem{BR97}
%
\begin{bbook}[vtex]
\bauthor{\bsnm{Bratteli},~\bfnm{Ola}\binits{O.}} \AND
\bauthor{\bsnm{Robinson},~\bfnm{Derek~W.}\binits{D.~W.}}
(\byear{1981}).
\btitle{Operator Algebras and Quantum-Statistical Mechanics II},
\bedition{2nd} ed. 
\bpublisher{Springer}, \baddress{New York}.
\bid{mr={611508}}
\end{bbook}
%
\endbibitem

\bibitem{C02}
%
\begin{bbook}[vtex]
\bauthor{\bsnm{Charalambides},~\bfnm{Charalambos~A.}\binits{C.~A.}}
(\byear{2002}).
\btitle{Enumerative Combinatorics}.
\bpublisher{Chapman and Hall}, \baddress{Boca Raton, FL}.
\bid{mr={1937238}}
\end{bbook}
%
\endbibitem

\bibitem{DZ98}
%
\begin{bbook}[mr]
\bauthor{\bsnm{Dembo},~\bfnm{Amir}\binits{A.}} \AND
\bauthor{\bsnm{Zeitouni},~\bfnm{Ofer}\binits{O.}}
(\byear{1998}).
\btitle{Large Deviations Techniques and Applications},
\bedition{2nd} ed.
\bseries{Applications of Mathematics (New York)}
\bvolume{38}.
\bpublisher{Springer}, \baddress{New York}.
\bid{mr={1619036}}
\end{bbook}
%
\endbibitem\vadjust{\goodbreak}

\bibitem{DMP05}
%
\begin{barticle}[mr]
\bauthor{\bsnm{Dorlas},~\bfnm{Teunis~C.}\binits{T.~C.}},
\bauthor{\bsnm{Martin},~\bfnm{Philippe~A.}\binits{P.~A.}} \AND
\bauthor{\bsnm{Pule},~\bfnm{Joseph~V.}\binits{J.~V.}}
(\byear{2005}).
\btitle{Long cycles in a perturbed mean field model of a boson gas}.
\bjournal{J. Stat. Phys.}
\bvolume{121}
\bpages{433--461}.
\bid{doi={10.1007/s10955-005-7582-0}, mr={2185335}}
\end{barticle}
%
\endbibitem

\bibitem{F53}
%
\begin{barticle}[vtex]
\bauthor{\bsnm{Feynman},~\bfnm{R.~P.}\binits{R.~P.}}
(\byear{1953}).
\btitle{Atomic theory of the $\lambda$ transition in Helium}.
\bjournal{Phys. Rev.}
\bvolume{91}
\bpages{1291--1301}.
\end{barticle}
%
\endbibitem

\bibitem{F91}
%
\begin{barticle}[mr]
\bauthor{\bsnm{Fichtner},~\bfnm{Karl-Heinz}\binits{K.-H.}}
(\byear{1991}).
\btitle{On the position distribution of the ideal {B}ose gas}.
\bjournal{Math. Nachr.}
\bvolume{151}
\bpages{59--67}.
\bid{doi={10.1002/mana.19911510105}, mr={1121197}}
\end{barticle}
%
\endbibitem

\bibitem{G88}
%
\begin{bbook}[mr]
\bauthor{\bsnm{Georgii},~\bfnm{Hans-Otto}\binits{H.-O.}}
(\byear{1988}).
\btitle{Gibbs Measures and Phase Transitions}.
\bseries{de Gruyter Studies in Mathematics}
\bvolume{9}.
\bpublisher{de Gruyter}, \baddress{Berlin}.
\bid{mr={956646}}
\end{bbook}
%
\endbibitem

\bibitem{GZ93}
%
\begin{barticle}[mr]
\bauthor{\bsnm{Georgii},~\bfnm{Hans-Otto}\binits{H.-O.}} \AND
\bauthor{\bsnm{Zessin},~\bfnm{Hans}\binits{H.}}
(\byear{1993}).
\btitle{Large deviations and the maximum entropy principle for marked point
random fields}.
\bjournal{Probab. Theory Related Fields}
\bvolume{96}
\bpages{177--204}.
\bid{doi={10.1007/BF01192132}, mr={1227031}}%
\end{barticle}%
%
\endbibitem%

\bibitem{G94}
%
\begin{barticle}[mr]
\bauthor{\bsnm{Georgii},~\bfnm{Hans-Otto}\binits{H.-O.}}
(\byear{1994}).
\btitle{Large deviations and the equivalence of ensembles for {G}ibbsian
particle systems with superstable interaction}.
\bjournal{Probab. Theory Related Fields}
\bvolume{99}
\bpages{171--195}.
\bid{doi={10.1007/BF01199021}, mr={1278881}}
\end{barticle}
%
\endbibitem

\bibitem{G70}
%
\begin{bincollection}[vtex]
\bauthor{\bsnm{Ginibre},~\bfnm{J.}\binits{J.}}
(\byear{1971}).
\btitle{Some applications of
functional integration in statistical mechanics}.
In \bbooktitle{Statistical Mechanics and Quantum Field Theory}
(\beditor{C. de Witt and R. Storaeds}, eds.)
\bpages{327--427}.
\bpublisher{Gordon and Breach},
\baddress{New York}.
\end{bincollection}
%
\endbibitem

\bibitem{LSSY05}
%
\begin{bbook}[mr]
\bauthor{\bsnm{Lieb},~\bfnm{Elliott~H.}\binits{E.~H.}},
\bauthor{\bsnm{Seiringer},~\bfnm{Robert}\binits{R.}},
\bauthor{\bsnm{Solovej},~\bfnm{Jan~Philip}\binits{J.~P.}} \AND
\bauthor{\bsnm{Yngvason},~\bfnm{Jakob}\binits{J.}}
(\byear{2005}).
\btitle{The Mathematics of the {B}ose Gas and Its Condensation}.
\bseries{Oberwolfach Seminars}
\bvolume{34}.
\bpublisher{Birkh\"auser}, \baddress{Basel}.
\bid{mr={2143817}}
\end{bbook}
%
\endbibitem

\bibitem{R09}
%
\begin{bmisc}[vtex]
\bauthor{\bsnm{Rafler},~\bfnm{M.}\binits{M.}}
(\byear{2009}).
\btitle{Gaussian Loop- and
Polya processes: A point process approach. Ph.D. thesis, Univ.
Potsdam}.
\end{bmisc}
%
\endbibitem

\bibitem{Rue69}
%
\begin{bbook}[mr]
\bauthor{\bsnm{Ruelle},~\bfnm{David}\binits{D.}}
(\byear{1969}).
\btitle{Statistical Mechanics: {R}igorous Results}.
\bpublisher{W. A. Benjamin, Inc.}, \baddress{New York}.
\bid{mr={0289084}}
\end{bbook}
%
\endbibitem

\bibitem{R71}
%
\begin{bbook}[vtex]
\bauthor{\bsnm{Robinson},~\bfnm{Derek~W.}\binits{D.~W.}}
(\byear{1971}).
\btitle{The Thermodynamic Pressure in Quantum Statistical Mechanics}.
\bseries{Lecture Notes in Physics}
\bvolume{9}.
\bpublisher{Springer}, \baddress{Berlin}.
\bid{mr={0432122}}
\end{bbook}
%
\endbibitem

\bibitem{S93}
%
\begin{barticle}[mr]
\bauthor{\bsnm{S{\"u}t{\H{o}}},~\bfnm{Andr{\'a}s}\binits{A.}}
(\byear{1993}).
\btitle{Percolation transition in the {B}ose gas}.
\bjournal{J. Phys. A}
\bvolume{26}
\bpages{4689--4710}.
\bid{mr={1241339}}%
\end{barticle}%
%
\endbibitem%

\bibitem{S02}
%
\begin{barticle}[mr]
\bauthor{\bsnm{S{\"u}t{\H{o}}},~\bfnm{Andr{\'a}s}\binits{A.}}
(\byear{2002}).
\btitle{Percolation transition in the {B}ose gas. {II}}.
\bjournal{J. Phys. A}
\bvolume{35}
\bpages{6995--7002}.
\bid{doi={10.1088/0305-4470/35/33/303}, mr={1945163}}
\end{barticle}
%
\endbibitem

\bibitem{Toth90}
%
\begin{barticle}[vtex]
\bauthor{\bsnm{T\'{o}th},~\bfnm{B\'{a}lint}\binits{B.}}
(\byear{1990}).
\btitle{Phase transition in an interacting {B}ose system. {A}n
application of
the theory of {V}entsel' and {F}reidlin}.
\bjournal{J. Stat. Phys.}
\bvolume{61}
\bpages{749--764}.
\bid{mr={1086297}}
\end{barticle}
%
\endbibitem

\bibitem{Ver96}
%
\begin{barticle}[mr]
\bauthor{\bsnm{Vershik},~\bfnm{A.~M.}\binits{A.~M.}}
(\byear{1996}).
\btitle{Statistical mechanics of combinatorial partitions, and their limit
configurations}.
\bjournal{Funktsional. Anal. i Prilozhen.}
\bvolume{30}
\bpages{19--39, 96}.
\bid{doi={10.1007/BF02509449}, mr={1402079}}
\end{barticle}
%
\endbibitem

\end{thebibliography}
\end{document}